\input ./style/arxiv-general.cfg
\documentclass[MSNbibl,number,citesort,seceqn,dvips]{arxbj}
\makeatletter
   \@ifpackageloaded{graphicx}{}{\usepackage{graphicx}}
\makeatother
\usepackage{mathbh}
\volume{22}
\issue{4}
\pubyear{2016}
\firstpage{2237}
\lastpage{2300}
\doi{10.3150/15-BEJ728}% Updated by VTEXPTS2LaTeX.exe, 23.06.2015 13:40
\docsubty{FLA}

\makeatletter
\DeclareFontFamily{U}{txsyc}{}
\DeclareFontShape{U}{txsyc}{m}{n}{
   <-> txsyc%
}{}
\DeclareFontShape{U}{txsyc}{bx}{n}{
   <-> txbsyc%
}{}
\DeclareFontShape{U}{txsyc}{l}{n}{<->ssub * txsyc/m/n}{}
\DeclareFontShape{U}{txsyc}{b}{n}{<->ssub * txsyc/bx/n}{}
\DeclareSymbolFont{symbolsC}{U}{txsyc}{m}{n}
\SetSymbolFont{symbolsC}{bold}{U}{txsyc}{bx}{n}
\DeclareFontSubstitution{U}{txsyc}{m}{n}
\newcommand{\Na}{N^{(\alpha)}}
\newcommand{\righttoleftarrow}{\circlearrowleft}
\DeclareFontFamily{U}{mathx}{\hyphenchar\font45}
\DeclareFontShape{U}{mathx}{m}{n}{
      <5> <6> <7> <8> <9> <10>
      <10.95> <12> <14.4> <17.28> <20.74> <24.88>
      mathx10
      }{}
\DeclareSymbolFont{mathx}{U}{mathx}{m}{n}
\DeclareFontSubstitution{U}{mathx}{m}{n}
\DeclareMathAccent{\widecheck}{0}{mathx}{"71}
\def\sfrac#1#2{#1/#2}
\def\vfrac#1#2{(#1)/#2}
\def\afrac#1#2{#1/(#2)}
\def\vafrac#1#2{(#1)/(#2)}

\newcommand{\rrvert}{\vert}

\newcommand{\rrVert}{\Vert}
\newcommand{\llvert}{\vert}
\newcommand{\llVert}{\Vert}
\newcommand{\eqref}[1]{(\ref{#1})}
\newcommand{\ffff}{:=} %\mathrel}{symbolsC}{"42}
\newcommand{\fd}{=:} %\mathrel}{symbolsC}{"43}
\newcommand{\cD}{\mathcal{D}}
\newcommand{\cI}{\mathcal{I}}
\newcommand{\cL}{\mathcal{L}}
\newcommand{\cN}{\mathcal{N}}
\newcommand{\cM}{\mathcal{M}}
\newcommand{\cO}{\mathcal{O}}
\newcommand{\cT}{\mathcal{T}}
\newcommand{\CC}{\mathbb{C}}
\newcommand{\LL}{\mathbb{L}}
\newcommand{\NN}{\mathbb{N}}
\newcommand{\PP}{\mathbb{P}}
\newcommand{\RR}{\mathbb{R}}
\newcommand{\TT}{\mathbb{T}}
\newcommand{\ZZ}{\mathbb{Z}}
\newcommand{\tr}{\triangle}
\newcommand{\osc}{\operatorname{osc}}
\newcommand{\sign}{\operatorname{sign}}
\newcommand{\st}{:}
\newcommand{\un}{\mathbh{1}}
\newcommand{\Var}{\operatorname{Var}}
\newcommand{\lin}{[\![}
\newcommand{\rin}{]\!]}
\newtheorem{pro}{Proposition}[section]
\newtheorem{lem}{Lemma}[section]
\newtheorem{theo}{Theorem}[section]
\newremark{con}{Conjecture}[section]
\newremark{rem}{Remark}[section]
\makeatother

\begin{document}
\begin{frontmatter}

\title{A stochastic algorithm finding $p$-means on~the circle}
\runtitle{A stochastic algorithm finding $p$-means on the circle}

\begin{aug}
%%%% inicialai - be tarpu
% Corresponding author: Laurent Miclo - miclo@math.univ-toulouse.fr% Updated by VTEXPTS2LaTeX.exe, 25.06.2015 08:48
%Updated by VTEXPTS2LaTeX.exe, 23.06.2015 13:40
\author[A]{\inits{M.}\fnms{Marc}~\snm{Arnaudon}\thanksref{A}\ead[label=e1]{marc.arnaudon@u-bordeaux1.fr}} \and
\author[B]{\inits{L.}\fnms{Laurent}~\snm{Miclo}\corref{}\thanksref{B}\ead[label=e2]{miclo@math.univ-toulouse.fr}}
%%\runauthor{} %% auto
%\dedicated{}
\address[A]{Institut de Math\'{e}matique de Bordeaux, UMR 5251,
Universit\'{e} de Bordeaux and CNRS,
351, Cours de la Lib\'{e}ration,
F-33405 TALENCE Cedex, France. \printead{e1}}
\address[B]{Institut de Math\'{e}matiques de Toulouse, UMR 5219,
Universit\'e Toulouse 3 and CNRS,
118, route de Narbonne,
31062 Toulouse Cedex 9, France. \printead{e2}}
\end{aug}

% HISTORY:
%
\received{\smonth{9} \syear{2013}}% Updated by VTEXPTS2LaTeX.exe,
%23.06.2015 13:40
%
\revised{\smonth{3} \syear{2015}}% Updated by VTEXPTS2LaTeX.exe,
%23.06.2015 13:40

% ABSTRACT
%
\begin{abstract}
A stochastic algorithm is proposed, finding some elements from the set
of intrinsic $p$-mean(s)
associated to a probability measure $\nu$ on a compact Riemannian
manifold and to $p\in[1,\infty)$.
It is fed sequentially with independent random variables $(Y_n)_{n\in
\mathbb{N}}$ distributed
according to $\nu$, which is often the only available knowledge of
$\nu$.
Furthermore, the algorithm is easy to implement, because it evolves
like a Brownian
motion between the random times when it jumps
in direction of one of the $Y_n$, $n\in\mathbb{N}$. Its principle is
based on simulated
annealing and homogenization, so that temperature and approximations
schemes must
be tuned up (plus a regularizing scheme if $\nu$ does not admit a H\"
olderian density).
The analysis of the convergence is restricted to the case where the
state space is a
circle. In its principle, the proof relies on the investigation of the
evolution of
a time-inhomogeneous $\mathbb{L}^2$ functional and on the
corresponding spectral gap estimates due to Holley, Kusuoka and Stroock.
But it requires new estimates on the discrepancies between the unknown
instantaneous invariant measures and some convenient Gibbs measures.
\end{abstract}

% KEYWORDS
% visi is mazosios raides ir pagal abecele
%
\end{frontmatter}

%s1 #&#
\section{Introduction}\label{sec1}

The purpose of this paper is to present a stochastic algorithm
finding some of the geometric \mbox{$p$-}means of probability measures defined
on compact Riemannian manifolds, for $p\in[1,\infty)$.
Its convergence is analyzed in the restricted case of the circle, as a
first step toward a more general result which is conjectured to be true.

%s1.1 #&#
\subsection{The general notion of $p$-means}\label{sec1.1}

The concepts of mean and median are well understood for real valued
random variables.
They can be extended to random variables taking values in metric spaces
in the following way.
Let
be given $\nu$ a probability measure on a metric space $M$, whose
distance is denoted $d$.
For $p\geq1$, consider the continuous mapping
%
%
%e1.1 #&#
\begin{eqnarray}
\label{U} U_p\st M\ni x\mapsto\int d^p(x,y) \nu(dy).
\end{eqnarray}
A global minimum of $U_p$ is called a $p$-mean of $\nu$,
at least if this function is not identically equal to $+\infty$
(equivalently, if all its values are finite, as it can be easily
deduced from the triangle inequality).
The set of $p$-means will be designated by
$\cM_p$, it is non-empty
as soon as $U_p$ goes to infinity at infinity (in the Alexandroff sense),
but in general it is not reduced to a singleton.
The notion of intrinsic mean and median correspond, respectively, to
$p=2$ and $p=1$.
If $M$ is $\RR$ endowed with its absolute value, one recovers the
usual mean and distance.

These extensions are justified by the increasing number of available
graph or manifold valued data samples in various scientific fields.
Examples of manifold valued data samples are given by sets of
parameters for families of laws endowed with Fisher information metric,
by Lie groups
(rotations, displacements) in control theory, by symmetric spaces in
imaging or signal processing.

For some applications (see, e.g., \cite{MR2254442}), it may be
important to find $\cM_p$ or at least some of its elements. In
practice, the knowledge
of $\nu$ is often given by a finite sequence $Y\ffff (Y_n)_{n\in\{1, 2,
\ldots, N\}}$ of independent random variables, identically distributed
according to $\nu$. Since $N\in\NN$ is in general large enough, we
will consider the limit situation where we have at our disposal
an infinite sequence $Y\ffff (Y_n)_{n\in\NN}$. One is then looking for
algorithms using this data and enabling to find some
elements of $\cM_p$.
In this paper, we will be mainly interested in the case where $M$ is
the circle, even if the proposed stochastic algorithm
can be considered more generally for compact Riemannian manifolds.

Algorithms for finding $p$-means or minimax centers have been
investigated in \cite{Le04,Sturm05,Groisser05,Groisser06,Badoiu-Clarkson03,Yang10,Bonnabel11,Afsari-Tron-Vidal11,Arnaudon-al12,Cardot-Cenac-Zitt12,Arnaudon-Nielsen12a}.
When possible a gradient descent algorithm is used. When the gradient
of the functional to minimize is difficult or impossible to compute,
a Robbins Monro-type algorithm is preferred. Either the functional to
minimize has only one local minimum which is also global,
or (Bonnabel \cite{Bonnabel11})
a local minimum is seeked. The case of Karcher means in the circle is
treated in \cite{2011arXiv11091986C} and \cite{2011arXiv11082141H}.
In this special situation, the global minimum of the functional can be
found by explicit formula.

For generalized means on compact manifolds, the situation is different
since the functional~\eqref{U} to minimize may have many local
minima, and no explicit formula for a global minimum can be expected.

%s1.2 #&#
\subsection{The case of the circle}\label{sec1.2}\label{tcotc}

In this subsection, we consider the case where $M$ is the circle $\TT
\ffff \RR/(2\pi\ZZ)$ endowed with its natural angular distance $d$.
As above, let $Y\ffff (Y_n)_{n\in\NN}$ be a sequence of independent
random variables distributed according to a fixed probability measure
$\nu$ on $\TT$.
Let $p\in[1,+\infty)$ be fixed, we present now a stochastic algorithm
finding some
elements of $\cM_p$ by using this data. It is based on simulated
annealing and homogenization procedures.
Thus, we will need, respectively, an inverse temperature evolution
$\beta
\st\RR_+\rightarrow\RR_+$ and an inverse speed up evolution
$\alpha\st\RR_+\rightarrow\RR_+^*$, where $\RR_+^*$ stands for
the set of
positive real numbers. Typically, they are, respectively, non-decreasing
and non-increasing and we have $\lim_{t\rightarrow+\infty}\beta
_t=+\infty$ and
$\lim_{t\rightarrow+\infty}\alpha_t=0$, but we are looking for more precise
conditions so that the
stochastic algorithm we describe below finds $\cM_p$ (namely, some
elements from this set).

Let $N\ffff (N_t)_{t\geq0}$ be a standard Poisson process: it starts at~0 at time 0 and has jumps of length~1 whose inter-arrival times are
independent and distributed according to exponential random variables
of parameter 1.
The process $N$ is assumed to be independent from the sequence~$Y$.
We define the speeded-up process $N^{(\alpha)}\ffff (N^{(\alpha
)}_t)_{t\geq0}$ via
%
%
%e1.2 #&#
\begin{eqnarray}
\label{Na} \forall t\geq0,\qquad N^{(\alpha)}_t&\ffff &
N_{\int_0^t \sfrac{1}{\alpha_s}
\,ds}.
\end{eqnarray}
Consider the time-inhomogeneous Markov process $X\ffff (X_t)_{t\geq0}$
which evolves in $M$ in the following heuristic way:
if $T>0$ is a jump time of $N^{(\alpha)}$, then $X$ jumps at the same time,
from $X_{T-}$
to $X_T$ which is obtained by following the shortest geodesic leading
from $X_{T-}$
to
$Y_{N^{(\alpha)}_{T}}$ at speed~1 during the time $(p/2)\beta_T\alpha
_Td^{p-1}(X_{T_-},Y_{N^{(\alpha)}_{T}})$.
Almost surely, the above shortest geodesic is unique
and there is no problem with its choice.
Indeed, by the end of the description below, $X_{T_-}$
will be independent of $Y_{N^{(\alpha)}_{T}}$ and
the law of $X_{T_-}$ will be absolutely continuous
with respect to the~Lebesgue measure $\lambda$ on $\TT$ renormalized
into a probability measure.
It ensures that almost surely, $Y_{N^{(\alpha)}_{T}}$ is not the
opposite point of
$X_{T_-}$ on $\TT$. The\vspace*{1pt} schemes $\alpha$ and $\beta$ will satisfy
$\lim_{t\rightarrow+\infty}\alpha_t\beta_t=0$,
so that for sufficiently large jump-times $T$, $X_T$ will be between
$X_{T-}$ and $Y_{N^{(\alpha)}_{T}}$
on the above geodesic and quite close to $X_{T-}$.\vspace*{1pt}

To proceed with the construction, we require that between consecutive
jump times (and between time 0 and
the first jump time), $X$ evolves as a Brownian motion on $\TT$ and
independently
of $Y$ and $N$.
Very informally, the evolution of the algorithm $X$ can be summarized
by the equation
\begin{eqnarray*}
\forall t\geq0,\qquad dX_t&=&d B_t+(p/2)
\alpha_t\beta_t d^{p-1}(X_{T_-},Y_{N^{(\alpha
)}_{T}})
\sigma(X_{t-},Y_{N^{(\alpha)}_t})\, dN^{(\alpha)}_t,
\end{eqnarray*}
where\vspace*{1pt} $(B_t)_{t\geq0}$ is a Brownian motion on $\TT$ and where
$\sigma(X_{t-},Y_{N^{(\alpha)}_t})$
is $1$ (resp., $-1$) if the shortest way from $X_{t-}$ to $Y_{N^{(\alpha)}_t}$
goes in the anti-clock wise
(resp., the clock-wise) direction,\vspace*{1pt} in the usual~representation of $\RR
/(2\pi\ZZ)$ in $\CC$.
In the above equation,
$(Y_{N^{(\alpha)}_t})_{t\geq0}$ should be interpreted as a fast
auxiliary process.
The law of $X$ is then entirely determined by the initial distribution
$m_0=\cL(X_0)$.
More generally at any time $t\geq0$, denote by $m_t$ the law of
$X_t$.

The first main result of this paper states that
at least if $\nu$ is sufficiently regular,
the above algorithm $X$ finds in probability at large times
the set $\cM_p$ of $p$-means:
%

%th1.1 #&#
\begin{theo}\label{t1}
Assume that $\nu$ admits a density with respect to $\lambda$ and that
this density is H\"older continuous with exponent $a\in(0,1]$. Then
there exist two constants $a_p>0$, depending on $p\geq1$ and $a$, and
$b_p\geq0$, depending on $p$,
such that for any scheme of the form
%
%
%e1.3 #&#
\begin{eqnarray}
\label{ab} \forall t\geq0,\qquad\cases{ \alpha_t \ffff
(1+t)^{-\sfrac{1}{a_p}},
\vspace*{3pt}\cr
\beta_t \ffff  b^{-1}\ln(1+t),}
\end{eqnarray}
where $b>b_p$, we have for any neighborhood $\cN$ of $\cM_p$ and for
any $m_0$,
%
%
%e1.4 #&#
\begin{eqnarray}
\label{mN1} \lim_{t\rightarrow+\infty}\PP[X_t\in\cN]&=&1.
\end{eqnarray}
\end{theo}

Thus, to find an element of $\cM_p$ with an important probability, one should
pick up the value of $X_t$ for sufficiently large times $t$.

The constant $a_p$ is the simplest to define, since it is given by
%
%
%e1.5 #&#
\begin{eqnarray}
\label{ap} a(p)&\ffff & \cases{ a, &\quad if $p=1$ or $p\geq2$,
\vspace*{3pt}\cr
\min(a,p-1), &
\quad if $p\in(1,2)$.}
\end{eqnarray}
The constant $b_p\geq0$ comes from the theory of simulated annealing
(see, e.g.,
\cite{MR995752}), which will be recalled in next section.
For the moment being, we just describe the constant $b_p$, in the
setting of a compact Riemannian manifold $M$,
since there is no extra difficulty and we will need it later on to
express a conjecture extending Theorem~\ref{t1}.
For any $x,y\in M$, let $\mathcal{C}_{x,y}$ be the set of continuous
paths $C\ffff (C(t))_{0\leq t\leq1}$
going from $C(0)=x$ to $C(1)=y$.
The elevation $U_p(C)$ of such a path $C$ relatively to $U_p$ is
defined by
\begin{eqnarray*}
U_p(C)&\ffff & \max_{t\in[0,1]} U_p\bigl(C(t)
\bigr)
\end{eqnarray*}
and the minimal elevation $U_p(x,y)$ between $x$ and $y$ is
given by
\begin{eqnarray*}
U_p(x,y)&\ffff & \min_{C\in\mathcal{C}_{x,y}} U_p(C).
\end{eqnarray*}
Then we consider
%
%
%e1.6 #&#
\begin{eqnarray}
\label{bU} b(U_p)&\ffff & \max_{x,y\in M}U_p(x,y)-U_p(x)-U_p(y)+
\min_M U_p.
\end{eqnarray}
This constant can also be seen as the largest depth of a well
not containing
a fixed global minimum of $U_p$. Namely, if $x_0\in\cM_p$,
then it is not difficult to see that
%
%
%e1.7 #&#
\begin{eqnarray}
\label{bU2} b(U_p)&=& \max_{y\in M}
U_p(x_0,y)-U_p(y),
\end{eqnarray}
independently of the choice of
$x_0\in\cM_p$ (cf. \cite{MR995752}).

Let us now describe a stochastic algorithm, derived from the previous
one, which enables one to find some of the $p$-means of any probability
measure $\nu$
on $\TT$.

For any $x\in\TT$ and $\kappa>0$, consider the probability measure $
K_{x,\kappa}$ whose density with respect to the Lebesgue measure
$\lambda(dy)$
is proportional to $(1-\kappa\llVert y-x\rrVert )_+$.
Assume next that we are given an evolution $\kappa\st\RR_+\ni
t\mapsto\kappa_t\in\RR_+^*$
and consider the process $Z\ffff (Z_t)_{t\geq0}$ evolving similarly to
$(X_t)_{t\geq0}$, except that
at the jump times $T$ of $N^{(\alpha)}$, the target $Y_{N^{(\alpha
)}_T}$ is replaced by
a point $W_T$
sampled from $K_{Y_{N^{(\alpha)}_T}, \kappa_T}$, independently from
the other
variables.

%
%
%th1.2 #&#
\begin{theo}\label{t2}
Let $\nu$ be an arbitrary probability measure
on $M=\TT$. For $p=2$, consider the schemes
\begin{eqnarray*}
\forall t\geq0,\qquad\cases{ \alpha_t \ffff (1+t)^{-c},
\vspace*{3pt}\cr
\beta_t \ffff  b^{-1}\ln(1+t),
\vspace*{3pt}\cr
\kappa_t \ffff
(1+t)^k,}
\end{eqnarray*}
with $b>b(U_2)$, $k>0$ and $c\geq2k+1$.
Then, for any neighborhood $\cN$ of $\cM_2$
and for any initial distribution
$\cL(Z_0)$, we get
\begin{eqnarray*}
\lim_{t\rightarrow+\infty}\PP[Z_t\in\cN]&=&1,
\end{eqnarray*}
where $\PP$ stands for the underlying probability.
\end{theo}

More generally, for any given $p\geq1$, it is possible to find similar schemes
(where $c$ depends furthermore on $p\geq1$)
enabling to find the set of $p$-means $\cM_p$ (see Remark~\ref{t2p1}).
Even if $\nu$ satisfies the condition of Theorem~\ref{t1}, it could
be more advantageous to
consider the alternative algorithm $Z$ instead of $X$ when the exponent
$a$ in (\ref{ab}) is too small.
%

%re1.1 #&#
\begin{rem}
The schemes $\alpha$, $\beta$ and $\kappa$ presented above are
simple examples of admissible evolutions;
they could be replaced, for instance, by
\begin{eqnarray*}
\forall t\geq0,\qquad\cases{ \alpha_t \ffff  C_1(r_1+t)^{-c},
\vspace*{3pt}\cr
\beta_t \ffff  b^{-1}\ln(r_2+t),
\vspace*{3pt}\cr
\kappa_t =C_2(r_3+t)^{k},}
\end{eqnarray*}
where $C_1, C_2>0$, $r_1,r_3>0$, $r_2\geq1$ and still under the
conditions $b>b(U_p)$, $k>0$ and $c\geq2k+1$. It is possible to deduce
more general conditions insuring the validity of the convergence
results of
Theorems~\ref{t1} and~\ref{t2} (see, e.g., Proposition~\ref{condab} below).
\end{rem}

How to choose in practice the exponents $c$ and $k$ satisfying $c\geq
2k+1$ in Theorem~\ref{t2}? We note that the larger $c$, the faster
$\alpha$ goes to zero and the faster the algorithm $Z$
is using the data $(Y_n)_{n\in\NN}$. In compensation, $k$ can be
chosen larger, which means
that $\nu$ is closer to its approximation by its transport through the
kernel $K_{\cdot,\kappa_t}(\cdot) $ (defined before the statement of
Theorem~\ref{t2}, for more details see Section~\ref{sec5}), namely the
convergence will be more precise. This is quite natural,
since more data have been required at some fixed time. So in practice a
trade-off has to be made
between the number of i.i.d. variables distributed according to $\nu$
one has at his disposal
and the quality of the approximation of $\cM_p$.

%s1.3 #&#
\subsection{Numerical illustration}\label{sec1.3}\label{ni}

The algorithm $X$ (and similarly for $Z$) is not so difficult to implement.
Let us identify $\TT$ with $(-\pi,\pi]$ and construct $X_t$ for some
fixed $t>0$.
Assume we are given $(Y_n)_{n\in\NN}$, $(\alpha_s)_{s\in[0,t]}$,
$(\beta_s)_{s\in[0,t]}$
and $X_0$ as in the \hyperref[sec1]{Introduction}. We need furthermore two independent
sequences $(\tau_n)_{n\in\NN}$ and $(V_n)_{n\in\NN}$, consisting
of i.i.d. random variables, respectively, distributed according to the
exponential law of parameter 1 and to the Gaussian law with mean 0 and
variance 1.
We begin by constructing the finite sequence $(T_n)_{n\in\lin0,N\rin
}$ corresponding to the jump times of $N^{(\alpha)}$: let $T_0\ffff 0$
and next by iteration, if $T_n$ was defined,
we take $T_{n+1}$ such that $\int_{T_n}^{T_{n+1}}1/\alpha_s \,ds=\tau_{n+1}$.
This is done until $T_N> t$, with $N\in\NN$, then we change the
definition of $T_N$ by imposing
$T_N=t$.
Next, we consider the sequence $(\widecheck{X}_n,\widehat X_n)_{n\in
\lin
0,N\rin}$ constructed through the following iteration (where the
variables are reduced modulo $2\pi$): starting from $\widecheck
{X}_0\ffff \widehat X_0\ffff  X_0$,
if $\widehat X_n$ was defined, with $n\in\lin0, N-1\rin$, we consider
%
%
%e1.8 #&#
\begin{eqnarray}
\label{widecheckX} \widecheck{X}_{n+1}\ffff \widehat X_n+
\sqrt{T_{n+1}-T_n} V_{n+1}.
\end{eqnarray}
Next, we define
%
%
%e1.9 #&#
\begin{eqnarray}
\label{witX} \widehat X_{n+1}&\ffff & \widecheck{X}_{n+1}+(p/2)
\alpha_{T_{n+1}}\beta_{T_{n+1}} \llvert W_{n+1}\rrvert
^{p-2}V_{n+1},
\end{eqnarray}
where
$W_{n+1}$ is the representative of $Y_{n+1}-\widecheck{X}_{n+1}$ in
$(-\pi,\pi]$ modulo $2\pi$.
Then $\widecheck{X}_{N}$ has the same law as $X_t$.

Theorems~\ref{t1} and~\ref{t2}
provide theoretical results at very large times,
but in practice, one has to work with a finite horizon $t$, for which
the best corresponding scheme $\beta$
may not be of the form of those given in these theorems
(see the lectures of \cite{MR2001g60167} for the classical simulated
annealing algorithm).
Thus, the previous theorems should only be seen as indications of what
could be tried in practice.
Let us illustrate that by some numerical simulations.
On the circle, still identified with $(-\pi,\pi]$, consider the
probability distribution $\nu=(\delta_0+\delta_\pi)/2$.
A priori we should resort to Theorem~\ref{t2}, but let us just
``apply'' Theorem~\ref{t1}
with $a=1$, namely with the scheme
\begin{eqnarray*}
\forall t\geq0,\qquad\alpha_t&\ffff & \frac{1}{1+t}.
\end{eqnarray*}
For $p=1$ the function $U_1$ is constant, meaning that the set of
medians $\cM_1$
is the whole circle. For $p>1$, the function $U_p$ admits two global
minima, $\cM_p=\{-\pi/2,\pi/2\}$,
and two global maxima, $0$ and $\pi$.
It is easy to see that $b(U_p)=\pi^p(1-2^{1-p})$, so that we can take,
for instance,
\begin{eqnarray*}
\forall t\geq0,\qquad\beta_t&\ffff & \frac{2}{\pi^p(1-2^{1-p})}\ln(1+t),\vadjust{\goodbreak}
\end{eqnarray*}
(for $p=1$, the factor in front of the logarithm can be chosen freely,
one could even choose
the scheme $\beta$ to be constant).
%, taken to be equal to 1 in the simulation below.
With the above notation, let $(Y_n)_{n\in\NN}$, $(\tau_n)_{n\in\NN
}$ and $(V_n)_{n\in\NN}$ be independent sequences consisting of
i.i.d. random variables, respectively, distributed according to the
uniform law
on $\{0,\pi\}$, to the
exponential law of parameter 1 and to the Gaussian law with mean 0 and
variance 1.
Let $t>0$ be fixed.
The finite sequence $(T_n)_{n\in\lin0,N\rin}$ is constructed through
the recurrence $T_0=100$
and
\begin{eqnarray*}
\forall n\in\lin0, N-1\rin,\qquad T_{n+1}&\ffff & \sqrt{(T_n+1)^2+
\tau_{n+1}}-1
\end{eqnarray*}
%

%
%f1 #&#
\begin{figure}[t]
\begin{tabular}{cc}

\includegraphics{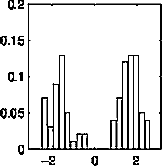}  & \includegraphics{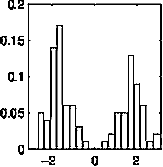}\\
\footnotesize{(a)} & \footnotesize{(b)}\\[6pt]

\includegraphics{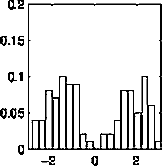}  & \includegraphics{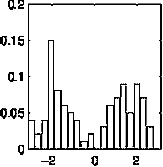}\\
\footnotesize{(c)} & \footnotesize{(d)}
\end{tabular}
\caption{\textup{(a)}~$p=2$ and $t=200$,
\textup{(b)}~$p=2$ and $t=400$,
\textup{(c)}~$p=1.1$ and $t=200$,
\textup{(d)}~$p=1.1$ and $t=400$.}\label{fig1}
\end{figure}

\noindent 
until $T_N>t$.
% Next we change the definition by imposing $T_N\ffff  t$.
Starting from $\widecheck{X}_0\ffff \widehat X_0\ffff 0$, we consider the
sequence $(\widecheck{X}_n,\widehat X_n)_{n\in\lin0,N\rin}$ defined via
\eqref{widecheckX} and \eqref{witX}.
The histograms of Figure~\ref{fig1} of the distribution of $\widecheck{X}_N$
correspond to $p=1.1$ and $p=2$ and $t=200$ and $t=400$ and they
are obtained with 100 samples of the procedure described above.

It appears that as time goes on, there is a tendency to concentrate on
the set of means $\{-\pi/2,\pi/2\}$,
but that this is more difficult to achieve for small $p>1$, due to the
fact that in the limit case $p=1$, one is
trying to sample according to the uniform distribution on $(-\pi,\pi
]$.

%
%f2 #&#
\begin{figure}[t]

\includegraphics{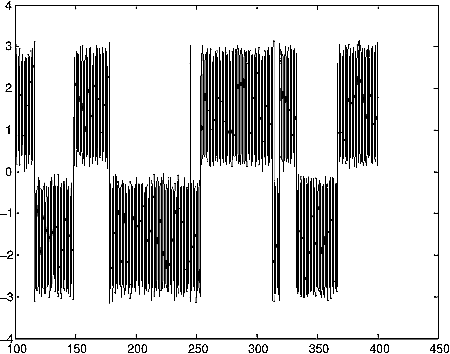}

\caption{A trajectory for $p=2$ and $t=400$.}\label{fig2}
\end{figure}

Figure~\ref{fig2} is plotting a typical trajectory (observed at the jump
times), with $p=2$, $t=400$ and for which the simulation gave $N={}$150,366
(close to $400^2-100^2$).
% (note that the time scale is not the natural one).
It should be emphasized that if instead of using 100 samples in a
Monte-Carlo procedure
as above, one rather resorts to the empirical measure generated by one
trajectory, one would get similar histograms.

%s1.4 #&#
\subsection{The conjecture for Riemannian manifolds}\label{sec1.4}

The description of the algorithm given in Section~\ref{tcotc} can be
extended to any compact Riemannian manifold $M$
endowed with its distance $d$. For general books on Riemannian
geometry, we refer to~\cite{Jost11}.

As above, let $Y\ffff (Y_n)_{n\in\NN}$ be a sequence of independent
random variables distributed according to a fixed probability measure
$\nu$ on $M$.
Let $p\in[1,+\infty)$ be fixed. We also
need an inverse temperature evolution $\beta\st\RR_+\rightarrow\RR
_+$ and
an inverse speed up evolution
$\alpha\st\RR_+\rightarrow\RR_+^*$, which typically will be non-decreasing
and non-increasing and satisfying $\lim_{t\rightarrow+\infty}\beta
_t=+\infty$ and
$\lim_{t\rightarrow+\infty}\alpha_t=0$.

We consider again the speeded-up process $N^{(\alpha)}\ffff (N^{(\alpha
)}_t)_{t\geq0}$ via
\begin{eqnarray*}
\forall t\geq0,\qquad N^{(\alpha)}_t&\ffff & N_{\int_0^t \sfrac{1}{\alpha_s}
\,ds},
\end{eqnarray*}
where $N\ffff (N_t)_{t\geq0}$ be a standard Poisson process independent
from $Y$.
The time-inhomoge\-neous Markov process $X\ffff (X_t)_{t\geq0}$ evolves in
$M$ in the following heuristic way:
if $T>0$ is a jump time of $N^{(\alpha)}$, then $X$ jumps at the same time,
from $X_{T-}$
to
\begin{eqnarray*}
X_T\ffff \exp_{X_{T-}}\bigl((p/2)\beta_T\alpha
_Td^{p-2}(X_{T_-},Y_{N^{(\alpha)}
_{T}})
\overrightarrow{X_{T_-}Y_{N^{(\alpha)}_{T}}}\bigr).
\end{eqnarray*}
By definition, the latter point is obtained by the
following  the shortest geodesic leading from $X_{T-}$ to
$Y_{\Na_{T}}$ at time 1, during a time
$s:=(p/2)\beta_T\alpha_Td^{p-2}(X_{T_-},Y_{\Na_{T}})$ (and\vspace*{1pt} thus may not really correspond
to an image of the exponential mapping if $s$ is not small enough).
The schemes $\alpha$ and $\beta$ will satisfy $\lim_{t\rightarrow
+\infty
}\alpha_t\beta_t=0$,
so that for sufficiently large jump-times $T$, $X_T$ will be between
$X_{T-}$ and $Y_{N^{(\alpha)}_{T}}$
on the above geodesic and quite close to $X_{T-}$.
Almost surely, the above shortest geodesics are unique
and there is no problem with their choices in the previous construction.
Indeed, by the end of the description below, $X_{T_-}$
will be independent of $Y_{N^{(\alpha)}_{T}}$ and
the law of $X_{T_-}$ will be absolutely continuous
with respect to the Riemannian probability $\lambda$, namely the
volume measure standardized to total volume one. It ensures that almost
surely, $Y_{N^{(\alpha)}_{T}}$ is not in the cut-locus of
$X_{T_-}$ (which is negligible with respect to $\lambda$) so that
there is only one shortest geodesic from $X_{T-}$
to
$Y_{N^{(\alpha)}_{T}}$.
To proceed with the construction, we require that between consecutive
jump times (and between time 0 and
the first jump time), $X$ evolves as a Brownian motion, relatively to
the Riemannian structure of
$M$ (see, e.g., the book of \cite{MR1011252}) and independently
of $Y$ and $N$.
Very informally, the evolution of the algorithm $X$ can be summarized
by the equation (in the tangent bundle $TM$)
\begin{eqnarray*}
\forall t\geq0,\qquad dX_t&=&d B_t+(p/2)
\alpha_t\beta_t d^{p-2}(X_{T_-},Y_{N^{(\alpha)}
_{T}})
\overrightarrow{X_{t-}Y_{N^{(\alpha)}_t}}\, dN^{(\alpha)}_t,
\end{eqnarray*}
where $(B_t)_{t\geq0}$ would be a Brownian motion on $M$ and where
$(Y_{N^{(\alpha)}_t})_{t\geq0}$ should be interpreted as a fast
auxiliary process.
The law of $X$ is then entirely determined by the initial distribution
$m_0=\cL(X_0)$.
We believe that the above algorithm $X$ finds in probability at large times
the set $\cM_p$ of $p$-means, at least if $\nu$ is sufficiently
regular, as in the case where $M=\TT$:

%
%co1.1 #&#
\begin{con}\label{conj}
Assume that $\nu$ admits a density with respect to $\lambda$ and that
this density is H\"older continuous with exponent $a\in(0,1]$. Then
there exist two constants $a_p>0$, depending on $p\geq1$ and $a$, and
$b_p\geq0$, depending on $p$ and $M$,
such that for any scheme of the form given in \eqref{ab},
where $b>b_p$, we have for any neighborhood $\cN$ of $\cM_p$ and for
any $m_0$,
\begin{eqnarray*}
\lim_{t\rightarrow+\infty}\PP[X_t\in\cN]&=&1.
\end{eqnarray*}
\end{con}

So as in Section~\ref{tcotc}, to find an element of $\cM_p$ with an
important probability, one should
pick up the value of $X_t$ for sufficiently large times $t$.

The constant $b_p\geq0$ should still coincide with the one defined in
\eqref{bU2}.

Let us now extend the stochastic algorithm $Z$, which should enable one
to find some of the $p$-means of any probability measure $\nu$
on the compact Riemannian manifold $M$.

For any $x\in M$ and $\kappa>0$, consider, on the tangent space $T_xM$,
the probability measure $\widetilde K_{x,\kappa}$ whose density with respect
to the Lebesgue measure $dv$
is proportional to $(1-\kappa\llVert v\rrVert )_+$ (where the Lebesgue
measure and the norm
are relative to the Euclidean structure on $T_xM$).
Denote $K_{x,\kappa}$ the image by the exponential mapping at $x$ of
$\widetilde K_{x,\kappa}$.
Assume next that we are given an evolution $\kappa\st\RR_+\ni
t\mapsto\kappa_t\in\RR_+^*$
and consider the process $Z\ffff (Z_t)_{t\geq0}$ evolving similarly to
$(X_t)_{t\geq0}$, except that
at the jump times $T$ of $N^{(\alpha)}$, the target $Y_{N^{(\alpha
)}_T}$ is replaced by
a point $W_T$
sampled from $K_{Y_{N^{(\alpha)}_T}, \kappa_T}$, independently from
the other
variables.
We believe that a variant of Theorem~\ref{t2} should hold more
generally on compact Riemannian manifolds $M$. But it seems that the
geometry of $M$ should play a role, especially through the behaviour of
the volume of small enlargements of the cut-locus of points.

Notice that a major difficulty for implementing an algorithm in a
high-dimensional manifold simulating the process $X_t$ is to compute the
logarithm map
$\overrightarrow{xy}=\exp_x^{-1}(y)$. Moreover, this logarithm can be
very instable around the cutlocus of~$x$. In~\cite
{arnaudonhal-00826532}, it is proposed to replace it by the gradient
of some cost function
and then to follow the flow of this gradient.

%s1.5 #&#
\subsection{Discussion}\label{sec1.5}

The purpose of this paper is to propose a stochastic algorithm finding
$p$-means by a sequential use of samples from the underlying probability
measure on a Riemannian manifold $M$, even if the formal proof of its
convergence is only shown for the circle, the first non-trivial
example.\vspace*{1pt}

When $\nu$ is an empirical measure $(\sum_{l=1}^N\delta_{x_l})/N$,
where the $x_l$, $l\in\lin1, N\rin$, are points on the circle,
Charlier \cite{2011arXiv11091986C,2011arXiv11082141H} and
McKilliam, Quinn and Clarkson
\cite{McKilliam}
proposed algorithms finding the 2-mean with complexities of order $N\ln
(N)$ and $N$ for the latter work.
Empirical measures can in practice be used to approximate more general
probability measures on the circle, but it seems this is not a very
efficient method, since
for each new point added to the empirical measure, the whole algorithm
finding the corresponding mean has to be started again from scratch.
To our knowledge,
the process of Theorem~\ref{t1} is the only algorithm finding
$p$-means for any $p\geq1$ and for any probability measure $\nu$
admitting H\"olderian densities,
even in the restricted situation of the circle.

Another strong motivation for this paper is the treatment of the jumps
of the algorithms $X$ and~$Z$, situation
which is not covered by the techniques of \cite{MR1425361} (to the
contrary of the jumps of
the auxiliary process, which can be more easily dealt with).

In \cite{arnaudonhal-00826532}, we extend the ideas of the present
paper to the situation were $d^p(x,y)$
in (\ref{U}) is replaced by a quantity $\kappa(x,y)$ depending
smoothly on the parameters $x$ and $y$ belonging to a compact
Riemannian manifold $M$.
Via convolutions with the underlying heat kernel, it leads to an
algorithm enabling to deal with mappings $\kappa$
which are only assumed to be continuous.
But due to this regularization procedure, the corresponding algorithm
is less straightforward to put in practice than the one presented here.
Of course, the direct implementability has a price, since it needs
precise information about a crucial object,
$L^*_{\alpha,\beta}[\un]$. It will be defined in Section~\ref{sec3} and its
investigation has to be divided in several cases depending on the value
of $p$.
This is hidden in~\cite{arnaudonhal-00826532}, because we were more
interested there in the generalization to general compact manifolds
than in practicality considerations.

More technical discussions of the results are partially scattered over
the manuscript, when it seems
more appropriate to introduce them; see, for instance, Remarks~\ref
{discu1},~\ref{discu2},~\ref{discu3} and~\ref{t2p1}.

The paper is constructed on the following plan.
In next section, we recall some results about simulated annealing which give
the heuristics for the above convergence. Another alternative algorithm
is presented,
in the same spirit as $X$ and $Z$, but without jumps. In Section~\ref{sec3}, we discuss
about the regularity of the function $U_p$, in terms of that of $\nu$.
It enables to see how close is the instantaneous invariant measure
associated to
the algorithm at large times $t\geq0$ to the Gibbs measures associated to
the potential $U_p$ and to the inverse temperature $\beta_t^{-1}$.
The proof of Theorem~\ref{t1} is given in Section~\ref{sec4}.
The fifth section is devoted to the extension presented in Theorem~\ref{t2}
and the \hyperref[append]{Appendix} deals with technicalities relative to the temporal
marginal laws of the algorithms.

%s2 #&#
\section{Principles underlying the proof}\label{sec2}

Here, some results about the classical simulated annealing are reviewed.
The algorithm $X$ described in the \hyperref[sec1]{Introduction} will then appear as a
natural modification.
This will also give us the opportunity to present another intermediate
algorithm.

%s2.1 #&#
\subsection{Simulated annealing}\label{sec2.1}

Consider again $M$ a compact Riemannian manifold and denote $
\langle\cdot,\cdot\rangle$, $\nabla$, $\tr$ and $\lambda$
the corresponding scalar product, gradient, Laplacian operator and
probability measure.
Let $U$ be a given smooth function on $M$ to which we associate the
constant $b(U)\geq0$
defined similarly as in (\ref{bU}). We denote by $\cM$ the set of
global minima of $U$.

A corresponding simulated annealing algorithm $\theta\ffff (\theta
_t)_{t\geq0}$ associated
to a measurable inverse temperature scheme $\beta\st\RR_+\rightarrow
\RR_+$
is defined through the
evolution equation
\begin{eqnarray*}
\forall t\geq0,\qquad d\theta_t&=& dB_t-
\frac{\beta_t}{2}\nabla U(\theta_t) \,dt.
\end{eqnarray*}
It is a shorthand meaning that $\theta$ is a time-inhomogeneous Markov
process whose
generator at any time $t\geq0$ is $L_{\beta_t}$, where
%
%
%e2.1 #&#
\begin{eqnarray}
\label{Lb} \forall\beta\geq0,\qquad L_{\beta} \cdot&\ffff &
\tfrac{1}{2}\bigl(\tr\cdot-\beta\langle\nabla U,\nabla\cdot\rangle
\bigr).
\end{eqnarray}
Holley, Kusuoka and Stroock
\cite{MR995752} have proven the following result.
%

%th2.1 #&#
\begin{theo}\label{HKS}
For any fixed $T\geq1$, consider the inverse temperature scheme
\begin{eqnarray*}
\forall t\geq0,\qquad\beta_t&=& b^{-1}\ln(T+t),
\end{eqnarray*}
with $b>b(U)$. Then for any neighborhood $\cN$ of $\cM$ and for any
initial distribution
$\cL(\theta_0)$, we have
\begin{eqnarray*}
\lim_{t\rightarrow+\infty}\PP[\theta_t\in\cN]&=&1.
\end{eqnarray*}
\end{theo}

A crucial ingredient of the proof of this convergence are the Gibbs
measures associated to the potential $U$.
They are defined as the probability measures $\mu_\beta$ given for
any $\beta\geq0$ by
%
%
%e2.2 #&#
\begin{eqnarray}
\label{mub} \mu_\beta(dx)&\ffff &\frac{\exp(-\beta U(x))}{Z_\beta} \lambda(dx),
\end{eqnarray}
where $Z_\beta\ffff \int\exp(-\beta U(x)) \lambda(dx)$ is the
normalizing factor.

Indeed, \cite{MR995752} show that $\cL(\theta_t)$ and $\mu_{\beta_t}$
become closer and closer as $t\geq0$ goes to infinity, for instance,
in the sense
of total variation:
%
%
%e2.3 #&#
\begin{eqnarray}
\label{tv} \lim_{t\rightarrow+\infty}\bigl\llVert\cL(
\theta_t)-\mu_{\beta
_t}\bigr\rrVert_{\mathrm{tv}} &=&0.
\end{eqnarray}
Theorem~\ref{HKS} is then an immediate consequence of the fact that
for any neighborhood $\cN$ of $\cM$,
\begin{eqnarray*}
\lim_{\beta\rightarrow+\infty}\mu_{\beta}[\cN]&=&1.
\end{eqnarray*}
The constant $b(U)$
is critical for the behaviour (\ref{tv}), in the sense that if
we take
\begin{eqnarray*}
\forall t\geq0,\qquad\beta_t&=& b^{-1}\ln(T+t),
\end{eqnarray*}
with $T\geq1$ and $b<b(U)$, then there exist initial distributions
$\cL(\theta_0)$ such that
(\ref{tv}) is not true.

But in general (see, e.g., \cite{MR1623480}), the constant
$b(U)$ is not critical for Theorem~\ref{HKS}, the corresponding critical
constant being, with the notation of the \hyperref[sec1]{Introduction},
\begin{eqnarray*}
b'(U)&\ffff &\min_{x_0\in\cM}\max_{y\in M}U(x_0,y)-U(y)
\leq b(U)
\end{eqnarray*}
(compare with (\ref{bU2}), where $U$ replaces $U_p$ and where a global
minimum $x_0\in\cM$ is fixed).
Note that it may happen that $b'(U)=b(U)$, for instance, if $\cM$ has
only one connected component.

Another remark about Theorem~\ref{HKS} is that the convergence in
probability of $\theta_t$ for large $t\geq0$ toward $\cM$
cannot be improved into an almost sure convergence. Denote by $A$ the
connected component of
$\{x\in M\st U(x)\leq\min_M U+b\}$ which contains $\cM$ (the
condition $b>b(U)$ ensures that $\cM$
is contained in only one connected component of the above set). Then
almost surely,
$A$~is the limiting set of the trajectory $(\theta_t)_{t\geq0}$ (see
\cite{MR1348382},
where the corresponding result is proven for a finite state space but
whose proof
could be extended to the setting of Theorem~\ref{HKS}). We believe
that all these remarks should also hold in the context of
Conjecture~\ref{conj}
and Theorem~\ref{t1}.

%s2.2 #&#
\subsection{Heuristic of the proof}\label{sec2.2}

Let us now heuristically put forward why a result such as Conjecture
\ref{conj}
should be true, in relation with Theorem~\ref{HKS}. For simplicity of
the exposition, assume that $\nu$
is absolutely continuous with respect to $\lambda$.
For almost every $x,y\in M$,
there exists a unique minimal geodesic with speed 1 leading from $x$ to $y$.
Denote it by $(\gamma(x,y,t))_{t\in\RR}$, so that $\gamma(x,y,0)=x$ and
$\gamma(x,y,d(x,y))=y$.
The process $(X_t)_{t\geq0}$ underlying Theorem~\ref{HKS} is
Markovian and its inhomogeneous family of generators
is $(L_{\alpha_t,\beta_t})_{t\geq0}$, where for any $\alpha>0$ and
$\beta\geq0$, $L_{\alpha,\beta}$ acts
on functions $f$ from $\mathcal{C}^2(M)$ via, for all $x\in M$,
%
%
%e2.4 #&#
\begin{eqnarray}
\label{Lab} L_{\alpha,\beta}[f](x)&\ffff & \frac{1}2\tr f(x)+
\frac{1}{\alpha}\int f\bigl(\gamma\bigl(x,y,(p/2)\beta\alpha d^{p-1}(x,y)
\bigr)\bigr)-f(x) \nu(dy)
\end{eqnarray}
(to simplify notation, we will try to avoid writing down explicitly the
dependence on $p\geq1$).
The $r$ is well-defined, due to the fact that $\nu\ll\lambda$ which
implies that the cut-locus of $x$ is negligible with respect to $\nu$.
Furthermore Fubini's theorem enables to see that the function
$L_{\alpha,\beta}[f]$ is at least measurable.
Next, we remark that as $\alpha$ goes to $0_+$, we have for any $f\in
\mathcal{C}^1(M)$, any $x\in M$ and any
$y\in M$ which is not in the cut-locus of $x$,
\begin{eqnarray*}
\lim_{\alpha\rightarrow0_+}\frac{f(\gamma(x,y,(p/2)\beta\alpha
d^{p-1}(x,y)))-f(x)}{\alpha} &=&\frac{1}2\beta
pd^{p-1}(x,y) \bigl\langle\nabla f(x),\dot{\gamma}(x,y,0) \bigr\rangle,
\end{eqnarray*}
for all $\beta\geq0$, so that for any $f\in\mathcal{C}^2(M)$ and
$x\in M$,
\begin{eqnarray*}
\forall\beta\geq0,\qquad\lim_{\alpha\rightarrow0_+}L_{\alpha,\beta}[f](x)&=&
\frac{1}2\tr f(x)+\frac{\beta}{2} p\int d^{p-1}(x,y) \bigl
\langle\nabla f(x),\dot{\gamma}(x,y,0) \bigr\rangle\nu(dy).
\end{eqnarray*}
Recall that the potential $U=U_p$ we are now interested in is given by
(\ref{U})
and that for almost every $(x,y)\in M^2$,
\begin{eqnarray*}
\nabla_x d^p(x,y)&=&-pd^{d-1}(x,y)\dot{
\gamma}(x,y,0)
\end{eqnarray*}
(problems occur for points $x$ in the cut-locus of $y$ and, if $p=1$,
for $x=y$), thus
%
%
%e2.5 #&#
\begin{eqnarray}
\label{Uprime} \nabla U_p(x)&=&-p\int d^{p-1}(x,y)\dot{
\gamma}(x,y,0) \nu(dy).
\end{eqnarray}
It follows that
or any $f\in\mathcal{C}^2(M)$ and $x\in M$,
\begin{eqnarray*}
\forall\beta\geq0,\qquad\lim_{\alpha\rightarrow0_+}L_{\alpha
,\beta
}[f](x)&=&L_{\beta}[f](x).
\end{eqnarray*}
Since $\lim_{t\rightarrow+\infty}\alpha_t=0$, it appears that at
least for
large times,
$(X_t)_{t\geq0}$ and $(\theta_t)_{t\geq0}$ should behave in a
similar way.
The validity of Theorem~\ref{HKS} for any $T\geq1$ and any initial
distribution $\cL(\theta_0)$
then suggests that Conjecture~\ref{conj} should hold.
But this rough explanation is not sufficient to understand the choice
of the scheme
$(\alpha_t)_{t\geq0}$, which will require more rigorous computations
relatively to the corresponding homogenization property.
The heuristics for Theorem~\ref{t2} are similar, since the underlying algorithm
$(Z_t)_{t\geq0}$ is Markovian and its inhomogeneous family of generators
$(L_{\alpha_t,\beta_t,\kappa_t})_{t\geq0}$ satisfies
\begin{eqnarray*}
\forall f\in\mathcal{C}^2(M),\qquad\lim_{t\rightarrow+\infty
}
\bigl\llVert L_{\alpha
_t,\beta_t,\kappa_t}[f]-L_{\beta_t}[f]\bigr\rrVert
_{\infty}&=&0.
\end{eqnarray*}
For any $\alpha>0$, $\beta\geq0$ and $\kappa>0$, the generator
$L_{\alpha,\beta,\kappa}$ acts on
functions $f\in\mathcal{C}^2(M)$ via, for all $x\in M$,
\begin{eqnarray*}
L_{\alpha,\beta,\kappa}[f](x)&\ffff & \frac{1}2\tr f(x)+\frac{1}{\alpha
}\int f
\bigl(\gamma\bigl(x,z,(p/2)\beta\alpha d^{p-1}(x,z)\bigr)\bigr)-f(x)
K_{y,\kappa}(dz)\nu(dy).
\end{eqnarray*}

The previous observations suggest another possible algorithm to find
the mean of a probability measure
$\nu$ on $M$. Consider the $M\times M$-valued inhomogeneous Markov\vspace*{-2pt}
process $(\widetilde X_t,Y_{N^{(\alpha)}_t+1})_{t\geq0}$
where $(N^{(\alpha)}_t)_{t\geq0}$ was defined in (\ref{Na}) and where
%
%
%e2.6 #&#
\begin{eqnarray}
\label{sdewiX} \forall t\geq0,\qquad d\widetilde X_t&=&dB_t+(p/2)
\beta_t d^{p-1}(\widetilde X_t,Y_{N^{(\alpha)}_t+1})
\dot{\gamma}(\widetilde X_t,Y_{N^{(\alpha)}_t+1},0) \,dt.
\end{eqnarray}
Again,\vspace*{1pt} up to appropriate choices of the schemes $(\alpha_t)_{t\geq0}$
and $(\beta_t)_{t\geq0}$, it can be expected that
for any neighborhood $\cN$ of $\cM$ and for any initial distribution
$\cL(\widetilde X_0)$,
\begin{eqnarray*}
\lim_{t\rightarrow+\infty}\PP[\widetilde X_t\in\cN]&=&1.
\end{eqnarray*}
Indeed, this can be obtained by following the line of arguments
presented in \cite{MR1425361}; see \cite{arnaudonhal-00717677}.

But the main drawback of the algorithm $(\widetilde X_t)_{t\geq0}$ is that
theoretically,
it is asking for the computation
of the unit vector $\dot{\gamma}(\widetilde X_t,Y_{N^{(\alpha
)}_t+1},0)$ and of the
distance $d(\widetilde X_t,Y_{N^{(\alpha)}_t+1})$, at any time $t\geq
0$. From a
practical point of view,
its complexity will be bad in comparison with that of the algorithm
$X\ffff (X_t)_{t\geq0}$, which is not so difficult to implement, as it
was seen in Section~\ref{ni}.

%s2.3 #&#
\subsection{Outline of the proof}\label{sec2.3}

Since the Gibbs measure $\mu_\beta$ defined in \eqref{mub}, with $U$
replaced by $U_p$, concentrates on $\cM_p$ for large $\beta$,
it will be sufficient to show that the law $m_t$ of $X_t$ becomes
closer and closer to $\mu_{\beta_t}$ for large $t$.
To measure this closeness, we use the $\LL^2$-discrepancy of $m_t$
with respect to $\mu_{\beta_t}$ defined~by
\begin{eqnarray*}
\forall t>0,\qquad I_t&\ffff & \int\biggl(\frac{m_t}{\mu_{\beta_t}}-1
\biggr)^2 \,d\mu_{\beta
_t}.
\end{eqnarray*}
(Alternatively, it would be interesting to see if the considerations
that follow could be extended to the case where
this quantity is replaced by the more natural relative entropy of $m_t$
with respect to~$\mu_{\beta_t}$.)
To show that this quantity goes to zero as $t$ becomes large,
we study its temporal evolution, by differentiating it. The fact that
$\mu_{\beta_t}$ is not the instantaneous invariant measure
(namely the probability measure left invariant by the generator at time
$t$), leads to supplementary term
with respect to what one usually gets by applying this approach (see,
e.g., \cite{MR1275365}). This term measures in some sense the distance
between $\mu_{\beta_t}$ and the instantaneous invariant measure at
time $t$ (which itself is not explicitly known).
A large part of the paper is devoted to estimate this supplementary
term, the final result being presented in Proposition~\ref{mbl}.
In Proposition~\ref{Iprime}, we deduce a bound on the evolution of the
quantity $I_t$.
To conclude in Proposition~\ref{condab} that the obtained ordinary
differential inequality is sufficient
to conclude that $\lim_{t\rightarrow+\infty} I_t=0$, we need an
estimate of the
spectral gap of the operator presented in \eqref{Lb} for large $\beta$.
For that, we resort to a result due to \cite{MR995752}
recalled in Proposition
\ref{HKS2}.

Let us emphasize that the resort to the object $L^*_{\alpha,\beta
}[\un]$ defined and investigated in Section~\ref{ri} to estimate the
discrepancy
between a well-known measure and an instantaneous invariant measure,
which is more difficult to apprehend, should be
of much broader use than the one presented here. Indeed, the function
$L^*_{\alpha,\beta}[\un]$ is constructed
by using directly only two objects which are supposed to be known: the
generator and the convenient measure we choose to
replace the instantaneous invariant measure, because $L^*_{\alpha
,\beta}$ is just the dual operator of $L_{\alpha,\beta}$
in $\LL^2[\mu_\beta]$ and $\un$ is the constant function taking the
value 1.

%s3 #&#
\section{Regularity issues}\label{sec3}\label{ri}

From this section on, we restrict ourselves to the case of the circle.
Here, we investigate the regularity of the potential $U_p$ introduced
in (\ref{U})
and use the obtained information to evaluate how far are the
instantaneous invariant
measures of the algorithm $X$ from the corresponding Gibbs measures, as
well as some other
preliminary bounds.

For any $x\in\TT$, we denote $x'$ the unique point in the cut-locus
of $x$, namely
the opposite point $x'=x+\pi$.
Recall that for $y\in\TT\setminus\{x'\}$, $(\gamma(x,y,t))_{t\in
\RR}$
denotes the unique minimal geodesic with speed 1 going from $x$ to $y$
and that
$\delta_x$
stands for the Dirac mass at $x$.

%
%
%le3.1 #&#
\begin{lem}\label{regU}
For any probability measure $\nu$ on $\TT$, we have for the potential
$U_p$ defined in (\ref{U}), in the distribution sense,
for $x\in\TT$,
\begin{eqnarray*}
U_p''(x)&=&\cases{ \displaystyle p(p-1)\int
_\TT d^{p-2}(y,x)-2p\pi^{p-1}
\delta_{y'}(x) \nu(dy), &\quad if $p>1$,
\vspace*{3pt}\cr
\displaystyle 2\int\bigl(
\delta_y(x)-\delta_{y'}(x)\bigr) \nu(dy), &\quad if $p=1$.}
\end{eqnarray*}
In particular if $\nu$ admits a continuous density with respect to
$\lambda$,
still denoted $\nu$, then we have that $U_p\in\mathcal{C}^2(\TT)$ and
\begin{eqnarray*}
\forall x\in\TT,\qquad U_p''(x)= \cases{
\displaystyle
p(p-1)\int_\TT d^{p-2}(y,x) \nu(dy)-p
\pi^{p-2} \nu\bigl(x'\bigr), &\quad if $p>1$,
\vspace*{3pt}\cr
\displaystyle\bigl(
\nu(x)-\nu\bigl(x'\bigr)\bigr)/\pi, &\quad if $p=1$.}
\end{eqnarray*}
\end{lem}

\begin{pf}
We begin by considering the case where $p>1$. Furthermore, we first
investigate the situation where $\nu=\delta_y$ for some fixed $y\in
\TT$.
Then
$U_p(x)=d^p(x,y)$ for any $x\in\TT$ and we have seen in (\ref{Uprime})
that
\begin{eqnarray*}
\forall x\neq y',\qquad U_p'(x)&=&-pd^{p-1}(x,y)
\dot{\gamma}(x,y,0).
\end{eqnarray*}
By continuity of $U_p$, this equality holds in the sense of
distributions on the whole set $\TT$.
To compute $U_p''$, consider a test function $\varphi\in\mathcal
{C}^\infty(\TT)$:
\begin{eqnarray*}
\int_\TT\varphi'(x) U_p'(x)
\,dx &=&p\int_y^{y+\pi}\varphi'(x)
(x-y)^{p-1} \,dx-p\int_{y-\pi
}^y
\varphi'(x) (y-x)^{p-1} \,dx
\\
&=&p\bigl[\varphi(x) (x-y)^{p-1}\bigr]_y^{y+\pi}-p(p-1)
\int_y^{y+\pi}\varphi(x) (x-y)^{p-2} \,dx
\\
&&{}-p\bigl[\varphi(x) (y-x)^{p-1}\bigr]_{y-\pi}^y-p(p-1)
\int_{y-\pi}^y\varphi(x) (y-x)^{p-2} \,dx
\\
&=&2p\pi^{p-1}\varphi\bigl(y'\bigr)-p(p-1)\int
_\TT\varphi(x)d^{p-2}(y,x) \,dx.
\end{eqnarray*}
So we get that for $x\in\TT$,
\begin{eqnarray*}
U_p''(x)&=& p(p-1)d^{p-2}(y,x)-2p
\pi^{p-1} \delta_{y'}(x).
\end{eqnarray*}
If $p=1$, starting again from
%
%
%e3.1 #&#
\begin{eqnarray}
\label{Uprime1} \forall x\neq y',\qquad U_1'(x)&=&-
\dot{\gamma}(x,y,0),
\end{eqnarray}
we rather get for any test function $\varphi\in\mathcal{C}^\infty
(\TT)$:
\begin{eqnarray*}
\int_\TT\varphi'(x) U_1'(x)
\,dx &=&\int_y^{y+\pi}\varphi'(x) \,dx-
\int_{y-\pi}^y\varphi'(x) \,dx
\\
&=&2\bigl(\varphi\bigl(y'\bigr)-\varphi(y)\bigr),
\end{eqnarray*}
so that
\begin{eqnarray*}
U_1''&=& 2(\delta_y-
\delta_{y'}).
\end{eqnarray*}
The general case of a probability measure $\nu$ follows
by integration with respect to $\nu(dy)$.

The second announced result follows from the observation that if $\nu$ admits
a density with respect to $\lambda$, we can
write for any $x\in\TT$,
\begin{eqnarray*}
\int\delta_{y'}(x) \nu(dy)&=&\int\delta_{x'}(y)\nu(y)
\frac
{dy}{2\pi}
\\
&=&\frac{\nu(x')}{2\pi}.
\end{eqnarray*}\upqed
\end{pf}

In particular, it appears that the potential $U_p$ belongs to $\mathcal
{C}^\infty(\TT)$, if
the density $\nu$ is smooth.
Let us come back to the case of a general probability measure $\nu$ on
$\TT$.
For any $\alpha>0$ and $\beta\geq0$, we are interested into the
generator $L_{\alpha,\beta}$
defined in (\ref{Lab}). Rigorously speaking, this definition is only
valid if $\nu$
is absolutely continuous. Otherwise, the right-hand side of (\ref
{Lab}) is not well-defined for $x\in\TT$ belonging to the union of
the cut-locus of the atoms of $\nu$. To get around this little
inconvenience, one can consider for
$x\in\TT$, $(\gamma_+(x,x+\pi,t))_{t\in\RR}$ and $(\gamma
_-(x,x+\pi,t))_{t\in\RR}$, the unique minimal geodesics with speed 1
leading from $x$ to $x+\pi$, respectively, in the anti-clockwise (namely
increasing in the cover $\RR$ of $\TT$) and
clockwise direction. If $y\in\TT\setminus\{x'\}$, we take as before
$(\gamma_+(x,y,t))_{t\in\RR}\ffff (\gamma(x,y,t))_{t\in\RR}\fd
(\gamma_-(x,y,t))_{t\in\RR}$.
Next, let $k$ be a Markov kernel from $\TT^2$ to $\{-,+\}$ and modify
the definition
(\ref{Lab}) by imposing that for any $f\in\mathcal{C}^2(\TT)$ and
any $x\in\TT$,
\begin{eqnarray*}
L_{\alpha,\beta}[f](x)&\ffff & \frac{1}2\partial^2 f(x)+
\frac{1}{\alpha
}\int f\bigl(\gamma_s\bigl(x,y,(p/2)\beta\alpha
d^{p-1}(x,y)\bigr)\bigr)-f(x) k\bigl((x,y),ds\bigr)\nu(dy),
\end{eqnarray*}
where $\partial$ stands for
the natural derivative on $\TT$.
Then
the function
$L_{\alpha,\beta}[f]$ is at least measurable.
But these considerations are not very relevant, since for any given
measurable evolutions
$\RR_+\ni t\mapsto\alpha_t\in\RR_+^*$ and $\RR_+\ni t\mapsto
\beta_t\in\RR_+$,
the solutions to the martingale problems associated to the inhomogeneous
family of generators $(L_{\alpha_t,\beta_t})_{t\geq0}$ (see, e.g., the
book of \cite{MR838085}) are all the same and
are described in a probabilistic way as the trajectory laws of the
processes $X$ presented in the \hyperref[sec1]{Introduction}. Indeed, this is a
consequence of the absolute continuity of the heat kernel
at any positive time (for arguments in the same spirit, see the \hyperref[append]{Appendix}).
So to simplify notation, we only consider the case where $k((x,y),-)=0$
for any $x,y\in\TT$,
this brought us back to the definition (\ref{Lab}),
where $(\gamma(x,y,t))_{t\in\RR}$ stands for $(\gamma
_+(x,y,t))_{t\in\RR}$, for any $x,y\in\TT$.

As it was mentioned for usual simulated annealing algorithms in the
previous section,
a traditional approach to prove Theorem~\ref{t1} would try to evaluate
at any time $t\geq0$,
how far is $\cL(X_t)$ from the instantaneous invariant probability
$\mu_{\alpha_t,\beta_t}$, namely that associated to
$L_{\alpha_t,\beta_t}$.
Unfortunately, for any $\alpha>0$ and $\beta\geq0$,
we have little information about the invariant probability $\mu
_{\alpha,\beta}$
of $L_{\alpha,\beta}$, even its existence cannot be deduced directly
from the compactness
of $\TT$, because the functions $L_{\alpha,\beta}[f]$ are not
necessarily continuous for $f\in\mathcal{C}^2(\TT)$. Indeed it will
be more convenient to use the Gibbs distribution $\mu_\beta$ defined
in (\ref{mub})
for $\beta\geq0$, where $U$ is replaced by $U_p$. It has the
advantage to be explicit and easy to work with, in particular it is
clear that for large $\beta\geq0$, $\mu_\beta$ concentrates
around $\cM_p$, the set of $p$-means of $\nu$.

The remaining part of this section is mainly devoted to a
quantification of what separates $\mu_\beta$ from
being an invariant probability of $L_{\alpha,\beta}$, for
$\alpha>0$ and $\beta\geq0$.
It will become clear in the next section that a practical
way to measure this discrepancy is through the evaluation
of $\mu_\beta[(L^*_{\alpha,\beta}[\un])^2]$, where $L^*_{\alpha
,\beta}$ is the dual operator
of $L_{\alpha,\beta}$ in $\LL^2(\mu_\beta)$ and where $\un$ is
the constant function
taking the value 1. Indeed, it can be seen that $L^*_{\alpha,\beta
}[\un]=0$ in $\LL^2(\mu_\beta)$
if and only if $\mu_\beta$ is invariant for $L_{\alpha,\beta}$.
We will also take advantage of the computations made in this direction
to provide some estimates on related quantities
which will be helpful later on.

Since the situation of the usual mean $p=2$ is important and is simpler
than the other cases, we
first treat it in detail in the following subsection. Next, we will
investigate the differences appearing in
the situation of the median. The third subsection will deal with the
cases $1<p<2$, whose computations
are technical and not very enlightening. We will only give some
indications about the remaining situation
$p\in(2,\infty)$, which is less involved.

Some other preliminaries about the regularity of the time marginal laws
of the considered algorithms will be treated in the \hyperref[append]{Appendix}. They are
of a more qualitative nature and will mainly serve to justify some
computations of the next sections, in some sense they are less relevant
than the estimates and proofs of Propositions~\ref{mbl0},~\ref{mbl1},
\ref{mbl2} and~\ref{mbl3} below, which are really at the heart of our
developments.

%s3.1 #&#
\subsection{Estimate of \texorpdfstring{$L^*_{\alpha,\beta}[\un]$}{$L^*_{alpha,beta}[1]$} in the case $p=2$}\label{sec3.1}\label{p0}

Before being more precise about the definition of $L^*_{\alpha,\beta
}$, we need an elementary result, where we will use the following
notation: for $y\in\TT$ and $\delta\geq0$, $B(y,\delta)$ stands
for the open ball centered at $y$ of radius $\delta$ and for any $s\in
\RR$, $T_{y,s}$ is the operator acting on measurable functions $f$
defined on $\TT$ via
%
%
%e3.2 #&#
\begin{eqnarray}
\label{Tys}\forall x\in\TT,\qquad T_{y,s}f(x)&\ffff & f\bigl(\gamma
\bigl(x,y,sd(x,y)\bigr)\bigr).
\end{eqnarray}
%

%le3.2 #&#
\begin{lem}\label{gf}
For any $y\in\TT$, any $s\in[0,1)$ and any measurable and bounded
functions $f,g$, we have
\begin{eqnarray*}
\int_\TT gT_{y,s}f \,d\lambda&=&\frac{1}{1-s}
\int_{B(y,(1-s)\pi)} fT_{y,-s/(1-s)}g \,d\lambda.
\end{eqnarray*}
\end{lem}

\begin{pf}
By definition, we have
\begin{eqnarray*}
2\pi\int_\TT gT_{y,s}f \,d\lambda&=&\int
^{y+\pi}_{y-\pi} g(x) f\bigl(x+s(y-x)\bigr) \,dx.
\end{eqnarray*}
In the right-hand side, consider the change of variables
$z\ffff  sy+(1-s)x$ to get that
it is equal to
\begin{eqnarray*}
\frac{1}{1-s} \int^{y+(1-s)\pi}_{y-(1-s)\pi} g \biggl(
\frac{z-sy}{1-s} \biggr) f(z) \,dz &=& \frac{ 2\pi}{1-s} \int
_{B(y,(1-s)\pi)} fT_{y,-s/(1-s)}f \,d\lambda,
\end{eqnarray*}
which corresponds to the announced result.
\end{pf}

This lemma has for consequence the next result, where $\cD$ is the
subspace of $\LL^2(\lambda)$
consisting of functions whose second derivative in the distribution
sense belongs
to $\LL^2(\lambda)$ (or equivalently to $\LL^2(\mu_\beta)$ for any
$\beta\geq0$).
%

%le3.3 #&#
\begin{lem}\label{Ldual0}
For $\alpha>0$ and $\beta\geq0$ such that $\alpha\beta\in[0,1)$,
the domain of the
maximal extension of $L_{\alpha,\beta}$ on $\LL^2(\mu_\beta)$ is
$\cD$.
Furthermore, the domain $\cD^*$ of its dual operator $L_{\alpha,\beta
}^*$ in $\LL^2(\mu_\beta)$
is the space $\{f\in\LL^2(\mu_\beta)\st\exp(-\beta U_2)f\in\cD\}
$ and we have for any $f\in\cD^*$,
\begin{eqnarray*}
L_{\alpha,\beta}^*f &=& \frac{1}2 \exp(\beta U_2)
\partial^2\bigl[\exp(-\beta U_2)f\bigr]
\\
&&{}+\frac
{\exp
(\beta U_2)}{\alpha(1-\alpha\beta)}\int\un_{B(y,(1-\alpha\beta
)\pi)} T_{y,-\alpha\beta/(1-\alpha\beta)}\bigl[\exp(-\beta
U_2)f \bigr] \nu(dy)-\frac{f}{\alpha}.
\end{eqnarray*}
In particular, if $\nu$ admits a continuous density, then $\cD^*=\cD
$ and the above formula holds for any $f\in\cD$.
\end{lem}

\begin{pf}
With the previous definitions, we can write for any $\alpha>0$ and
$\beta\geq0$,
\begin{eqnarray*}
L_{\alpha,\beta}&=&\frac{1}2\partial^2+\frac{1}{\alpha}
\int T_{y,\alpha
\beta} \nu(dy)-\frac{I}{\alpha},
\end{eqnarray*}
where $I$ is the identity operator.
Note furthermore that the identity operator is bounded from $\LL
^2(\lambda)$
to $\LL^2(\mu_\beta)$ and conversely.
Thus, to get the first assertion, it is sufficient to show that
$\int T_{y,\alpha\beta} \nu(dy)$ is bounded from $\LL^2(\lambda)$
to itself,
or even only that $\llVert T_{y,\alpha\beta}\rrVert _{\LL
^2(\lambda )\righttoleftarrow}$
is uniformly bounded in $y\in\TT$.
To see that this is true, consider a bounded and measurable function
$f$ and assume that
$\alpha\beta\in[0,1)$.
Since $(T_{y,\alpha\beta}f)^2=T_{y,\alpha\beta}f^2$, we can apply
Lemma~\ref{gf}
with $s=\alpha\beta$, $g=\un$ and $f$ replaced by $f^2$ to get
that
\begin{eqnarray*}
\int(T_{y,\alpha\beta}f)^2 \,d\lambda&=& \frac{1}{1-\alpha\beta}\int
_{B(y,(1-s)\pi)} f^2T_{y,-\alpha\beta
/(1-\alpha\beta)}\un \,d\lambda
\\
&=& \frac{1}{1-\alpha\beta}\int_{B(y,(1-s)\pi)} f^2 \,d\lambda
\\
&\leq&\frac{1}{1-\alpha\beta}\int f^2 \,d\lambda.
\end{eqnarray*}
Next to see that for any $f,g\in\mathcal{C}^2(\TT)$,
%
%
%e3.3 #&#
\begin{eqnarray}
\label{Ladjoint} \int g L_{\alpha,\beta}f \,d\mu_\beta&=& \int f
L_{\alpha,\beta
}^*g \,d\mu_\beta,
\end{eqnarray}
where $L_{\alpha,\beta}^*$ is the operator defined in the statement
of the lemma,
we note that, on one hand,
\begin{eqnarray*}
\int g \partial^2f \,d\mu_\beta&=&Z_\beta^{-1}
\int\exp(-\beta U_2)g\partial^2f \,d\lambda
\\
&=& \int f \exp(\beta U_2)\partial^2\bigl[\exp(-\beta
U_2)g\bigr] \,d\mu_\beta
\end{eqnarray*}
and on the other hand, for any $y\in\TT$,
\begin{eqnarray*}
\int gT_{y,\alpha\beta}f \,d\mu_\beta&=& Z_\beta^{-1}
\int\exp(-\beta U_2)gT_{y,\alpha\beta}f \,d\lambda,
\end{eqnarray*}
so that we can use again Lemma~\ref{gf}.
After an additional integration with respect to $\nu(dy)$,
(\ref{Ladjoint}) follows without difficulty.
To conclude, it is sufficient to see that for any $f\in\LL^2(\mu
_\beta)$,
$L^*_{\alpha,\beta}f\in\LL^2(\mu_\beta)$ (where $L^*_{\alpha
,\beta}f$ is first interpreted as a distribution) if and only if $\exp
(-\beta U_2)f\in\cD$.
This is done by adapting the arguments given in the first part of the
proof, in particular we get
that
\begin{eqnarray*}
&& \biggl\llVert\frac{\exp(\beta U_2)}{\alpha(1-\alpha\beta)}\int\un
_{B(y,(1-\alpha\beta)\pi)} T_{y,-\alpha\beta/(1-\alpha\beta
)}\bigl[
\exp(-\beta U_2) \cdot\bigr] \nu(dy)\biggr\rrVert
^2_{\LL
^2(\lambda
)\righttoleftarrow}
\\
&&\quad \leq\frac{\exp(2\beta\osc(U_2))}{\alpha^2(1-\alpha
\beta)}.
\end{eqnarray*}\upqed
\end{pf}

%
%re3.1 #&#
\begin{rem}
By working in a similar spirit, the previous lemma,
except for the expression of $L_{\alpha,\beta}^*$, is valid for any
$\alpha>0$ and $\beta\geq0$ such that $\alpha\beta\neq1$.
The case $\alpha\beta=1$ can be different: it follows from
\begin{eqnarray*}
L_{\alpha,1/\alpha}&=&\frac{1}2\partial^2+\frac{1}{\alpha}(
\nu-I),
\end{eqnarray*}
that if $\nu$ does not admit a density with respect to $\lambda$
which belongs to $\LL^2(\lambda)$, then the domain of definition of
$L_{\alpha,1/\alpha}^*$ is
$
\cD^*\cap\{f\in\LL^2(\mu_\beta)\st\mu_\beta[f]=0\}$, subspace
which is not dense
in $\LL^2(\lambda)$ and worse for our purposes, which does not
contain $\un$.
Anyway, this degenerate situation is not very interesting for us,
because the evolutions
$(\alpha_t)_{t\geq0}$ and $(\beta_t)_{t\geq0}$ we consider satisfy
$\alpha_t\beta_t\in(0,1)$
for $t$ large enough. Furthermore, we will consider probability measures
$\nu$ admitting a continuous density, in particular
belonging to $\LL^2(\lambda)$. In this case, $L_{\alpha,1/\alpha}$
and $L_{\alpha,1/\alpha}^*$
admit $\cD$ for natural domain, as in fact $L_{\alpha,\beta}$ and
$L^*_{\alpha,\beta}$ for any $\beta\geq0$.
\end{rem}

For any $\alpha>0$ and $\beta\geq0$ such that $\alpha\beta\in[0,1)$,
denote $\eta=\alpha\beta/(1-\alpha\beta)$.
As seen from the previous lemma, a consequence of the assumption that
$U_2$ is $\mathcal{C}^2$ is that for any $x\in\TT$,
\begin{eqnarray}
\label{Lstarun}
L_{\alpha,\beta}^*\un(x)
&=& \frac{1}2 \exp\bigl(\beta U_2(x)\bigr)
\partial^2\exp\bigl(-\beta U_2(x)\bigr)-\frac
{1}{\alpha
}\nonumber
\\
&&{} +
\nonumber
\frac{\exp(\beta U_2(x))}{\alpha(1-\alpha\beta)}\int\un_{B(y,(1-\alpha
\beta)\pi)}(x) T_{y,-\eta}
\bigl[\exp(-\beta U_2)\bigr](x) \nu(dy)
\nonumber\\[-8pt]\\[-8pt]\nonumber
&=& \frac{\beta^2}2\bigl(U_2'(x)
\bigr)^2-\frac{\beta}2U_2''(x)-
\frac
{1}{\alpha}
\\
&&{} +\frac{1}{\alpha(1-\alpha\beta)}\int_{B(x,(1-\alpha\beta)\pi)} \exp
\bigl(\beta
\bigl[U_2(x)-U_2\bigl(\gamma\bigl(x,y,-\eta d(x,y)\bigr)
\bigr)\bigr]\bigr) \nu(dy).\nonumber
\end{eqnarray}
It appears that $L_{\alpha,\beta}^*\un$ is defined and continuous if
$\nu$ has a continuous density (with respect to~$\lambda$). The next result
evaluates the uniform norm of this function under a little stronger
regularity assumption.
Despite it may seem quite plain, we would like to emphasize that the
use of an estimate of
$L_{\alpha,\beta}^*\un$ to replace the invariant measure of
$L_{\alpha,\beta}$
by the more tractable $\mu_\beta$ is a key to all the results
presented in the \hyperref[sec1]{Introduction}.
%

%pr3.1 #&#
\begin{pro}\label{mbl0}
Assume that $\nu$ admits a density with respect to $\lambda$ which is
H\"older continuous, that is, there exists $a\in(0,1]$ and $A>0$ such that
%
%
%e3.4 #&#
\begin{eqnarray}
\label{aA} \forall x,y\in\TT,\qquad\bigl\llvert\nu(y)-\nu(x)\bigr
\rrvert&
\leq& Ad^a(x,y).
\end{eqnarray}
Then there exists a constant $C(A)>0$, only depending on $A$, such that
for any $\beta\geq1$ and $\alpha\in(0,1/(2\beta^2))$, we have
\begin{eqnarray*}
\bigl\llVert L_{\alpha,\beta}^*\un\bigr\rrVert_\infty&\leq&C(A)\max
\bigl(\alpha\beta^4,\alpha^a\beta^{1+a} \bigr).
\end{eqnarray*}
\end{pro}

\begin{pf}
In view of the expression of $ L_{\alpha,\beta}^*\un(x)$ given before
the statement of the proposition, we want to estimate for any fixed
$x\in\TT$,
the quantity
\begin{eqnarray*}
&& \int_{B(x,(1-\alpha\beta)\pi)} \exp\bigl(\beta\bigl[U_2(x)-U_2
\bigl(\gamma\bigl(x,y,-\eta d(x,y)\bigr)\bigr)\bigr]\bigr) \nu(dy)
\\
&&\quad = \int_{x-(1-\alpha\beta)\pi}^{x+(1-\alpha\beta)\pi} \exp\bigl(\beta
\bigl[U_2(x)-U_2\bigl(x-\eta(y-x)\bigr)\bigr]\bigr)
\nu(dy).
\end{eqnarray*}
Lemma~\ref{regU} and the continuity of the density $\nu$ ensure that
$U_2\in\mathcal{C}^2(\TT)$. Furthermore, since this density takes
the value 1 somewhere on $\TT$,
we get that
%
%
%e3.5 #&#
\begin{eqnarray}
\label{Usec} \bigl\llVert U_2''\bigr
\rrVert_\infty\leq2A\pi^a \leq2\pi A.
\end{eqnarray}
Since $U_2'$ vanishes somewhere on $\TT$, we can deduce from
this bound that $\llVert U_2'\rrVert _\infty\leq4\pi^2 A$,
but for
$A>1/(2\pi)$, it is better
to use (\ref{Uprime}), which gives directly $\llVert U_2'\rrVert
_\infty\leq
2\pi$.

Expanding the function $U_2$ around $x$, we see
that for any $y\in(x-(1-\alpha\beta)\pi,x+(1-\alpha\beta)\pi)$
and $\eta\in(0,1]$
(this is satisfied because the assumptions on $\alpha$ and $\beta$
ensure that $\alpha\beta\in(0,1/2)$),
we can find $z\in(x-(1-\alpha\beta)\pi,x+(1-\alpha\beta)\pi)$
such that
\begin{eqnarray*}
\beta\bigl[U_2(x)-U_2\bigl(x-\eta(y-x)\bigr)\bigr]&=&
\beta\eta U_2'(x) (y-x)-\beta\eta^2
U_2''(z)\frac{(y-x)^2}{2}.
\end{eqnarray*}
The last term can be written under the form $\cO_A(\alpha^2\beta
^3)$, where for any $\epsilon>0$, $\cO_A(\epsilon)$
designates a quantity which is bounded by $K(A)\epsilon$, where $K(A)$
is a constant
depending only on $A$ (as usual $\cO$ has a similar meaning, but with
a universal constant).
Note that we also have $\beta\eta U_2'(x)(y-x)=\cO(\alpha\beta^2)$.
Observing that for any $r,s\in\RR$, we can find $u,v\in(0,1)$ such that
$\exp(r+s)=(1+r+r^2\exp(ur)/2)(1+s\exp(vs))$
and in conjunction with the assumption $\alpha\beta^2\leq1/2$, we
can write that
%
%
%e3.6 #&#
\begin{eqnarray}
\label{expUU} \exp\bigl(\beta\bigl[U_2(x)-U_2\bigl(x-
\eta(y-x)\bigr)\bigr]\bigr)&=&1+ \beta\eta U_2'(x)
(y-x)+\cO_A\bigl(\alpha^2\beta^4\bigr).
\end{eqnarray}
Integrating this expression, we get that
\begin{eqnarray*}
&& \int_{B(x,(1-\alpha\beta)\pi)} \exp\bigl(\beta\bigl[U_2(x)-U_2
\bigl(\gamma(x,y,-\eta)\bigr)\bigr]\bigr) \nu(dy)
\\
&&\quad = \nu\bigl[B\bigl(x,(1-\alpha\beta)\pi\bigr)\bigr]+ \beta\eta
U_2'(x) \int_{x-(1-\alpha\beta)\pi}^{x+(1-\alpha\beta
)\pi}y-x
\nu(dy) +\cO_A\bigl(\alpha^2\beta^4\bigr).
\end{eqnarray*}
Recalling that $\nu$ has no atom, the first term
is equal to $1-\nu(B(x',\alpha\beta\pi)$.
Taking into account
(\ref{Uprime}),
we have $U_2'(x)=-2\int_{x-\pi}^{x+\pi}y-x \nu(dy)$, so that
the second term is
equal to
\begin{eqnarray*}
&& \beta\eta U_2'(x)\int
_{x-\pi}^{x+\pi}y-x \nu(dy) - \beta\eta
U_2'(x)\int_{x'-\alpha\beta\pi}^{x'+\alpha\beta\pi
}y-x
\nu(dy)
= -\frac{\beta\eta}2 \bigl(U_2'(x)
\bigr)^2+\cO_A\bigl(\alpha^2
\beta^3\bigr)
\end{eqnarray*}
(in the last term of the left-hand side, $y-x$ is to be interpreted as
its representative in $(-\pi,\pi]$ modulo $2\pi$).
We can now return to (\ref{Lstarun}) and recalling the expression for $U_2''$
given in Lemma~\ref{regU}, we obtain that for any $x\in\TT$,
\begin{eqnarray*}
L_{\alpha,\beta}^*\un(x)&=& \frac{\beta^2}2\bigl(U_2'(x)
\bigr)^2-\beta\bigl(1-\nu\bigl(x'\bigr)\bigr)-
\frac{1}{\alpha}
\\
&&{} +\frac{1}{\alpha(1-\alpha\beta)} \biggl(1-\nu
\bigl(B\bigl(x',\alpha\beta\pi
\bigr)\bigr)-\frac{\beta\eta}2 \bigl(U_2'(x)
\bigr)^2 +\cO_A\bigl(\alpha^2
\beta^4\bigr) \biggr)
\\
&=&\frac{1}{\alpha(1-\alpha\beta)}-\beta-\frac{1}{\alpha} +\frac{\beta
^2}{2} \biggl(1-
\frac{1}{(1-\alpha\beta)^2} \biggr) \bigl(U_2'(x)
\bigr)^2
\\
&&{}+\beta\biggl(\nu\bigl(x'\bigr)-\frac{\nu(B(x',\alpha\beta\pi))}{\alpha
\beta(1-\alpha\beta)} \biggr)+
\cO_A\bigl(\alpha\beta^4\bigr)
\\
&=&\beta\biggl(\nu\bigl(x'\bigr)-\frac{\nu(B(x',\alpha\beta\pi))}{\alpha
\beta(1-\alpha\beta)} \biggr) +
\cO_A\bigl(\alpha\beta^4\bigr)
\\
&=& \frac{\beta}{1-\alpha\beta} \biggl(\nu\bigl(x'\bigr)-\frac{\nu
(B(x',\alpha
\beta\pi))}{\alpha\beta}
\biggr) -\frac{\alpha\beta^2}{1-\alpha\beta}\nu\bigl(x'\bigr) +\cO_A
\bigl(\alpha\beta^4\bigr)
\\
&=&\frac{\beta}{1-\alpha\beta} \frac{1}{2\pi\alpha\beta}\int_{x'-\alpha
\beta\pi}^{x'+\alpha
\beta\pi}
\nu\bigl(x'\bigr)-\nu(y) \,dy +\cO_A\bigl(\alpha
\beta^4\bigr).
\end{eqnarray*}
The justification of the H\"older continuity comes above all from the evaluation
of the latter integral:
\begin{eqnarray*}
\biggl\llvert\int_{x'-\alpha\beta\pi}^{x'+\alpha\beta\pi} \nu
\bigl(x'\bigr)-\nu(y) \,dy\biggr\rrvert&\leq& A\int
_{x'-\alpha\beta\pi}^{x'+\alpha\beta\pi} \bigl\llvert x'-y\bigr
\rrvert^a \,dy
\\
&=& 2A\frac{(\alpha\beta\pi)^{1+a}}{1+a}
\\
&\leq&2A(\alpha\beta\pi)^{1+a}.
\end{eqnarray*}
The bound announced in the lemma follows at once.
\end{pf}

To finish this subsection, let us present a related but more
straightforward preliminary bound.
%

%le3.4 #&#
\begin{lem}\label{Tstarpa0}
There exists a constant $k>0$ such that
for any $s>0$ and $\beta\geq1$ with $\beta s\leq1/2$, we have,
for any $y\in\TT$ and $f\in\mathcal{C}^1(\TT)$,
%
%
%e3.7 #&#
\begin{eqnarray}
\label{borned0}
\int_{B(y,(1-s)\pi)} \bigl(T^*_{y,s}[g_y](x)-g_y(x)
\bigr)^2 \mu_\beta(dx)&\leq& ks^2
\beta^2 \biggl(\int(\partial f)^2 \,d\mu_{\beta}+
\int f^2 \,d\mu_\beta\biggr),\quad
\end{eqnarray}
where $T^*_{y,s}$
is the adjoint operator of $T_{y,s}$ in $\LL^2(\mu_\beta)$ and
where for any fixed $y\in\TT$,
\begin{eqnarray*}
\forall x\in\TT\setminus\bigl\{y'\bigr\},\qquad
g_y(x)&\ffff & f(x)d(x,y)\dot{\gamma}(x,y,0)
\end{eqnarray*}
(neglecting the cut-locus point $y'$ of $y$).
\end{lem}

\begin{pf}
Since the problem is clearly invariant by translation of $y\in\TT$,
we can work with a fixed value of $y$, the most convenient to simplify
the notation
being $y=0\in\RR/(2\pi\ZZ)$.
Then the function $g\equiv g_0$ is given by $g(x)=-xf(x)$ for $x\in
(-\pi,\pi)$.

Due to the above assumptions, $s\in(0,1/2)$ and we are in position to
use Lemma~\ref{gf} to see that for $s\in(0,1/2)$ and for
a.e. $x\in(-(1-s)\pi,(1-s)\pi)$,
\begin{eqnarray*}
T^*_{s}[g](x)&=&\frac{1}{1-s}\exp\bigl(\beta U_2(x)
\bigr)T_{-\eta}\bigl[\exp(-\beta U_2)g\bigr](x),
\end{eqnarray*}
with $\eta\ffff  s/(1-s)$ and where we simplified notation by replacing
$T^*_{0,s}$ and $T_{0,-\eta}$ by $T^*_s$ and
$T_{-\eta}$.
This observation induces us to introduce on $(-(1-s)\pi,(1-s)\pi)$
the decomposition
\begin{eqnarray*}
T^*_{s}[g]-g&=&T^*_{s}[g]-\frac{1}{1-s}T_{-\eta}[g]+
\frac
{1}{1-s}\bigl(T_{-\eta}[g]-g\bigr) +\frac{s}{1-s}g,
\end{eqnarray*}
leading to
%
%
%e3.8 #&#
\begin{eqnarray}
\label{J1230} \int\bigl(T^*_{s}[g](x)-g(x)\bigr)^2
\mu_\beta(dx) &\leq&\frac{3}{(1-s)^2}J_1+
\frac{3}{(1-s)^2}J_2+\frac{3s^2}{(1-s)^2}J_3,
\end{eqnarray}
where
\begin{eqnarray*}
J_1&\ffff & \int_{-(1-s)\pi}^{(1-s)\pi} \bigl(\exp
\bigl(\beta\bigl[U_2(x)-U_2\bigl((1+\eta)x\bigr)\bigr]
\bigr)-1\bigr)^2\bigl(T_{-\eta}[g] \bigr)^2
\mu_\beta(dx),
\\
J_2&\ffff & \int_{-(1-s)\pi}^{(1-s)\pi}
\bigl(T_{-\eta}[g] -g\bigr)^2 \,d\mu_\beta,
\\
J_3&\ffff & \int_{-(1-s)\pi}^{(1-s)\pi}
g^2 \,d\mu_\beta.
\end{eqnarray*}

The simplest term to treat is $J_3$: we just bound it above by $\int
g^2 \,d\mu_\beta$.
Recalling that $g\leq\pi^2 f^2$, we end up with a bound which goes in
the direction of (\ref{borned0}), due to the factor ${3s^2}/{(1-s)^2}$
in~(\ref{J1230}) and the fact that $\beta\geq1$.

Next, we estimate the term $J_1$.
Via the change of variable $z\ffff (1+\eta)x$, Lemma~\ref{gf} enables
to write it down under the form
\begin{eqnarray*}
&& (1-s)\int_{\TT} \bigl(\exp\bigl(\beta
\bigl[U_2\bigl((1-s)z\bigr)-U_2(z)\bigr]\bigr)-1
\bigr)^2g^2(z)\exp(\beta\bigl[U_2(z)-U_2
\bigl((1-s)z\bigr)\bigr] \mu_\beta(dz)
\\
&&\quad=4(1-s) \int_{\TT}\sinh^2\bigl(\beta
\bigl[U_2\bigl((1-s)z\bigr)-U_2(z)\bigr]/2
\bigr)g^2(z) \mu_\beta(dz).
\end{eqnarray*}
Since $\beta s\leq1/2$, we are assured of the bounds
%
%
%e3.9 #&#
\begin{eqnarray}
\label{diffU2}
\nonumber
\bigl\llvert\beta\bigl[U_2\bigl((1-s)z
\bigr)-U_2(z)\bigr]\bigr\rrvert&\leq&\beta\bigl\llVert
U_2'\bigr\rrVert_\infty\pi s
\\
&\leq& 4\pi^2\beta s
\\
&\leq&2\pi^2\nonumber
\end{eqnarray}
and we deduce that
\begin{eqnarray*}
J_1&\leq&16\pi^4 \cosh^2\bigl(
\pi^2\bigr)\beta^2s^2\int g^2 \,d\mu_\beta.
\end{eqnarray*}
Again this bound is going in the direction of (\ref{borned0}).

We are thus left with the task of finding a bound on $J_2$
and this is where the Dirichlet type quantity $\int(f')^2 \,d\mu_\beta
$ will be needed.
Of course, its origin is to be found in the fundamental theorem of
calculus, which enables to
write for any $x\in(-(1-s)\pi, (1-s)\pi)$,
\begin{eqnarray*}
T_{-\eta}[g] (x)-g(x)&=&-\eta\int_0^1g'
\bigl((1+\eta v)x\bigr)x \,dv.
\end{eqnarray*}
It follows that
%
%
%e3.10 #&#
\begin{eqnarray}
\label{J1} J_2&\leq& \pi^2\eta^2 \int
_{-(1-s)\pi}^{(1-s)\pi}\mu_\beta(dx)\int
_0^1\,dv \bigl( g'\bigl((1+\eta v)x
\bigr) \bigr)^2.
\end{eqnarray}
Recalling the definition of $g$, we have for any $z\in(-\pi,\pi)$,
\begin{eqnarray*}
\bigl(g'(z)\bigr)^2&\leq& 2\bigl(\pi^2
\bigl(f'(z)\bigr)^2+f^2(z)\bigr),
\end{eqnarray*}
where we used again that $\llVert U_2'\rrVert _\infty\leq
2\pi$ and that
$\beta\geq1$.
Next, we deduce from a computation similar to (\ref{diffU2}) and from
$\eta\leq2s$ that
\begin{eqnarray*}
\frac{\mu_\beta(x)}{\mu_\beta((1+\eta v)x)} &\leq& \exp\bigl(4\pi^2\bigr),
\end{eqnarray*}
so it appears that there exists a universal constant $k_1>0$ such that
\begin{eqnarray*}
\int_{-(1-s)\pi}^{(1-s)\pi}\mu_\beta(dx)\int
_0^1\,dv \bigl( g'\bigl((1+\eta v)x
\bigr) \bigr)^2&\leq&k_1 \int_0^1 dv \int_{-(1-s)\pi}^{(1-s)\pi}\lambda(dx) T_{-\eta
v}[h](x),
\end{eqnarray*}
where
\begin{eqnarray*}
\forall x\in\TT,\qquad h(x)&\ffff &\bigl[\bigl(f'(x)
\bigr)^2+f^2(x)\bigr]\mu_\beta(x).
\end{eqnarray*}
The proof of Lemma~\ref{gf} shows that for any fixed $v\in[0,1]$,
\begin{eqnarray*}
\int_{-(1-s)\pi}^{(1-s)\pi}T_{-\eta v}[h](x) \lambda(dx)
&\leq& \frac{1}{1+v\eta}\int_\TT h(x) \lambda(dx)
\\
&\leq& \int_\TT h(x) \lambda(dx)
\\
&=&\int_\TT\bigl(f'\bigr)^2 \,d\mu_\beta+\int_\TT f^2 \,d
\mu_\beta.
\end{eqnarray*}
Coming back to (\ref{J1}) and recalling that $\eta=s/(1-s)$, we
obtain that
\begin{eqnarray*}
J_2&\leq& k_2s^2 \biggl(\int
_\TT\bigl(f'\bigr)^2 \,d
\mu_\beta+\int_\TT f^2 \,d
\mu_\beta\biggr),
\end{eqnarray*}
for another universal constant $k_2>0$. This ends the proof of (\ref{borned0}).
\end{pf}

%s3.2 #&#
\subsection{Estimate of \texorpdfstring{$L^*_{\alpha,\beta}[\un]$}{$L^*_{alpha,beta}[1]$} in the case $p=1$}\label{sec3.2}\label{p1}

When we are interested in finding medians, the definition (\ref{Tys})
must be modified into
%
%
%e3.11 #&#
\begin{eqnarray}
\label{Tys1}\forall x\in\TT,\qquad T_{y,s}f(x)&\ffff & f\bigl(\gamma
(x,y,s)\bigr).
\end{eqnarray}
Similarly to what we have done in Lemma~\ref{gf}, we begin by
computing the adjoint $T^{\dagger}_{y,s}$ of $T_{y,s}$ in $\LL
^2(\lambda)$,
for any fixed $y\in\TT$ and $s\in\RR_+$ small enough.
%

%le3.5 #&#
\begin{lem}\label{dagger1}
Assume that $s\in[0,\pi/2)$.
Then for any bounded and measurable function $g$, we have, for almost
every $x\in\TT$ (identified with its representative in $(y-\pi,y+\pi)$),
\begin{eqnarray*}
T^{\dagger}_{y,s}[g](x) &=& \un_{(y-\pi+s,y-s)}(x)g(x-s)+
\un_{(y-s,y+s)}(x) \bigl(g(x-s)+g(x+s)\bigr)
\\
&&{}+ \un_{(y+s,y+\pi-s)}(x)g(x+s).
\end{eqnarray*}
\end{lem}

\begin{pf}
By definition, we have, for any bounded and measurable functions $f,g$,
\begin{eqnarray*}
2\pi\int_\TT gT_{y,s}f \,d\lambda&=&\int
^{y+\pi}_{y-\pi} g(x) f\bigl(x+\sign(y-x)s\bigr) \,dx.
\end{eqnarray*}
Let us first consider the integral
\begin{eqnarray*}
\int^{y+\pi}_{y} g(x) f\bigl(x+\sign(y-x)s\bigr) \,dx
&=&\int^{y+\pi}_{y} g(x) f(x-s) \,dx
\\
&=&\int^{y+\pi-s}_{y-s} g(x+s) f(x) \,dx
\\
&=&\int^{y+\pi-s}_{y+s} g(x+s) f(x) \,dx+\int
^{y+s}_{y-s} g(x+s) f(x) \,dx.
\end{eqnarray*}
The symmetrical computation on $(y-\pi,y)$ leads to the announced result.
\end{pf}

It is not difficult to adapt the proof of Lemma~\ref{Ldual0}, to get,
with the same notation,
%

%le3.6 #&#
\begin{lem}\label{Ldual1}
For $\alpha>0$ and $\beta\geq0$ such that $\alpha\beta\in[0,\pi
)$, the domain of the
maximal extension of $L_{\alpha,\beta}$ on $\LL^2(\mu_\beta)$ is
$\cD$.
Furthermore, the domain of its dual operator $L_{\alpha,\beta}^*$ in
$\LL^2(\mu_\beta)$
is $\cD^*$ and we have for any $f\in\cD^*$,
\begin{eqnarray*}
L_{\alpha,\beta}^*f &=& \frac{1}2 \exp(\beta U_1)
\partial^2\bigl[\exp(-\beta U_1)f\bigr]+
\frac{1}\alpha\int T_{y,\vfrac{\alpha\beta}{2}}^{*}[f] \nu(dy)-
\frac{f}{\alpha},
\end{eqnarray*}
where
\begin{eqnarray*}
T_{y,\vfrac{\alpha\beta}{2}}^{*}[f]&=& \exp(\beta U_1)
T_{y,\vfrac{\alpha\beta}{2}}^{\dagger}\bigl[\exp(-\beta U_1)f \bigr].
\end{eqnarray*}
In particular, if $\nu$ admits a continuous density, then $\cD^*=\cD
$ and the above formula holds for any $f\in\cD$.
\end{lem}

To be able to consider $L_{\alpha,\beta}^*\un$, we have thus to
assume that $\nu$ admits a continuous density, so that
$\un\in\cD^*=\cD$.
Furthermore, we obtain then that for almost every $x\in\TT$,
\begin{eqnarray*}
L_{\alpha,\beta}^*\un(x)&=&\frac{\beta^2}{2}\bigl(U_1'(x)
\bigr)^2-\frac
{\beta}{2} U_1''(x)+
\frac{1}{\alpha} \biggl(\int T_{y,\vfrac{\alpha\beta}{2}}^{*}[\un](x) \nu(dy)-1
\biggr).
\end{eqnarray*}
By expanding the various terms of the right-hand side, we are to show
the equivalent of Proposition~\ref{mbl0}.

%
%pr3.2 #&#
\begin{pro}\label{mbl1}
Assume that $\nu$ admits a density with respect to $\lambda$
satisfying (\ref{aA}).
Then there exists a constant $C(A)>0$, only depending on $A$, such that
for any $\beta\geq1$ and $\alpha\in(0,\pi\beta^{-2})$, we have
\begin{eqnarray*}
\bigl\llVert L_{\alpha,\beta}^*\un\bigr\rrVert_\infty&\leq&C(A)\max
\bigl(\alpha\beta^4,\alpha^a\beta^{1+a} \bigr).
\end{eqnarray*}
\end{pro}

\begin{pf}
From (\ref{Uprime}) and Lemma~\ref{regU}, we deduce, respectively, that
for all $x\in\TT$,
%
%
%e3.12 #&#
\begin{eqnarray}
\nonumber
U_1'(x)&=&-\int\dot{
\gamma}(x,y,0) \nu(dy)
\nonumber\\[-8pt]\label{Uprime1a} \\[-8pt]\nonumber
&=&\nu\bigl((x-\pi,x)\bigr)-\nu\bigl((x,x+\pi)\bigr),
\\
\label{Usecond1} U_1''(x)&=&\bigl(\nu(x)-
\nu\bigl(x'\bigr)\bigr)/\pi.
\end{eqnarray}
On the other hand, from Lemma~\ref{dagger1} we get that for all $s\in
[0,\pi/2)$ and for almost every $x\in\TT$,
\begin{eqnarray*}
&& \int T_{y,s}^{*}[\un](x) \nu(dy)
\\
&&\quad = \nu\bigl((x+s,x+\pi-s)\bigr)\exp\bigl(\beta\bigl(U_1(x)-U_1(x-s)
\bigr)\bigr)
\\
&&\qquad{} +\nu\bigl((x-s,x+s)\bigr)\bigl[\exp\bigl(\beta\bigl(U_1(x)-U_1(x-s)
\bigr)\bigr)+\exp\bigl(\beta\bigl(U_1(x)-U_1(x+s)\bigr)
\bigr)\bigr]
\\
&&\qquad{} +\nu\bigl((x-\pi+s,x-s)\bigr)\exp\bigl(\beta\bigl(U_1(x)-U_1(x+s)
\bigr)\bigr)
\\
&&\quad = \nu\bigl((x,x+\pi)\bigr)\exp\bigl(\beta\bigl(U_1(x)-U_1(x-s)
\bigr)\bigr)+\nu\bigl((x-\pi,x)\bigr)\exp\bigl(\beta\bigl(U_1(x)-U_1(x+s)
\bigr)\bigr)
\\
&&\qquad{} +\nu\bigl((x-s,x)\bigr)\exp\bigl(\beta\bigl(U_1(x)-U_1(x-s)
\bigr)\bigr)
\\
&&\qquad {}+\nu\bigl((x,x+s)\bigr)\exp\bigl(\beta\bigl(U_1(x)-U_1(x+s)
\bigr)\bigr)
\\
&&\qquad{} -\nu\bigl(\bigl(x'-s,x'\bigr)\bigr)\exp\bigl(
\beta\bigl(U_1(x)-U_1(x-s)\bigr)\bigr)-\nu\bigl(
\bigl(x',x'+s\bigr)\bigr)
\\
&&\qquad{}\times \exp\bigl(\beta
\bigl(U_1(x)-U_1(x+s)\bigr)\bigr).
\end{eqnarray*}
This leads us to define $s=\alpha\beta/2\in(0,\pi/2)$, so that we
can decompose
\begin{eqnarray*}
\frac{2}{\beta}L_{\alpha,\beta}^*\un(x)&=& I_1(x,s)+I_2(x,s)+I_3(x,s),
\end{eqnarray*}
with
\begin{eqnarray*}
I_1(x,s)&\ffff & \frac{1}\pi\biggl(\pi\frac{\nu((x-s,x+s))}{s}-
\nu(x) \biggr)-\frac{1}\pi\biggl(\pi\frac{\nu((x'-s,x'+s))}{s}-\nu
\bigl(x'\bigr) \biggr),
\\
I_2(x,s)&\ffff & \frac{\nu((x-s,x))-\nu((x'-s,x'))}{s}\bigl[\exp\bigl(\beta
\bigl(U_1(x)-U_1(x-s)\bigr)\bigr)-1\bigr]
\\
&&{}+ \frac{\nu((x,x+s))-\nu((x',x'+s))}{s}\bigl[\exp\bigl(\beta\bigl
(U_1(x)-U_1(x+s)
\bigr)\bigr)-1\bigr],
\\
I_3(x,s)&\ffff & \nu\bigl((x,x+\pi)\bigr)\frac{\exp(\beta
(U_1(x)-U_1(x-s)))-1-s\beta U_1'(x)}{s}
\\
&&{}+\nu\bigl((x-\pi,x)\bigr)\frac{\exp
(\beta(U_1(x)-U_1(x+s)))-1+s\beta U_1'(x)}{s}.
\end{eqnarray*}
Assumption (\ref{aA}) enables to evaluate $I_1(x,s)$, because we have
for any $x\in\TT$ and $s\in(0,\pi/2)$,
\begin{eqnarray*}
\biggl\llvert\pi\frac{\nu((x-s,x+s))}{s}-\nu(x)\biggr\rrvert&= & \frac
{1}{2 s}
\biggl\llvert\int_{(x-s,x+s)}\nu(z)-\nu(x) \,dz\biggr\rrvert
\\
&\leq& \frac{A}{2 s} \int_{(x-s,x+s)} \llvert z-x\rrvert
^a \,dz
\\
&=&\frac{As^a}{1+a}
\\
&\leq&As^a.
\end{eqnarray*}
By considering the Taylor's expansion with remainder at the first order
of the mapping $s\mapsto\exp(\beta[U_1(x)-U_1(x-s)])$ at $s=0$ and
by taking into account (\ref{Uprime1a}), we get for any $x\in\TT$
and $s\in(0,\pi/(2\beta))$,
\begin{eqnarray*}
\bigl\llvert I_2(x,s)\bigr\rrvert&\leq& 2\frac{\llVert \nu\rrVert
_{\infty}}{2\pi} \exp
\bigl(\beta\bigl\llVert U'_1\bigr\rrVert
_{\infty}s\bigr)\beta\bigl\llVert U'_1\bigr
\rrVert_{\infty}s
\\
&\leq&\frac{\llVert \nu\rrVert _{\infty}}{\pi} \exp(\beta s)\beta s
\\
&\leq&2 \frac{1+\pi A}{\pi}\exp(\pi/2)\beta s.
\end{eqnarray*}
The term $I_3(x,s)$ is bounded in a similar manner, rather expanding at
the second order the previous mapping
and using (\ref{Usecond1}) to see that $\llVert U_1''\rrVert
_\infty\leq A$.
\end{pf}

We finish this subsection with the a variant of Lemma~\ref{Tstarpa0}.
%

%le3.7 #&#
\begin{lem}\label{Tstarpa1}
There exists a universal constant $k>0$,
such that
for any $s>0$ and $\beta\geq1$ with $\beta s\leq1$, we have,
for any $f\in\mathcal{C}^1(\TT)$,
\begin{eqnarray*}
\int_{B(y,\pi-s)} \bigl(T^*_{y,s}[\widetilde
g_y](x)-g_y(x)\bigr)^2 \mu_\beta
(dx)&\leq& ks^2\beta^2 \biggl(\int(\partial
f)^2 \,d\mu_{\beta}+\int f^2 \,d\mu_\beta
\biggr),
\end{eqnarray*}
where $T^*_{y,s}$
is the adjoint operator of $T_{y,s}$ in $\LL^2(\mu_\beta)$ and
where for any fixed $y\in\TT$,
\begin{eqnarray*}
\forall x\in\TT\setminus\bigl\{y'\bigr\},&\qquad& \cases{
g_y(x) \ffff  f(x)\dot{\gamma}(x,y,0),
\vspace*{3pt}\cr
\widetilde g_y(x)
\ffff \un_{(y-\pi,y-s)\sqcup(y+s,y+\pi)}(x)g_y(x).}
\end{eqnarray*}
\end{lem}

\begin{pf}
As remarked at the beginning of the proof of Lemma~\ref{Tstarpa0}, it
is sufficient to deal with the case $y=0$.
To simplify the notation, we remove $y=0$ from the indices, in
particular we consider the mappings $g$ and $\widetilde g$ defined by
$g(x)= -\sign(x)f(x)$ and $\widetilde g(x)= \un_{(-\pi,-s)\sqcup
(s,\pi)}(x)g(x)$.

Taking into account that $\widetilde g$ vanishes on $(-s,s)$, we deduce from
Lemmas~\ref{dagger1} and
\ref{Ldual1} that for a.e. $x\in(-\pi+s,\pi-s)$,
\begin{eqnarray*}
T^*_{s}[\widetilde g](x)&=&\exp\bigl(\beta U_2(x)
\bigr)T_{-s}\bigl[\exp(-\beta U_2)\widetilde g\bigr](x).
\end{eqnarray*}
This observation leads us to consider the upper bound
\begin{eqnarray*}
\int_{-\pi+s}^{\pi-s} \bigl(T^*_{s}[
\widetilde g](x)-g(x)\bigr)^2 \mu_\beta(dx) &
\leq&2J_1+ 2J_2,
\end{eqnarray*}
where
\begin{eqnarray*}
J_1&\ffff & \int_{-\pi+s}^{\pi-s} \bigl(\exp
\bigl(\beta\bigl[U_2(x)-U_2\bigl(x+\sign(x)s\bigr)\bigr]
\bigr)-1\bigr)^2\bigl(T_{-s}[\widetilde g]
\bigr)^2 \mu_\beta(dx),
\\
J_2&\ffff & \int_{-\pi+s}^{\pi-s}
\bigl(T_{-s}[\widetilde g] -g\bigr)^2 \,d\mu_\beta.
\end{eqnarray*}
The arguments used in the proof of Lemma~\ref{Tstarpa0} to deal with
$J_1$ and $J_2$ can now
be easily adapted (even simplified) to obtain the wanted bounds.
For instance, one would have noted that
\begin{eqnarray*}
J_2&=&\int_{-\pi+s}^0\bigl(g(x-s)-g(x)
\bigr)^2 \mu_\beta(dx)+ \int^{\pi-s}_0
\bigl(g(x+s)-g(x)\bigr)^2 \mu_\beta(dx).
\end{eqnarray*}\upqed
\end{pf}

%s3.3 #&#
\subsection{Estimate of \texorpdfstring{$L^*_{\alpha,\beta}[\un]$}{$L^*_{alpha,beta}[1]$} in the cases $1<p<2$}\label{sec3.3}\label{p2}

\begin{figure}

\includegraphics{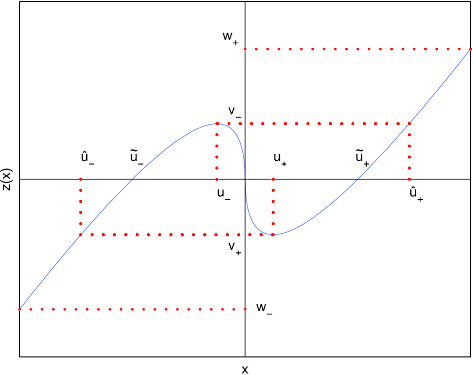}

\caption{The function $z$.}\label{fig3}
\end{figure}

In this situation, for any fixed $y\in\TT$ and $s\geq0$, the
definition (\ref{Tys}) must be replaced by
%
%
%e3.13 #&#
\begin{eqnarray}
\label{Tys2}\forall x\in\TT,\qquad T_{y,s}f(x)&\ffff & f\bigl(\gamma
\bigl(x,y,sd^{p-1}(x,y)\bigr)\bigr).
\end{eqnarray}
It leads us to introduce the function $z$ defined on $(y-\pi,y+\pi)$ by
%
%
%e3.14 #&#
\begin{eqnarray}
\label{z} z(x)&\ffff &\cases{ x-s(x-y)^{p-1}, &\quad if $x\in
[y,y+\pi)$,
\vspace*{3pt}\cr
x+s(y-x)^{p-1}, &\quad if $x\in(y-\pi,y]$.}
\end{eqnarray}
To study the variations of this function, by symmetry, it is sufficient
to consider its restriction to $(y,y+\pi)$. We need the following definitions,
all of them depending on $y\in\TT$, $s\geq0$ and $p\in(1,2)$:
\begin{eqnarray*}
u_+&\ffff & y+(p-1)^{\afrac{1}{2-p}}s^{\afrac{1}{2-p}},
\\
\widetilde u_+&\ffff & y+s^{\afrac{1}{2-p}},
\\
v_+&\ffff & y- \bigl((p-1)^{\vafrac{p-1}{2-p}}-(p-1)^{\afrac{1}{2-p}} \bigr)
s^{\afrac{1}{2-p}},
\\
w_+&\ffff & y+\pi-\pi^{p-1}s.
\end{eqnarray*}
Let $\sigma(p)$ be the largest positive real number in $(0,1/2)$ such that
for $s\in(0,\sigma(p))$, we have $u_+<y+\pi$, $v_+>y-\pi$ and $w_+-y>y-v_+$.
One checks that for $s\in(0,\sigma(p))$, the function $z$ is
decreasing on $(y,u_+)$ and increasing on $(u_+,y+\pi)$.
Furthermore $v_+=z(u_+)$, $w_+=z(y+\pi)$\vadjust{\goodbreak} and $\widetilde u_+$ is the unique
point in $(u_+,y+\pi)$ such that $z(\widetilde u_+)=y$.
Let us also introduce $\widehat u_+$ the unique point in $(\widetilde
u_+,y+\pi)$
such that
and $z(\widehat u_+)=-v_+$.
All these definitions, as well as the symmetric notions with respect to
$(y,y)$, where the indices $+$ are replaced by $-$, are summarized in
the following picture (see Figure \ref{fig3}).

%
%f3 #&#

Thus for $s\in(0,\sigma(p))$, we can consider $\varphi_+ \st
[v_+,y]\rightarrow[y,u_+]$ and
$\psi_+\st[v_+,w_+]\rightarrow[u_+,y+\pi]$ the inverses of $z$,
respectively, restricted to $[y,u_+]$ and $[u_+,y+\pi]$.
The mappings $\varphi_-$ and $\psi_-$ are defined in a symmetrical
manner on $[y,v_-]$ and $[w_-,v_-]$.
These quantities were necessary to compute the adjoint $T^{\dagger
}_{y,s}$ of $T_{y,s}$ in $\LL^2(\lambda)$,
for any fixed $y\in\TT$ and $s>0$ small enough.

%
%
%le3.8 #&#
\begin{lem}\label{dagger2}
Assume that $s\in(0,\sigma(p))$.
Then for any bounded and measurable function $g$, we have, for almost
every $x\in\TT$ (identified with its representative in $(y-\pi,y+\pi)$),
\begin{eqnarray*}
T^{\dagger}_{y,s}[g](x) &=&\un_{(w_-,v_+)}(x)
\psi_-'(x)g\bigl(\psi_-(x)\bigr)+\un_{(v_-,w_+)}(x)\psi
_+'(x)g\bigl(\psi_+(x)\bigr)
\\
&&{} +\un_{(v_+,y)}(x)\bigl[\psi_-'(x)g\bigl(\psi_-(x)\bigr)+
\psi_+'(x)g\bigl(\psi_+(x)\bigr)+\bigl\llvert\varphi_+'(x)
\bigr\rrvert g\bigl(\varphi_+(x)\bigr)\bigr]
\\
&&{} +\un_{(y,v_-)}(x)\bigl[\psi_-'(x)g\bigl(\psi_-(x)\bigr)+
\psi_+'(x)g\bigl(\psi_+(x)\bigr)+\bigl\llvert\varphi'_-(x)
\bigr\rrvert g\bigl(\varphi_-(x)\bigr)\bigr].
\end{eqnarray*}
\end{lem}

\begin{pf}
The above formula is based on straightforward applications of the
change of variable formula.
For instance one can write for any bounded and measurable functions
$f,g$ defined on~$(y-\pi,y+\pi)$,
\begin{eqnarray*}
\int_{(y,u_+)} g(x)f\bigl(T_{y,s}(x)\bigr) \,dx&=&\int
_{(v_+,y)} f(z) g\bigl(\varphi_+(z)\bigr) \bigl\llvert
\varphi_+'(z)\bigr\rrvert \,dz.
\end{eqnarray*}\upqed
\end{pf}

Since we are more interested in adjoint operators in $\LL^2(\mu_\beta
)$, let us define for any fixed $y\in\TT$, $s\in(0,\sigma(p))$ and
any bounded and measurable function $f$
defined on $(y-\pi,y+\pi)$,
%
%
%e3.15 #&#
\begin{eqnarray}
\label{Tysstar} T_{y,s}^{*}[f]&\ffff & \exp(\beta
U_p) T_{y,s}^{\dagger}\bigl[\exp(-\beta
U_p)f \bigr].
\end{eqnarray}
Then we get the equivalent of Lemmas~\ref{Ldual0} and~\ref{Ldual1}.

%
%
%le3.9 #&#
\begin{lem}\label{Ldual2}
For $\alpha>0$ and $\beta> 0$ such that $s\ffff  p\alpha\beta/2\in
(0,\sigma(p))$, the domain of the
maximal extension of $L_{\alpha,\beta}$ on $\LL^2(\mu_\beta)$ is
$\cD$.
Furthermore, the domain of its dual operator $L_{\alpha,\beta}^*$ in
$\LL^2(\mu_\beta)$
is $\cD^*$ and we have for any $f\in\cD^*$,
\begin{eqnarray*}
L_{\alpha,\beta}^*f &=& \frac{1}2 \exp(\beta U_p)
\partial^2\bigl[\exp(-\beta U_p)f\bigr]+
\frac{1}\alpha\int T_{y,s}^{*}[f] \nu(dy)-
\frac{f}{\alpha}.
\end{eqnarray*}
In particular, if $\nu$ admits a continuous density, then $\cD^*=\cD
$ and the above formula holds for any $f\in\cD$.
\end{lem}

Once again, the assumption that $\nu$ admits a continuous density
enables us to consider
$L_{\alpha,\beta}^*\un$, which is given, under the conditions of the
previous lemma, for almost every $x\in\TT$, by
%
%
%e3.16 #&#
\begin{eqnarray}
\label{Lstarp} L_{\alpha,\beta}^*\un(x)&=&\frac{\beta^2}{2}
\bigl(U_p'(x)\bigr)^2-\frac
{\beta}{2}
U_p''(x)+\frac{1}{\alpha} \biggl(\int
T_{y,\vfrac{p\alpha\beta}{2}}^{*}[\un](x) \nu(dy)-1 \biggr).
\end{eqnarray}
We deduce the following.

%
%
%pr3.3 #&#
\begin{pro}\label{mbl2}
Assume that $\nu$ admits a density with respect to $\lambda$
satisfying (\ref{aA}).
Then there exists a constant $C(A,p)>0$, only depending on $A>0$ and
$p\in(1,2)$, such that
for any $\beta\geq1$ and $\alpha\in(0,\sigma(p)/\beta^{2})$, we have
\begin{eqnarray*}
\bigl\llVert L_{\alpha,\beta}^*\un\bigr\rrVert_\infty&\leq&C(A,p)\max
\bigl(\alpha\beta^4,\alpha^{p-1}\beta^{1+p},
\alpha^a\beta^{1+a} \bigr).
\end{eqnarray*}
\end{pro}

\begin{pf}
We first keep in mind that from (\ref{Uprime}) and Lemma~\ref{regU},
we have for all $x\in\TT$,
%
%
%e3.17 #&#
%e3.18 #&#
\begin{eqnarray}
\label{Uprimep} U_p'(x)&=&p \biggl(\int
_{x-\pi}^x(x-y)^{p-1} \nu(dy)-\int
_{x}^{x+\pi
} (y-x)^{p-1} \nu(dy) \biggr),
\\
\label{Usecondp} U_p''(x)&=&p(p-1)\int
_\TT d^{p-2}(y,x) \nu(dy)-p\pi^{p-2} \nu
\bigl(x'\bigr).
\end{eqnarray}
Taking into account (\ref{Lstarp}), our goal is to see how the terms
$\beta(U_p'(x))^2$ and $-U_p''(x)$ cancel
with some parts of the
integral
\begin{eqnarray*}
\frac{p}{s}\int T_{y,s}^{*}[\un](x)-1 \nu(dy),
\end{eqnarray*}
where $s\ffff  p\alpha\beta/2\in(0,\sigma(p)/\beta)\subset(0,\sigma(p))$,
and to bound what remains
by a quantity of the form $C'(A,p)( \beta^2 s+\beta s^{p-1}+s^a)$, for
another constant
$C'(A,p)>0$, only depending on $A>0$ and $p\in(1,2)$.

We decompose the domain of integration of $\nu(dy)$ into six essential
parts (with the convention that $-\pi\leq y-x<\pi$ and
remember that the points $w_-, v_+, v_-$ and $w_+$ depend on $y$):
\begin{eqnarray*}
J_1&\ffff & \{y\in\TT\st y-\pi<x<w_-\},
\\
J_2&\ffff & \{y\in\TT\st w_-<x<v_+\},
\\
J_3&\ffff & \{y\in\TT\st v_+<x<y\},
\\
J_4&\ffff & \{y\in\TT\st y<x<v_-\},
\\
J_5&\ffff & \{y\in\TT\st v_-<x<w_+\},
\\
J_6&\ffff & \{y\in\TT\st w_+<x<y+\pi\}.
\end{eqnarray*}

The cases of $J_1$ and $J_6$ are the simplest to treat. For instance,
for $J_6$, we write
that
\begin{eqnarray*}
\frac{p}{s}\int_{J_6} T_{y,s}^{*}[
\un](x)-1 \nu(dy)&=&- \frac
{p}{s}\int_{x'}^{x'+\pi^{p-1}s}
\un\nu(dy)
\\
&=& -\frac{p}{s}\int_{x'}^{x'+\pi^{p-1}s} \nu(y)
\frac{dy}{2\pi
}
\\
&=& - \frac{p\pi^{p-2}}{2}\nu\bigl(x'\bigr)- \frac{p}{2\pi s}\int
_{x'}^{x'+\pi^{p-1}s} \nu(y) -\nu\bigl(x'\bigr)
\,dy.
\end{eqnarray*}
A similar computation for $J_1$ and the use of assumption (\ref{aA})
lead to the bound
\begin{eqnarray}
\label{J6}
\nonumber
\biggl\llvert\frac{p}{s}\int_{J_1\sqcup J_6}
T_{y,s}^{*}[\un](x)-1 \nu(dy)+p\pi^{p-2}\nu
\bigl(x'\bigr)\biggr\rrvert&\leq& Ap\frac{\pi
^{(1+a)(p-1)-1}}{1+a}s^a
\nonumber\\[-8pt]\\[-8pt]\nonumber
&\leq&2\pi A s^{a}.
\end{eqnarray}

The most important parts correspond to $J_2$ and $J_5$. For example,
considering $J_5$, which can be written down as the segment
$(x_-,x_+)$, with
\begin{eqnarray*}
x_-&\ffff &x-\pi+\pi^{p-1}s,
\\
x_+&\ffff &x- \bigl((p-1)^{\vafrac{p-1}{2-p}}-(p-1)^{\afrac{1}{2-p}} \bigr)
s^{\afrac{1}{2-p}},
\end{eqnarray*}
we have to evaluate the
integral
%
%
%e3.19 #&#
\begin{eqnarray}
\label{J51} \frac{p}{s}\int_{x_-}^{x_+}
\psi_+'(x)\exp\bigl(\beta\bigl[U_p(x)-U_p
\bigl(\psi_+(x)\bigr)\bigr]\bigr)-1 \nu(dy)
\end{eqnarray}
($y$ is present in the integrand through $\psi_+(x)$ and $\psi
_+'(x)$). Indeed, in view of (\ref{Uprimep}) and (\ref{Usecondp}),
we would like to compare it to
%
%
%e3.20 #&#
\begin{eqnarray}
\label{J52} -\beta U_p'(x)\int_{x_-}^{x_+}
(x-y)^{p-1} \nu(dy)+p(p-1)\int_{x_-}^{x_+}
(x-y)^{p-2} \nu(dy).
\end{eqnarray}
To do so, we will expand the terms $\psi_+'(x)$ and $\exp(\beta
[U_p(x)-U_p(\psi_+(x))])$ as functions of the (hidden) parameter $s>0$.
Fix $y\in J_5$ and recall that it amounts to
$x\in(v_-,w_+)$.
Due to (\ref{z}) and to the definition of $\psi_+$, we have for such $x$,
%
%
%e3.21 #&#
\begin{eqnarray}
\label{psiprime} \psi_+'(x)&=&\frac{1}{1-s(p-1)(\psi_+(x)-y)^{p-2}}.
\end{eqnarray}

Let us begin by working heuristically, to outline why the quantities
(\ref{J51}) and (\ref{J52}) should be close.
From the above expression, we get
\begin{eqnarray*}
\psi_+'(x)&\simeq&1+s(p-1) \bigl(\psi_+(x)-y\bigr)^{p-2}.
\end{eqnarray*}
By definition of $\psi_+$, we have
\begin{eqnarray}
\label{xypsi}
\nonumber
x-y&=&\psi_+(x)-y-s\bigl(\psi_+(x)-y
\bigr)^{p-1}
\nonumber\\[-8pt]\\[-8pt]\nonumber
&=&\bigl(\psi_+(x)-y\bigr) \bigl(1-s\bigl(\psi_+(x)-y\bigr)^{p-2}
\bigr),
\end{eqnarray}
so that $x-y\simeq\psi_+(x)-y$ and
\begin{eqnarray*}
\psi_+'(x)&\simeq&1+s(p-1) (x-y)^{p-2}.
\end{eqnarray*}
On the other hand,
\begin{eqnarray*}
\exp\bigl(\beta\bigl[U_p(x)-U_p\bigl(\psi_+(x)\bigr)
\bigr]\bigr)&\simeq&1+\beta\bigl[U_p(x)-U_p\bigl(
\psi_+(x)\bigr)\bigr]
\\
&\simeq&1+\beta U_p'(x) \bigl(x-\psi_+(x)\bigr)
\\
&=&1-s\beta U_p'(x) \bigl(\psi_+(x)-y
\bigr)^{p-1}
\\
&\simeq& 1-s\beta U_p'(x) (x-y)^{p-1}.
\end{eqnarray*}
Putting together these approximations, we end up with
\begin{eqnarray*}
\psi_+'(x)\exp\bigl(\beta\bigl[U_p(x)-U_p
\bigl(\psi_+(x)\bigr)\bigr]\bigr)-1&\simeq& s\bigl[(p-1) (x-y)^{p-2}-
\beta U_p'(x) (x-y)^{p-1}\bigr],
\end{eqnarray*}
suggesting the proximity of (\ref{J51}) and (\ref{J52}), after
integration with respect to $\nu(dy)$ on $(x_-,x_+)$.

To\vspace*{1pt} justify and quantify these computations, we start by remarking that
$\psi_+(x)-y$ is bounded below by $\widehat u_+-y$, itself bounded below
by $\widetilde u_+-y=s^{\afrac{1}{2-p}}$.
But this lower bound will not be sufficient in (\ref{xypsi}),
so let us
improve it a little. By definition of $\widehat u_+$, we have
\begin{eqnarray*}
v_--y&=&\widehat u-y-s(\widehat u-y)^{p-1},
\end{eqnarray*}
so that $\widehat u_+-y=k_ps^{\afrac{1}{2-p}}$ where $k_p$ is the unique
solution larger than 1 of the equation
%
%
%e3.22 #&#
\begin{eqnarray}
\label{kp} k_p-k_p^{p-1}&=&(p-1)^{\vafrac{p-1}{2-p}}-(p-1)^{\afrac{1}{2-p}}.
\end{eqnarray}
It follows that for any $y\in J_5$,
%
%
%e3.23 #&#
\begin{eqnarray}
\label{borneJ5}
\nonumber
1 &\leq&\frac{1}{1-s(\psi_+(x)-y)^{p-2}}\leq\frac{1}{1-s(\widehat u_+-y)^{p-2}}
\\
&=&\frac{\widehat u_+-y}{v_--y}
\\
&=&K_p,\nonumber
\end{eqnarray}
where the latter quantity only depends on $p\in(1,2)$ and is given by
\begin{eqnarray*}
K_p&\ffff & \frac{k_p}{(p-1)^{\vafrac{p-1}{2-p}}-(p-1)^{\afrac{1}{2-p}}}.
\end{eqnarray*}

In particular,
coming back to (\ref{psiprime}) and taking into account (\ref
{xypsi}), we get that for $y\in J_5'$,
\begin{eqnarray*}
\bigl\llvert\psi_+'(x)-1-s(p-1) \bigl(\psi_+(x)-y
\bigr)^{p-2}\bigr\rrvert&= & \frac{(s(p-1)(\psi
_+(x)-y)^{p-2})^2}{1-s(p-1)(\psi_+(x)-y)^{p-2}}
\\
&\leq&(p-1)^2s^2\frac{(\psi_+(x)-y)^{2(p-2)}}{1-s(\psi
_+(x)-y)^{p-2}}
\\
&=&(p-1)^2s^2\frac{(x-y)^{2(p-2)}}{(1-s(\psi
_+(x)-y)^{p-2})^{1+2(p-2)}}
\\
&\leq& (p-1)^2K_p^{(2p-3)_+} s^2(x-y)^{2(p-2)}.
\end{eqnarray*}
To complete this estimate, we note that in a similar way, still for
$y\in J_5$,
\begin{eqnarray*}
\bigl\llvert\bigl(\psi_+(x)-y\bigr)^{p-2}-(x-y)^{p-2}\bigr
\rrvert&=&(x-y)^{p-2}\bigl\llvert1-\bigl(1-s\bigl(\psi_+(x)-y
\bigr)^{p-2}\bigr)^{2-p}\bigr\rrvert
\\
&\leq&(x-y)^{p-2}\bigl\llvert1-\bigl(1-s\bigl(\psi_+(x)-y
\bigr)^{p-2}\bigr)\bigr\rrvert
\\
&=&s(x-y)^{p-2}\bigl(\psi_+(x)-y\bigr)^{p-2}
\\
&= &s(x-y)^{2(p-2)}\bigl(1-s\bigl(\psi_+(x)-y\bigr)^{p-2}
\bigr)^{2-p}
\\
&\leq&s(x-y)^{2(p-2)},
\end{eqnarray*}
so that in the end,
%
%
%e3.24 #&#
\begin{eqnarray}
\label{fact1}
&& \bigl\llvert\psi_+'(x)-1-s(p-1)
(x-y)^{p-2}\bigr\rrvert
\leq\bigl[(p-1)^2K_p^{(2p-3)_+}
+p-1\bigr]s^2(x-y)^{2(p-2)}.\qquad\quad
\end{eqnarray}

We now come to the term $\exp(\beta[U_p(x)-U_p(\psi_+(x))])$.
First we remark that
\begin{eqnarray*}
\bigl\llvert U_p(x)-U_p\bigl(\psi_+(x)\bigr)\bigr\rrvert
&\leq& \bigl\llVert U'_p\bigr\rrVert_\infty
\bigl\llvert x-\psi_+(x)\bigr\rrvert
\\
&\leq& p\pi^{p-1}s\bigl(\psi_+(x)-y\bigr)^{p-1}
\\
&\leq& p\pi^{2(p-1)}s
\\
&\leq&2\pi^2 s.
\end{eqnarray*}
It follows, recalling our assumption $\beta s\leq\sigma(p)$, that
\begin{eqnarray*}
&& \bigl\llvert\exp\bigl(\beta\bigl[U_p(x)-U_p
\bigl(\psi_+(x)\bigr)\bigr]\bigr)-1-\beta\bigl[U_p(x)-U_p
\bigl(\psi_+(x)\bigr)\bigr]\bigr\rrvert
\\
&&\quad \leq\frac{ \beta^2 [U_p(x)-U_p(\psi_+(x))]^2}{2}\exp\bigl(2\pi^2
\beta s\bigr)
\\
&&\quad \leq 2\pi^4\beta^2\exp\bigl(2\pi^2
\sigma(p)\bigr)s^2.
\end{eqnarray*}
In addition, we have
\begin{eqnarray*}
\bigl\llvert U_p(x)-U_p\bigl(\psi_+(x)
\bigr)-U_p'(x) \bigl(x-\psi_+(x)\bigr)\bigr\rrvert&
\leq& \frac{\llVert U''_p\rrVert _\infty}{2}\bigl(x-\psi_+(x)\bigr)^2.
\end{eqnarray*}
In view of (\ref{Usecondp}) and taking into account that $\int U_p''
\,d\lambda=0$, we have
\begin{eqnarray*}
\bigl\llVert U''_p\bigr\rrVert
_\infty&\leq& 2p(p-1)\llVert\nu\rrVert_\infty\int
_0^{\pi} u^{p-2} \frac{du}{2\pi}
\\
&=&2p\pi^{p-1}(1+\pi A).
\end{eqnarray*}
So we get
\begin{eqnarray*}
\bigl\llvert U_p(x)-U_p\bigl(\psi_+(x)
\bigr)-U_p'(x) \bigl(x-\psi_+(x)\bigr)\bigr\rrvert&
\leq& 2\pi(1+\pi A) \bigl(x-\psi_+(x)\bigr)^2
\\
&\leq& 2\pi(1+\pi A) s^2\bigl(\psi_+(x)-y\bigr)^{2(p-1)}
\\
&\leq&2\pi^3 (1+\pi A)s^2,
\end{eqnarray*}
namely
\begin{eqnarray*}
\bigl\llvert U_p(x)-U_p\bigl(\psi_+(x)
\bigr)+sU_p'(x) \bigl(\psi_+(x)-y\bigr)^{p-1}
\bigr\rrvert&\leq&2\pi^3 s^2.
\end{eqnarray*}
Finally, using the inequality
\begin{eqnarray*}
\forall u,v\geq0, \forall p\in(1,2),\qquad\bigl\llvert u^{p-1}-v^{p-1}
\bigr\rrvert&\leq& \llvert u-v\rrvert^{p-1},
\end{eqnarray*}
it appears that
%
%
%e3.25 #&#
\begin{eqnarray}
\label{pmoins1}
\nonumber
\bigl\llvert\bigl(\psi_+(x)-y\bigr)^{p-1}-(x-y)^{p-1}
\bigr\rrvert&\leq& \bigl\llvert\psi_+(x)-x\bigr\rrvert^{p-1}
\\
&= &\bigl\llvert\psi_+(x)-y\bigr\rrvert^{(p-1)^2}s^{p-1}
\\
&\leq&\pi^{(p-1)^2}s^{p-1},\nonumber
\end{eqnarray}
so we can deduce that
\begin{eqnarray*}
&& \bigl\llvert\exp\bigl(\beta\bigl[U_p(x)-U_p
\bigl(\psi_+(x)\bigr)\bigr]\bigr)-1+\beta sU_p'(x)
(x-y)^{p-1}\bigr\rrvert
\\
&&\quad \leq p\pi^p K_p\beta s^{p}+2
\pi^3\beta\bigl(1+\pi A+\pi\exp\bigl(2\pi^2 \sigma(p)
\bigr)\beta\bigr) s^2.
\end{eqnarray*}

From the latter bound and (\ref{fact1}), we obtain a constant
$K(p,A)>0$ depending only on $p\in(1,2)$ and $A>0$,
such that
\begin{eqnarray}
\label{fact2}
\nonumber
&& \frac{p}{s}\biggl\llvert\int
_{x_-}^{x_+} \psi_+'(x)\exp\bigl(\beta
\bigl[U_p(x)-U_p\bigl(\psi_+(x)\bigr)\bigr]\bigr) \nonumber
\\
&&\qquad{} -
\biggl(1+\frac{s(p-1)}{(x-y)^{2-p}}\biggr) \bigl(1-\beta sU_p'(x)
(x-y)^{p-1}\bigr) \nu(dy)\biggr\rrvert
\\
&&\quad \leq K(p,A) \biggl(\beta s^{p-1}+\beta^2 s+ s\int
_{x_-}^{x_+} (x-y)^{2(p-2)} \nu(dy) \biggr).\nonumber
\end{eqnarray}
This leads us to upper bound
\begin{eqnarray*}
\int_{x_-}^{x_+} (x-y)^{2(p-2)} \nu(dy)&\leq&
\frac{\llVert \nu
\rrVert
_\infty}{2\pi}\int_{x_-}^{x_+}
(x-y)^{2(p-2)} \,dy
\\
&\leq&\frac{1+A\pi}{2\pi}\int_{\kappa_ps^{\afrac{1}{2-p}}}^{\pi
-\pi^{p-1}s}
y^{2(p-2)} \,dy,
\end{eqnarray*}
with
%
%
%e3.26 #&#
\begin{eqnarray}
\label{kappap}\kappa_p&\ffff & (p-1)^{\vafrac{p-1}{2-p}}-(p-1)^{\afrac{1}{2-p}}.
\end{eqnarray}
An immediate computation gives, for $p\in(1,2)$, a constant $\kappa
_p'>0$ such that for any $s\in(0,\sigma(p))$,
%
%
%e3.27 #&#
\begin{eqnarray}
\label{int3cas} \int_{\kappa_ps^{\afrac{1}{2-p}}}^{\pi-\pi^{p-1}s}
y^{2(p-2)} \,dy
&\leq& \kappa_p' \cases{
1, &\quad if $p> 3/2$,
\vspace*{3pt}\cr
\displaystyle \ln
\bigl(\bigl(1+\sigma(p)\bigr)/s\bigr), &\quad if $p= 3/2$,
\vspace*{3pt}\cr
\displaystyle s^{\vafrac{2p-3}{2-p}}, &
\quad if $p< 3/2$.}
\end{eqnarray}
Since $1+\frac{2p-3}{2-p}>p-1$, $\beta\geq1$ and $s\in(0,\sigma(p))$,
we can find another constant $K'(p,A)>0$ such that
the right-hand side of (\ref{fact2}) can be replaced by
$K'(p,A)(\beta s^{p-1}+\beta^2 s)$.
It is now easy to see that such an expression, up to a new change of
the factor $K'(p,A)$,
bounds the difference between (\ref{J51}) and (\ref{J52}). Indeed,
just use that
\begin{eqnarray*}
\int_{\kappa_ps^{\afrac{1}{2-p}}}^{\pi-\pi^{p-1}s} y^{2p-3} \,dy&\leq&\pi
\int
_{\kappa_ps^{\afrac{1}{2-p}}}^{\pi-\pi^{p-1}s} y^{2(p-2)} \,dy
\end{eqnarray*}
and resort to (\ref{int3cas}).

There is no more difficulty in checking that the cost of replacing
$x_-$ and $x_+$, respectively, by
$x-\pi$ and $x$
in (\ref{J52}) is also bounded by $K''(p,A)(\beta
s^{p/(2-p)}+s^{(p-1)/(2-p)})\leq2K''(p,A)\beta s^{p-1}$, for an
appropriate choice
of the factor $K''(p,A)$ depending on $p\in(1,2)$ and $A>0$.

Symmetrical computations for $J_2$ and remembering (\ref{J6}) lead to the
existence of a constant $K'''(p,A)>0$, depending only on $p\in(1,2)$
and $A>0$, such that for
$\beta\geq1$ and $s\in(0,\sigma(p)/\beta)$, we have
\begin{eqnarray*}
&& \biggl\llvert\beta\bigl(U_p'(x)\bigr)^2-U_p''(x)+
\frac{p}{s} \biggl(\int_{J_1\sqcup
J_2\sqcup J_5\sqcup J_6} T_{y,s}^{*}[
\un](x) \nu(dy)-1 \biggr)\biggr\rrvert
\\
&&\quad \leq K'''(p,A)
\bigl(s^a +\beta s^{p-1}+\beta^2 s\bigr).
\end{eqnarray*}

It remains to treat the segments $J_3$ and $J_4$ and again by symmetry,
let us deal with $J_4$ only: it is sufficient to exhibit
a constant $K^{(4)}(p,A)>0$, depending on $p\in(1,2)$ and $A>0$, such
that for $\beta\geq1$ and $s\in(0,\sigma(p)/\beta)$,
\begin{eqnarray*}
\frac{p}{s}\biggl\llvert\int_{J_4}
T_{y,s}^{*}[\un](x)-1 \nu(dy)\biggr\rrvert&\leq&
K^{(4)}(p,A)s^{\vafrac{p-1}{2-p}}
\end{eqnarray*}
(since the right-hand side is itself bounded by $K^{(4)}(p,A)(\sigma
(p))^{\afrac{(p-1)^2}{2-p}}s^{p-1}$),
or equivalently
%
%
%e3.28 #&#
\begin{eqnarray}
\label{J4} \biggl\llvert\int_{J_4} T_{y,s}^{*}[
\un](x)-1 \nu(dy)\biggr\rrvert&\leq& \frac{K^{(4)}(p,A)}{p}s^{\afrac{1}{2-p}}.
\end{eqnarray}
The constant part is immediate to bound:
\begin{eqnarray*}
\int_{J_4} 1 \nu(dy)&\leq& \frac{\llVert \nu\rrVert
_\infty}{2\pi}\int
_{J_4}1 \,dy
\\
&\leq& \frac{1+\pi A}{2\pi}\int_{x-\kappa_ps^{1/(2-p)}}^{x}1 \,dy
\\
&=& \frac{(1+\pi A)\kappa_p}{2\pi}s^{\afrac{1}{2-p}}.
\end{eqnarray*}
For the other part, we first remark that for $y\in J_4$, we have
\begin{eqnarray*}
%\begin{array} {rcccl}
y&<&x<y+\kappa_ps^{\afrac{1}{2-p}},
\\
y+s^{\afrac{1}{2-p}}&<&\psi_+(x)<y+k_ps^{\afrac{1}{2-p}},
\\
y-(p-1)^{\afrac{1}{2-p}}s^{\afrac{1}{2-p}}&<&\varphi_-(x)<y,
\\
y-s^{\afrac{1}{2-p}}&<&\psi_-(x)<y-(p-1)^{\afrac{1}{2-p}}s^{\afrac{1}{2-p}}
%\end{array}
%
\end{eqnarray*}
(recall that $\widehat u_+=y+k_ps^{\afrac{1}{2-p}}$ with $k_p$ defined in
(\ref{kp})).
It follows that we can find a constant $\kappa_p''>0$, depending only
on $p\in(1,2)$,
such that for $s\in(0,\sigma(p))$,
\begin{eqnarray*}
&& \max\bigl(\bigl\llvert U_p(x)-U_p\bigl(\psi_+(x)\bigr)
\bigr\rrvert,\bigl\llvert U_p(x)-U_p\bigl(\psi_-(x)\bigr)
\bigr\rrvert,\bigl\llvert U_p(x)-U_p\bigl(\varphi_-(x)
\bigr)\bigr\rrvert\bigr)
\\
&&\quad \leq \kappa_p''
s^{\afrac{1}{2-p}}
\\
&&\quad \leq \kappa_p''\bigl( \sigma(p)
\bigr)^{\vafrac{p-1}{2-p}}s.
\end{eqnarray*}
In particular, we can find another constant $\kappa_p'''>0$, such that
under the conditions that $\beta\geq1$ and $\beta s\in(0,\sigma(p))$,
\begin{eqnarray*}
\exp\bigl(\beta\max\bigl(\bigl\llvert U_p(x)-U_p\bigl(
\psi_+(x)\bigr)\bigr\rrvert,\bigl\llvert U_p(x)-U_p\bigl(
\psi_-(x)\bigr)\bigr\rrvert,\bigl\llvert U_p(x)-U_p\bigl(
\varphi_-(x)\bigr)\bigr\rrvert\bigr) \bigr) &\leq& \kappa_p'''.
\end{eqnarray*}
Thus, denoting $\psi$ one of the functions $\psi_+$, $\varphi_-$ or
$\psi_-$, and remembering the bound
$\llVert \nu\rrVert _\infty\leq1+\pi A$, it is sufficient
to exhibit another
constant $\kappa_p^{(4)}>0$ such that
%
%
%e3.29 #&#
\begin{eqnarray}
\label{J41} \int_{J_4} \bigl\llvert\psi'(x)
\bigr\rrvert \,dy&\leq&\kappa_p^{(4)}s^{\afrac{1}{2-p}}.
\end{eqnarray}
Let us consider the case $\psi=\psi_+$, the other functions admit a
similar treatment.
We begin by making the dependence of $\psi_+(x)$ more explicit by
writing it $\psi_+(x,y)$.
From the definition of this quantity (see the first line of (\ref
{xypsi})) and from (\ref{psiprime}), we get
\begin{eqnarray*}
\partial_y\psi_+(x,y)&=&-\frac{s(p-1)(\psi
_+(x,y)-y)^{p-2}}{1-s(p-1)(\psi_+(x,y)-y)^{p-2}}
\\
&=&-s(p-1) \bigl(\psi_+(x,y)-y\bigr)^{p-2}\partial_x
\psi_+(x,y),
\end{eqnarray*}
so that the left-hand side of (\ref{J41}) can be rewritten
\begin{eqnarray*}
&& \frac{1}{s(p-1)}\int_{J_4} \bigl\llvert\bigl(
\psi_+(x,y)-y\bigr)^{2-p}\partial_y\psi_+(x,y)\bigr\rrvert
\,dy
\\
&&\quad \leq \frac{1}{s(p-1)}\int_{J_4} \bigl(
k_ps^{\afrac{1}{2-p}}\bigr)^{2-p}\bigl\llvert
\partial_y\psi_+(x,y)\bigr\rrvert \,dy
\\
&&\quad \leq \frac{ k_p^{2-p}}{(p-1)}\int_{J_4} \bigl\llvert\partial
_y\psi_+(x,y)\bigr\rrvert \,dy.
\end{eqnarray*}
Checking that $J_4=(x-\kappa_ps^{\afrac{1}{2-p}},x)$, the last integral
is equal to
$|\psi_+(x,x)-\psi_+(x,x- \kappa_ps^{\afrac{1}{2-p}})|$.
By definition of $\psi_+$,
we have $\psi_+(x,x)=x$ and it appears that the quantity $\zeta\ffff
\psi_+(x,x-\kappa_ps^{\afrac{1}{2-p}})-x$
is a positive solution to the equation
\begin{eqnarray*}
\zeta=s\bigl(\zeta+\kappa_p s^{\afrac{1}{2-p}}\bigr)^{p-1}.
\end{eqnarray*}
It follows that $\zeta=k_p's^{\afrac{1}{2-p}}$ where $k_p'$
is the unique positive solution of $k_p'=(k_p'+\kappa
_p)^{p-1}$.

Thus, (\ref{J41}) is proven and we can conclude to the validity of
(\ref{J4}).
\end{pf}

To finish this subsection, here is a version of Lemma~\ref{Tstarpa1}
for $p\in(1,2)$,
which is a little weaker, since we need a preliminary integration with
respect to $\nu(y)$.
%

%le3.10 #&#
\begin{lem}\label{Tstarpa2}
Under the assumption (\ref{aA}),
there exists a universal constant $k(p,A)>0$, depending only on $p\in
(1,2)$ and $A>0$, such that
for any $s>0$ and $\beta\geq1$ with $\beta s\leq\sigma(p)$, we have,
for any $f\in\mathcal{C}^1(\TT)$,
\begin{eqnarray}
\label{borned2}
\nonumber
&& \int_\TT\nu(dy)\int
_{B(y,\pi-\pi^{p-1}s)} \bigl(T^*_{y,s}[\widetilde
g_y](x)-g_y(x)\bigr)^2 \mu_\beta(dx)
\nonumber\\[-8pt]\\[-8pt]\nonumber
&&\quad \leq k(p,A) \bigl( s^{2(p-1)} +\beta^2s^2\bigr)
\biggl(\int(\partial f)^2 \,d\mu_{\beta}+\int f^2 \,d
\mu_\beta\biggr),
\end{eqnarray}
where $T^*_{y,s}$
is the adjoint operator of $T_{y,s}$ in $\LL^2(\mu_\beta)$ and
where for any fixed $y\in\TT$,
\begin{eqnarray*}
\forall x\in\TT\setminus\bigl\{y'\bigr\},&\qquad&\cases{
\displaystyle g_y(x) \ffff  f(x)d^{p-1}(x,y)\dot{\gamma}(x,y,0),
\vspace*{3pt}\cr
\displaystyle\widetilde g_y(x) \ffff \un_{(y-\pi,y-s^{\afrac{1}{2-p}})\sqcup(y+s^{\afrac{1}{2-p}},y+\pi)}(x)g_y(x).}
\end{eqnarray*}
\end{lem}

\begin{pf}
We begin by fixing $y\in\TT$ and by remembering the notation of the
proof of Proposition~\ref{mbl2} (see Figure~\ref{fig3}).
Due to fact that $\widetilde g_y$ vanishes on $(\widetilde
u_-,\widetilde u_+)=(y-
s^{\afrac{1}{2-p}},y+
s^{\afrac{1}{2-p}})$, we deduce from Lemma~\ref{dagger2}
and
(\ref{Tysstar}) that for a.e. $x\in(y-\pi+\pi^{p-1}s, y+\pi-\pi^{p-1}s)$,
\begin{eqnarray*}
T^*_{y,s}[\widetilde g_y](x)&=&\psi_\varepsilon'(x)
\exp\bigl(\beta\bigl[ U_p(x)-U_p\bigl(
\psi_\varepsilon(x)\bigr)\bigr]\bigr)\widetilde g_y\bigl(
\psi_\varepsilon(x)\bigr),
\end{eqnarray*}
where $\varepsilon\in\{-,+\}$ stands for the sign of $x-y$ with the
conventions of the proof of Proposition~\ref{mbl2}.
Thus, we are led to the decomposition
\begin{eqnarray*}
\int_{B(y,\pi-\pi^{p-1}s)} \bigl(T^*_{y,s}[\widetilde
g_y](x)-g_y(x)\bigr)^2 \mu
_\beta(dx)&\leq& 3J_1(y) +3J_2(y)+3J_3(y),
\end{eqnarray*}
where
\begin{eqnarray*}
J_1(y)&\ffff &\int_{B(y,\pi-\pi^{p-1}s)} \bigl(\exp\bigl(\beta
\bigl[ U_p(x)-U_p\bigl(\psi_\varepsilon(x)\bigr)\bigr]
\bigr)-1\bigr)^2 \bigl(\psi_\varepsilon'(x)
\widetilde g_y\bigl(\psi_\varepsilon(x)\bigr)\bigr)^2
\mu_\beta(dx),
\\
J_2(y)&\ffff &\int_{B(y,\pi-\pi^{p-1}s)} \bigl(
\psi_\varepsilon'(x)\bigr)^2\bigl(\widetilde
g_y\bigl(\psi_\varepsilon(x)\bigr)-g_y(x)
\bigr)^2 \mu_\beta(dx),
\\
J_3(y)&\ffff &\int_{B(y,\pi-\pi^{p-1}s)} \bigl(
\psi_\varepsilon'(x)-1\bigr)^2g_y^2(x)
\mu_\beta(dx).
\end{eqnarray*}

We begin by dealing with $J_1(y)$, or rather with just half of it, by
symmetry and to avoid the consideration of $\varepsilon$:
\begin{eqnarray*}
\int_{y}^{y+\pi-\pi^{p-1}s} \bigl(\exp\bigl(\beta\bigl[
U_p(x)-U_p\bigl(\psi_+(x)\bigr)\bigr]\bigr)-1
\bigr)^2\bigl(\psi_+'(x)\widetilde g_y\bigl(
\psi_+(x)\bigr)\bigr)^2 \mu_\beta(dx).
\end{eqnarray*}
Let us recall that $x=\psi_+(x)-s(\psi_+(x)-y)^{p-1}$ and that $\psi
_+(x)-y\geq s^{\afrac{1}{2-p}}$. From (\ref{psiprime}),
we deduce that for $x\in(y,y+\pi-\pi^{p-1}s)$, $1\leq\psi_+(x)\leq
1/(2-p)$. Thus, it is sufficient to bound
\begin{eqnarray*}
\int_{y}^{y+\pi-\pi^{p-1}s} \bigl(\exp\bigl(\beta\bigl[
U_p(x)-U_p\bigl(\psi_+(x)\bigr)\bigr]\bigr)-1
\bigr)^2\bigl(\widetilde g_y\bigl(\psi_+(x)\bigr)
\bigr)^2 \mu_\beta(dx).
\end{eqnarray*}
Furthermore, for $x\in(y,y+\pi-\pi^{p-1}s)$, we have
%
%
%e3.30 #&#
\begin{eqnarray}
\label{xpsix} \bigl\llvert x-\psi_+(x)\bigr\rrvert&\leq& s \pi^{p-1},
\end{eqnarray}
so under the assumption that $s\beta\in(0,1/2)$, we can bound $ (\exp
(\beta[ U_p(x)-U_p(\psi_+(x))])-1)^2$
by a term of the form $k\beta^2 s^2$ for a universal constant $k>0$.
It remains to use $\widetilde g_y^2(x)\leq\pi^2f^2(x)$ to get an
upper bound
going in the direction of (\ref{borned2}).

We now come to $J_2(y)$ and again only to half of it:
\begin{eqnarray*}
\int_{y}^{y+\pi-\pi^{p-1}s} \bigl(\psi_+'(x)
\bigr)^2\bigl(\widetilde g_y\bigl(\psi_+(x)
\bigr)-g_y(x)\bigr)^2 \mu_\beta(dx).
\end{eqnarray*}
Due to the upper bound on $\psi_+$ seen just above, it is sufficient
to deal with
\begin{eqnarray*}
\int_{y}^{y+\pi-\pi^{p-1}s} \bigl(\widetilde g_y
\bigl(\psi_+(x)\bigr)-g_y(x)\bigr)^2 \mu_\beta(dx).
\end{eqnarray*}
But for $x\in(y,y+\pi-\pi^{p-1}s)$, we have $\psi_+(x)\in(y+
s^{\afrac{1}{2-p}},y+\pi)$, so that $\widetilde g_y(\psi_+(x))=g_y(\psi
_+(x))$ and the above expression is equal to
\begin{eqnarray*}
\int_{y}^{y+\pi-\pi^{p-1}s} \bigl( g_y\bigl(
\psi_+(x)\bigr)-g_y(x)\bigr)^2 \mu_\beta(dx).
\end{eqnarray*}
Coming back to the definition of $g_y$, it appears that for $x\in
(y,y+\pi-\pi^{p-1}s)$, both $\psi_+(x)$ and $x$ belong to
the same hemicircle obtained by cutting $\TT$ at $y$ and $y'$, so
\begin{eqnarray*}
&& \bigl( g_y\bigl(\psi_+(x)\bigr)-g_y(x)
\bigr)^2
\\
&&\quad =\bigl(d^{p-1}\bigl(y,\psi_+(x)\bigr)f\bigl(\psi_+(x)
\bigr)-d^{p-1}(y,x)f(x)\bigr)^2
\\
&&\quad \leq 2d^{2(p-1)}\bigl(y,\psi_+(x)\bigr) \bigl(f\bigl(\psi_+(x)
\bigr)-f(x)\bigr)^2+ 2f^2(x) \bigl(d^{p-1}
\bigl(y,\psi_+(x)\bigr)-d^{p-1}(y,x)\bigr)^2
\\
&&\quad \leq 2\pi^{2(p-1)}\bigl(f\bigl(\psi_+(x)\bigr)-f(x)\bigr)^2+2
\pi^{2(p-1)^2}s^{2(p-1)}f^2(x),
\end{eqnarray*}
where we have used (\ref{pmoins1}) to majorize the last term.
From (\ref{xpsix}), we deduce that
\begin{eqnarray*}
\bigl(f\bigl(\psi_+(x)\bigr)-f(x)\bigr)^2&\leq& 2s
\pi^{p-1}\int_{x-s \pi^{p-1}}^{x+s
\pi^{p-1}}
\bigl(f'(z)\bigr)^2 \,dz.
\end{eqnarray*}
As usual, the assumption $0<s\beta\leq1/2$ enables to find a
universal constant $k>0$
such that for any $z\in({x-s \pi^{p-1}},{x+s \pi^{p-1}})$, we have
$\mu_\beta(x)\leq k\mu_\beta(z)$. From the above computations, it
follows there exists another
universal constant $k'>0$ such that for any $y\in\TT$,
\begin{eqnarray*}
J_2(y) &\leq& k' \biggl(s^{2(p-1)}\int
f^2 \,d\mu_\beta+s^2\int\bigl(f'
\bigr)^2 \,d\mu_\beta\biggr)
\\
&\leq& k's^{2(p-1)} \biggl(\int f^2 \,d
\mu_\beta+\int\bigl(f'\bigr)^2 \,d\mu
_\beta\biggr).
\end{eqnarray*}

Finally, we come to $J_3(y)$, which will need to be integrated with
respect to $\nu(dy)$.
From (\ref{psiprime}), we first get that
\begin{eqnarray*}
J_3(y)&=& \int_{B(y,\pi-\pi^{p-1}s)} \biggl(\frac{s(p-1)d^{p-2}(\psi
_{\varepsilon
}(x),y)}{1-s(p-1)d^{p-2}(\psi_{\varepsilon}(x),y)}
\biggr)^2g_y^2(x) \mu_\beta(dx)
\\
&\leq&\frac{(p-1)^2}{(2-p)^2} s^2\int_{B(y,\pi-\pi^{p-1}s)}
d^{2(p-2)}\bigl(\psi_{\varepsilon}(x),y\bigr)g_y^2(x)
\mu_\beta(dx)
\\
&\leq& \frac{\pi^{2(p-1)}(p-1)^2}{(2-p)^2}s^2 \int_{B(y,\pi-\pi^{p-1}s)}
d^{2(p-2)}\bigl(\psi_{\varepsilon}(x),y\bigr)f^2(x)
\mu_\beta(dx).
\end{eqnarray*}
Next, recalling that $\llVert \nu\rrVert _\infty\leq1+\pi
A$ and that $d(\psi
_{\varepsilon}(x),y)\geq s^{\afrac{1}{2-p}}$ for any
$x\in B(y,\pi-\pi^{p-1}s)$, it appears that
\begin{eqnarray*}
&& \int_\TT J_3(y) \nu(dy)
\\
&&\quad \leq \frac{1+\pi A}{2\pi} \frac{\pi^{2(p-1)}(p-1)^2}{(2-p)^2}s^2\int _\TT \,dy \int_{B(y,\pi-\pi^{p-1}s)} d^{2(p-2)}
\bigl(\psi_{\varepsilon}(x),y\bigr)f^2(x) \mu_\beta(dx)
\\
&&\quad \leq \frac{1+\pi A}{2\pi} \frac{\pi^{2(p-1)}(p-1)^2}{(2-p)^2}s^2\int
_\TT\mu_\beta(dx) f^2(x)
\\
&&\qquad{}\times \int
_\TT\un_{ \{d(\psi_{\varepsilon}(x),y)\geq s^{\afrac{1}{2-p}} \}}
d^{2(p-2)}\bigl(
\psi_{\varepsilon}(x),y\bigr) \,dy.
\end{eqnarray*}
But for any fixed $z\in\RR/(2\pi\ZZ)$, we compute that
\begin{eqnarray*}
\int_\TT\un_{ \{d(z,y)\geq s^{\afrac{1}{2-p}} \}} d^{2(p-2)}(z,y) \,dy&=&2
\int_{ s^{\afrac{1}{2-p}}}^{\pi}\frac
{1}{y^{2(2-p)}} \,dy
\\
&\leq&k_p'' \cases{ 1, &\quad if $p>
3/2$,
\vspace*{3pt}\cr
\ln(1/s), &\quad if $p= 3/2$,
\vspace*{3pt}\cr
s^{\vafrac{2p-3}{2-p}}, &\quad if $p< 3/2$,}
\end{eqnarray*}
for $s\in(0,1/2)$ and for an appropriate constant $k_p''>0$ depending
only on $p\in(1,2)$.
It is not difficult to check that as $s\rightarrow0_+$, we have
\begin{eqnarray*}
s^{2(p-1)}&\gg&\cases{ s^2, &\quad if $p> 3/2$,
\vspace*{3pt}\cr
s^2\ln(1/s), &\quad if $p= 3/2$,
\vspace*{3pt}\cr
s^2s^{\vafrac{2p-3}{2-p}},
&\quad if $p< 3/2$.}
\end{eqnarray*}
It follows that for any $p\in(1,2)$, we can find a constant
$k'(p,A)>0$, depending only on $p\in(1,2)$
and $A>0$, such that
\begin{eqnarray*}
\int_\TT J_3(y) \nu(dy)&\leq&k'(p,A)
s^{2(p-1)} \int_\TT f^2(x)
\mu_\beta(dx).
\end{eqnarray*}
This ends the proof of the estimate (\ref{borned2}).
\end{pf}

%s3.4 #&#
\subsection{Estimate of \texorpdfstring{$L^*_{\alpha,\beta}[\un]$}{$L^*_{alpha,beta}[1]$} in the cases $p>2$}\label{sec3.4}\label{p3}

This situation is simpler than the one treated in the previous subsection
and is similar to the case $p=2$, because for $y\in T$ fixed and $s\geq
0$ small enough, the mapping $z$ defined in (\ref{z})
is injective when $p>2$. Again for any fixed $y\in\TT$ and $s\geq
0$, the definition (\ref{Tys}) has to be replaced by (\ref{Tys2}),
namely,
%
%
%e3.31 #&#
\begin{eqnarray}
\label{Tys3}\forall x\in\TT,\qquad T_{y,s}f(x)&\ffff & f\bigl(z(x)
\bigr).
\end{eqnarray}
With the previous subsections in mind, the computations are quite
straightforward, so
we will just outline them.

The first task is to determine the adjoint $ T_{y,s}^\dagger$ of
$T_{y,s}$ in $\LL^2(\lambda)$.
An immediate change of variable gives that for any $s\in(0,\sigma)$,
for any bounded and measurable function $g$, we have, for almost every
$x\in\TT$ (identified with its representative in $(y-\pi,y+\pi)$),
\begin{eqnarray*}
T^{\dagger}_{y,s}[g](x) &=&\un_{(y,z(y))}(x)
\psi'(x)g\bigl(\psi(x)\bigr),
\end{eqnarray*}
where $\sigma\ffff \pi^{2-p}/(p-1)$ and $\psi\st(z(y-\pi),z(y+\pi
))\rightarrow(y-\pi,y+\pi)$ is the inverse mapping of $z$ (with the slight
abuses of notation:
$z(y-\pi)\ffff  x-\pi+\pi^{p-1}s$, $z(y+\pi)\ffff  x+\pi-\pi^{p-1}s$).
The adjoint $T_{y,s}^{*}$ of $T_{y,s}$ in $\LL^2(\mu_\beta)$ is
still given by (\ref{Tysstar}).
As in the previous subsections,
this operator is bounded in $\LL^2(\mu_\beta)$.
It follows, if $\nu$ admits a continuous density with respect to
$\lambda$ and at least for $\alpha>0$ and $\beta\geq0$ such that
$s\ffff (p/2)\alpha\beta\in[0,\sigma)$, that
the adjoint $L_{\alpha,\beta}^*$ of $L_{\alpha,\beta}$ in $\LL
^2(\mu_\beta)$ is defined on $\cD$.
In particular, we can consider $L_{\alpha,\beta}^*\un$, which is
given, for almost every $x\in\TT$, by
%
%
%e3.32 #&#
\begin{eqnarray}
\label{Lstarp2} L_{\alpha,\beta}^*\un(x)&=&\frac{\beta^2}{2}
\bigl(U_p'(x)\bigr)^2-\frac
{\beta}{2}
U_p''(x)+\frac{p\beta}{2s} \biggl(\int
T_{y,s}^{*}[\un](x) \nu(dy)-1 \biggr).
\end{eqnarray}
From this formula, we deduce the following.
%

%pr3.4 #&#
\begin{pro}\label{mbl3}
Assume that $\nu$ admits a density with respect to $\lambda$
satisfying (\ref{aA}).
Then there exists a constant $C(A,p)>0$, only depending on $A>0$ and
$p>2$, such that
for any $\beta\geq1$ and $\alpha\in(0,\sigma/(p\beta^{2}))$, we have
\begin{eqnarray*}
\bigl\llVert L_{\alpha,\beta}^*\un\bigr\rrVert_\infty&\leq&C(A,p)\max
\bigl(\alpha\beta^4,\alpha^a\beta^{1+a} \bigr).
\end{eqnarray*}
\end{pro}

\begin{pf}
The arguments are similar to those of the case $J_5$ in the proof of
Proposition~\ref{mbl2}, but are less involved,
because the omnipresent term $1-s(p-1)(\psi(x)-y)^{p-2}$ is now easy
to bound:
for any $s\in[0,\sigma/2]$, we have for any $y\in\TT$ and $x\in
(z(y-\pi),z(y+\pi))$,
\begin{eqnarray*}
\tfrac{1}2 \leq1-(p-1)\bigl\llvert\psi(x)-y\bigr\rrvert
^{p-2}s \leq1.
\end{eqnarray*}
In particular, we have under these conditions,
\begin{eqnarray*}
\psi'(x)&=&\frac{1}{1-(p-1)\llvert \psi(x)-y\rrvert
^{p-2}s}\in[1,2].
\end{eqnarray*}
Following the arguments of the previous subsection,
one finds a constant $K(p,A)$, depending only on $p>2$ and $A>0$,
such that for any $\beta\geq1$, $s\in[0,\sigma/(2\beta)]$ and
$x\in(z(y-\pi),z(y+\pi))$,
\begin{eqnarray*}
\bigl\llvert\psi_+'(x)-1-(p-1)\bigl\llvert\psi_+(x)-y\bigr
\rrvert^{p-2}s\bigr\rrvert&\leq& K(p,A)s^2,
\\
\bigl\llvert\exp\bigl(\beta\bigl[U_p(x)-U_p\bigl(
\psi_+(x)\bigr)\bigr]\bigr)-1+\beta\sign(x-y)U_p'(x)
\llvert x-y\rrvert^{p-1}s\bigr\rrvert&\leq& K(p,A)\beta^2s^2.
\end{eqnarray*}
This bound enables us to approximate $T^*_{\alpha,\beta}\un(x)-1$ up
to a term $\cO_{p,A}(\beta^2s^2)$ (recall that
this
designates a quantity which is bounded by an expression of the form
$K'(p,A)\beta^2s^2$ for a constant $K'(p,A)>0$
depending on $p>2$ and $A>0$), by
\begin{eqnarray*}
\bigl((p-1)\bigl\llvert\psi_+(x)-y\bigr\rrvert^{p-2}-\beta\sign
(x-y)U_p'(x)\llvert x-y\rrvert^{p-1} \bigr)s.
\end{eqnarray*}
Next, we consider
\begin{eqnarray}
\label{Jxprime}
\nonumber
J&\ffff & \bigl\{y\in\TT\st x\in\bigl(z(y-\pi),z(y+\pi)
\bigr)\bigr\}
\nonumber\\[-8pt]\\[-8pt]\nonumber
&=&\TT\setminus\bigl[x'-s\pi^{p-1},x'+s
\pi^{p-1}\bigr],
\end{eqnarray}
in order to decompose
\begin{eqnarray}
\label{decompJ}
\nonumber
&& \frac{p\beta}{2s}\int_\TT
T^*_{y,s}[\un](x)-1 \nu(dy)
\nonumber\\[-8pt]\\[-8pt]\nonumber
&&\quad =\frac{p\beta}{2s}\int_J T^*_{y,s}[\un](x)-1
\nu(dy)-\frac{p\beta}{2s}\nu\bigl(\bigl[x'-s\pi^{p-1},x'+s
\pi^{p-1}\bigr]\bigr).
\end{eqnarray}
According to the previous estimate, up to a term $\cO_{p,A}(\beta
^3s^2)$ the first integral is equal to
\begin{eqnarray*}
\frac{p(p-1)\beta}{2}\int_J d^{p-2}(y,x) \nu(dy)-
\frac{p\beta^2}{2} U_p'(x)\int_{J}
\sign(x-y)d^{p-1}(x,y) \nu(dy).
\end{eqnarray*}
In view of (\ref{Jxprime}), up to an additional term $\cO_{p,A}(\beta
^2s)$, we can replace $J$ in the above integrals by $\TT$.
Thus putting together (\ref{Lstarp2}) and (\ref{decompJ}) with (\ref
{Uprimep}) and (\ref{Usecondp}) (which are also valid here),
it remains to estimate
\begin{eqnarray*}
\frac{p\beta}{2}\biggl\llvert\pi^{p-2}\nu\bigl(x'
\bigr)-\frac{1}{s}\nu\bigl[x'-s\pi^{p-1},x'+s
\pi^{p-1}\bigr]\biggr\rrvert
\end{eqnarray*}
and this is easily done through the assumption (\ref{aA}).
\end{pf}

We finish this subsection with the equivalent of Lemma~\ref{Tstarpa0}.
%

%le3.11 #&#
\begin{lem}\label{Tstarpa3}
For $p>2$,
there exists a constant $k(p)>0$, depending only on $p>2$, such that
for any $s\in(0,\sigma)$, with $\sigma\ffff \pi^{2-p}/(p-1)$, and
$\beta\geq1$ with $\beta s\leq1$, we have,
for any $y\in\TT$ and $f\in\mathcal{C}^1(\TT)$,
\begin{eqnarray*}
\int_{B(y,\pi-s\pi^{p-1})} \bigl(T^*_{y,s}[g_y](x)-g_y(x)
\bigr)^2 \mu_\beta(dx)&\leq& k(p)s^2
\beta^2 \biggl(\int(\partial f)^2 \,d\mu_{\beta}+
\int f^2 \,d\mu_\beta\biggr),
\end{eqnarray*}
where $T^*_{y,s}$
is the adjoint operator of $T_{y,s}$ in $\LL^2(\mu_\beta)$ and
where for any fixed $y\in\TT$,
\begin{eqnarray*}
\forall x\in\TT\setminus\bigl\{y'\bigr\},\qquad
g_y(x)&\ffff & f(x)d^{p-1}(x,y)\dot{\gamma}(x,y,0).
\end{eqnarray*}
\end{lem}

\begin{pf}
We only sketch the arguments, which are just an adaptation of those of
the proof of Lemma~\ref{Tstarpa0}.
Again it is sufficient to deal with the case $y=0$, which is removed
from the notation, and consequently with the function $g(x)=-\sign
(x)\llvert x\rrvert ^{p-1}f(x)$.
As seen previously in this subsection, we have for $s\in(0,\sigma)$
and $x\in(-\pi,\pi)$,
\begin{eqnarray*}
T^*_s[g](x) &=&\un_{(-\pi+s\pi^{p-1}, \pi-s\pi^{p-1})}(x)\exp\bigl(\beta
\bigl[U_p(x)-U_p\bigl(\psi(x)\bigr)\bigr]\bigr)
\psi'(x)g\bigl(\psi(x)\bigr),
\end{eqnarray*}
where $\psi$ is the inverse mapping of
$(-\pi,\pi)\ni x\mapsto x-\sign(x)\llvert x\rrvert ^{p-1}$.
Recall that for $x\in(-\pi+s\pi^{p-1}, \pi-s\pi^{p-1})$,
%
%
%e3.33 #&#
\begin{eqnarray}
\label{fastoche} \psi'(x)&=&\frac{1}{1-(p-1)\llvert \psi(x)\rrvert
^{p-2}s}\in[1,2].
\end{eqnarray}
Considering the decomposition
\begin{eqnarray*}
T^*_s[g](x)-g(x)&=& \bigl(\exp\bigl(\beta\bigl[U_p(x)-U_p
\bigl(\psi(x)\bigr)\bigr]\bigr)-1\bigr)\psi'(x)g\bigl(\psi(x)\bigr)
\\
&&{}+\psi'(x) \bigl(g\bigl(\psi(x)\bigr)-g(x)\bigr)+ \bigl(
\psi'(x)-1\bigr)g(x),
\end{eqnarray*}
we are led, after integration with respect to $\un_{(-\pi+s\pi
^{p-1}, \pi-s\pi^{p-1})}(x) \mu_\beta(dx)$, to computations
similar to those of Sections~\ref{p0} and~\ref{p2},
and indeed simpler than in the latter one, due to the boundedness
property described in (\ref{fastoche}).
\end{pf}

Let us summarize the Propositions~\ref{mbl0},~\ref{mbl1},~\ref{mbl2}
and~\ref{mbl3} of the previous subsections into the statement.
%

%pr3.5 #&#
\begin{pro}\label{mbl}
Assume that (\ref{aA}) is satisfied and for $p\geq1$, consider the
constant $a(p)>0$ defined in (\ref{ap}).
Then there exists two constants $\sigma(p)\in(0,1/2)$ and $C(A,p)>0$,
depending only on the quantities inside the parentheses,
such that for any $\alpha>0$ and $\beta>1$ such that $\alpha\beta
<\sigma(p)$, we have
\begin{eqnarray*}
\sqrt{\mu_\beta\bigl[\bigl(L^*_{\alpha,\beta}\un\bigr)^2
\bigr]}&\leq&C(A,p)\alpha^{a(p)}\beta^4.
\end{eqnarray*}
\end{pro}

Despite this bound is very rough, since we have replaced an essential
norm by a $\LL^2$ norm,
it will be sufficient
in the next section, when $\alpha^{a(p)}\beta^4$
is small, as a measure of the discrepancy between $\mu_\beta$
and the invariant measure for $L_{\alpha,\beta}$.

%s4 #&#
\section{Proof of convergence}\label{sec4}\label{proof}

This is the main part of the paper: we are going to prove Theorem~\ref{t1}
by the investigation of the evolution of a $\LL^2$ type
functional.

On $\TT$ consider the algorithm $X\ffff (X_t)_{t\geq0}$ described in
the \hyperref[sec1]{Introduction}.
We require that the underlying probability measure $\nu$ admits a
density with respect to $\lambda$
which is H\"older continuous: $a\in(0,1]$ and $A>0$ are constants such that
(\ref{aA}) is satisfied.
For the time being,
the schemes $\alpha\st\RR_+\rightarrow\RR_+^*$ and
$\beta\st\RR_+\rightarrow\RR_+$ are assumed to be, respectively, continuous
and continuously differentiable.
Only later on, in Proposition~\ref{condab}, will we present the
conditions insuring the wanted convergence (\ref{mN1}).
On the initial distribution $m_0$, the last ingredient necessary to specify
the law of $X$, no hypothesis is made.
We also denote $m_t$ the law of $X_t$, for any $t>0$. From the lemmas
given in the \hyperref[append]{Appendix},
we have that $m_t$ admits a $\mathcal{C}^1$ density with respect to
$\lambda$,
which is equally written $m_t$.
As it was mentioned in the previous section, we want to compare these
temporal marginal laws
with the corresponding instantaneous Gibbs measures, which were defined
in (\ref{mub})
with respect to the potential $U_p$ given in (\ref{U}).
A convenient way to quantify this discrepancy is to consider the
variance of the density of
$m_t$ with respect to $\mu_{\beta_t}$ under the probability measure
$\mu_{\beta_t}$:
%
%
%e4.1 #&#
\begin{eqnarray}
\label{It} \forall t>0,\qquad I_t&\ffff & \int\biggl(
\frac{m_t}{\mu_{\beta_t}}-1 \biggr)^2 \,d\mu_{\beta
_t}.
\end{eqnarray}
Our goal here is to derive a differential inequality satisfied by this quantity,
which implies its convergence to zero under appropriate conditions on the
schemes $\alpha$ and $\beta$.
More precisely, our purpose is to obtain the following.
%

%pr4.1 #&#
\begin{pro}\label{Iprime}
There exists two constants $c_1(p,A), c_2(p,A)>0$, depending on $p\geq
1$ and $A>0$, and a constant $\varsigma(p)\in(0,1/2)$, depending on
$p\geq1$, such that
for any $t>0$ with $\beta_t\geq1$ and $0<\alpha_t\beta_t^2\leq
\varsigma(p)$, we have
\begin{eqnarray*}
I_t'&\leq& -c_1(p,A) \bigl(
\beta_t^{-3}\exp\bigl(-b(U_p)\beta_t
\bigr)-\alpha_t^{\widetilde a(p)}\beta_t^3-
\bigl\llvert\beta_t'\bigr\rrvert\bigr)I_t+c_2(p,A)
\bigl(\alpha_t^{a(p)}\beta_t^4+\bigl
\llvert\beta_t'\bigr\rrvert\bigr)\sqrt{I_t},
\end{eqnarray*}
where $b(U_p)$ was defined in (\ref{bU}), $a(p)$ in Proposition~\ref
{mbl} and
\begin{eqnarray*}
\widetilde a(p) &\ffff & \cases{ 1, &\quad if $p=1$ or $p\geq3/2$,
\vspace*{3pt}\cr
2(p-1), &
\quad if $p\in(1,3/2)$.}
\end{eqnarray*}
\end{pro}

At least formally, there is no difficulty to differentiate the quantity $I_t$
with respect to the time $t>0$. But we postpone the rigorous
justification of the following computations
to the end of the \hyperref[append]{Appendix}, where the regularity of the temporal
marginal laws is discussed in detail.
Thus, we get at any time $t>0$,
\begin{eqnarray*}
I'_t&=&2\int\biggl(\frac{m_t}{\mu_{\beta_t}}-1 \biggr)
\frac
{\partial_t
m_t}{\mu_{\beta_t}} \,d\mu_{\beta_t} -2\int\biggl(\frac{m_t}{\mu_{\beta
_t}}-1 \biggr)
\frac{m_t}{\mu
_{\beta
_t}}\partial_t \ln(\mu_{\beta_t}) \,d
\mu_{\beta_t}
\\
&&{}+\int\biggl(\frac
{m_t}{\mu_{\beta_t}}-1 \biggr)^2\partial_t
\ln(\mu_{\beta_t}) \,d\mu_{\beta_t}
\\
&=& 2\int\biggl(\frac{m_t}{\mu_{\beta_t}}-1 \biggr)\partial_t
m_t \,d\lambda-\int\biggl(\frac{m_t}{\mu_{\beta_t}}-1 \biggr)^2
\partial_t \ln(\mu_{\beta
_t}) \,d\mu_{\beta_t}
-2\int\biggl(\frac{m_t}{\mu_{\beta_t}}-1 \biggr)\partial_t \ln(\mu
_{\beta_t}) \,d\mu_{\beta_t}
\\
&\leq& 2\int\biggl(\frac{m_t}{\mu_{\beta_t}}-1 \biggr)\partial_t
m_t \,d\lambda+\bigl\llVert\partial_t \ln(
\mu_{\beta_t})\bigr\rrVert_{\infty
} \biggl( \int\biggl(
\frac{m_t}{\mu_{\beta_t}}-1 \biggr)^2 \,d\mu_{\beta_t}+ 2\int\biggl
\llvert\frac{m_t}{\mu_{\beta_t}}-1\biggr\rrvert \,d\mu_{\beta
_t} \biggr)
\\
&\leq& 2\int\biggl(\frac{m_t}{\mu_{\beta_t}}-1 \biggr)\partial_t
m_t \,d\lambda+\bigl\llVert\partial_t \ln(
\mu_{\beta_t})\bigr\rrVert_{\infty
} ( I_t+2
\sqrt{I_t} ),
\end{eqnarray*}
where we used the Cauchy--Schwarz inequality.
The last term is easy to deal with.
%

%le4.1 #&#
\begin{lem}\label{palnmu}
For any $t\geq0$, we have
\begin{eqnarray*}
\bigl\llVert\partial_t \ln(\mu_{\beta_t})\bigr\rrVert
_{\infty
}&\leq& \pi^p\bigl\llvert\beta_t'
\bigr\rrvert.
\end{eqnarray*}
\end{lem}

\begin{pf}
Since for any $t\geq0$, we have
\begin{eqnarray*}
\forall x\in\TT,\qquad\ln(\mu_{\beta_t})&=&- \beta_tU_p
(x)-\ln\biggl(\int\exp\bigl(-\beta_t U_p(y)\bigr)
\lambda(dy) \biggr),
\end{eqnarray*}
it appears that
\begin{eqnarray*}
\forall x\in\TT,\qquad\partial_t \ln(\mu_{\beta_t})&=&
\beta_t'\int U_p(y)-U_p(x)
\mu_{\beta
_t}(dy),
\end{eqnarray*}
so that
\begin{eqnarray*}
\bigl\llVert\partial_t \ln(\mu_{\beta_t})\bigr\rrVert
_{\infty
}&\leq& \osc(U_p)\bigl\llvert\beta_t'
\bigr\rrvert.
\end{eqnarray*}
The bound $\osc(U_p)\leq\pi^p$ is an immediate consequence of the
definition (\ref{U})
of $U_p$ and of the fact that the (intrinsic) diameter of $\TT$ is
$\pi$.
\end{pf}

Denote for any $t>0$, $f_t\ffff  m_t/\mu_{\beta_t}$.
If this function was to be $\mathcal{C}^2$, we would get, by the
martingale problem satisfied by the
law of $X$, that
\begin{eqnarray*}
\int\biggl(\frac{m_t}{\mu_{\beta_t}}-1 \biggr)\partial_t m_t \,d
\lambda&=& \int L_{\alpha_t,\beta_t} [ f_t-1 ] \,dm_t
\\
&=&\int L_{\alpha_t,\beta_t} [ f_t-1 ] f_t \,d\mu_{\beta_t},
\end{eqnarray*}
where $ L_{\alpha_t,\beta_t}$, described in the previous section, is
the instantaneous generator at time $t\geq0$ of~$X$. The interest of the estimate of Proposition~\ref{mbl} comes from
the decomposition of the previous term
into
\begin{eqnarray*}
&& \int L_{\alpha_t,\beta_t} [ f_t-1 ] (f_t-1) \,d\mu_{\beta_t}+\int L_{\alpha_t,\beta_t} [ f_t-1 ] \,d\mu_{\beta_t}
\\
&&\quad = \int L_{\alpha_t,\beta_t} [ f_t-1 ] (f_t-1) \,d\mu_{\beta_t}+\int( f_t-1)L_{\alpha_t,\beta_t}^*[\un] \,d\mu_{\beta_t}
\\
&&\quad \leq \int L_{\alpha_t,\beta_t} [ f_t-1 ] (f_t-1) \,d\mu_{\beta_t}+\sqrt{I_t}\sqrt{\mu_{\beta_t}\bigl[
\bigl(L_{\alpha
_t,\beta_t}^*[\un]\bigr)^2\bigr]}.
\end{eqnarray*}
It follows that to prove Proposition~\ref{Iprime}, it
remains to treat the first term in the above right-hand side. A first step is the
following.
%

%le4.2 #&#
\begin{lem}\label{Lffm}
There exist a constant $c_3(p,A)>0$, depending on $p\geq1$ and $A>0$
and a constant $\widetilde\sigma(p)\in(0,1/2)$, such that
for any $\alpha>0$ and $\beta\geq1$ such that $\alpha\beta^2\leq
\widetilde\sigma(p)$, we have,
for any $f\in\mathcal{C}^2(\TT)$,
\begin{eqnarray*}
&& \int L_{\alpha,\beta} [ f-1 ] (f-1) \,d\mu_{\beta}
\\
&&\quad \leq- \biggl(\frac{
1}2-c_3(p,A)\alpha^{\widetilde a(p)}
\beta^3 \biggr)\int(\partial f)^2 \,d\mu_{\beta}
+c_3(p,A)\alpha^{\widetilde a(p)}\beta^3
\int(f-1)^2 \,d\mu_{\beta},
\end{eqnarray*}
where $\widetilde a(p)$ is defined in Proposition~\ref{Iprime}.
\end{lem}

\begin{pf}
For any $\alpha>0$ and $\beta\geq0$, we begin by decomposing
the generator $L_{\alpha,\beta}$ into
%
%
%e4.2 #&#
\begin{eqnarray}
\label{LLR} L_{\alpha,\beta}&=&L_{\beta}+R_{\alpha,\beta},
\end{eqnarray}
where $L_\beta\ffff (\partial^2-\beta U_p'\partial)/2$ was defined in
(\ref{Lb})
(recall that $U_p'$ is well-defined, since $\nu$ has no atom)
and where $R_{\alpha,\beta}$ is the remaining operator.
An immediate integration by parts leads to
\begin{eqnarray*}
\int L_{\beta} [ f-1 ] (f-1) \,d\mu_{\beta}&= &-
\frac{
1}{2}\int\bigl(\partial(f-1)\bigr)^2 \,d\mu_{\beta}
\\
&=&-\frac{
1}{2}\int(\partial f)^2 \,d\mu_{\beta}.
\end{eqnarray*}
Thus, our main task is to find
constants $c_3(p,A)>0$ and $\widetilde\sigma(p)\in(0,1/2)$ such that
for any $\alpha>0$ and $\beta\geq1$ with $\alpha\beta^2\leq
\widetilde
\sigma(p)$, we have,
for any $f\in\mathcal{C}^2(\TT)$,
%
%
%e4.3 #&#
\begin{eqnarray}
\label{borned}
\hspace*{-20pt}&& \biggl\llvert\int R_{\alpha,\beta} [ f-1 ] (f-1) \,d\mu_{\beta}\biggr\rrvert
\leq {c_3(p,A)
\alpha^{\widetilde a(p)}\beta^3} \biggl(\int(\partial f)^2 \,d\mu_{\beta
}+\int(f-1)^2 \,d\mu_\beta\biggr).
\end{eqnarray}
By definition, we have for any $f\in\mathcal{C}^2(\TT)$ (but what
follows is valid for $f\in\mathcal{C}^1(\TT)$),
%
%e4.4 #&#
\begin{eqnarray}
R_{\alpha,\beta}[f](x)&=& \frac{1}\alpha\int f\bigl(
\gamma\bigl(x,y,(p/2)\alpha\beta d^{p-1}(x,y)\bigr)\bigr)-f(x) \nu(dy)+
\frac{\beta}2 U_p'(x)f'(x)\nonumber
\\
\eqntext{\forall x\in\TT.}
\end{eqnarray}
To evaluate this quantity, on one hand, recall that we have for any
$x\in\TT$,
\begin{eqnarray*}
U_p'(x)&=&-p\int_{\TT}
d^{p-1}(x,y)\dot{\gamma}(x,y,0) \nu(dy)
\end{eqnarray*}
and on the other hand, write that for any $x\in\TT$ and $y\in\TT
\setminus\{ x\}$,
\begin{eqnarray*}
&& f\bigl(\gamma\bigl(x,y,(p/2)\alpha\beta d^{p-1}(x,y)\bigr)
\bigr)-f(x)
\\
&&\quad =\frac{p}{2}\alpha\beta\int_0^1f'
\bigl(\gamma\bigl(x,y,(p/2)\alpha\beta d(x,y) u\bigr)\bigr)d^{p-1}(x,y)
\dot{\gamma}(x,y,0) \,du.
\end{eqnarray*}
Writing $s\ffff (p/2)\alpha\beta$ and considering again the operators
introduced in (\ref{Tys2}) (now for any $p\geq1$),
it follows that
\begin{eqnarray*}
&& \int R_{\alpha,\beta} [ f-1 ] (f-1) \,d\mu_{\beta}
\\
&&\quad = \frac{p\beta}{2}\int_0^1 du\int\nu(dy)
\int\mu_\beta(dx) \bigl(T_{y,s
u}\bigl[f'
\bigr](x)-f'(x)\bigr) \bigl(f(x)-1\bigr)d^{p-1}(x,y)\dot{
\gamma}(x,y,0)
\\
&&\quad = \frac{p\beta}{2}\int_0^1 du\int\nu(dy)
\int\mu_\beta(dx) \bigl(T_{y,s
u}\bigl[f'
\bigr](x)-f'(x)\bigr)g_y(x),
\end{eqnarray*}
where for any fixed $y\in\TT$,
%
%
%e4.5 #&#
\begin{eqnarray}
\label{gy} \forall x\in\TT\setminus\{y\},\qquad g_y(x)&\ffff &
\bigl(f(x)-1\bigr)d^{p-1}(x,y)\dot{\gamma}(x,y,0)
\end{eqnarray}
(with, e.g., the
convention that $g_y(y')\ffff 0$).
Let us also fix the variable $u\in[0,1]$ for a while.

We begin by considering the case where $p\geq2$.
By definition of $T^*_{y,s u}$ (discussed in Section~\ref{ri}), we have
%
%
%e4.6 #&#
\begin{eqnarray}
\label{Tfprimeg} \int\bigl(T_{y,s u}\bigl[f'
\bigr](x)-f'(x)\bigr)g_y(x) \mu_\beta(dx)&=&
\int f'(x) \bigl(T_{y,s u}^*[g_y](x)-g_y(x)
\bigr) \mu_\beta(dx)
\nonumber\\[-8pt]\\[-8pt]\nonumber
&=& I_1(y,u)+I_2(y,u),
\end{eqnarray}
where for any $y\in\TT$,
\begin{eqnarray}
\label{I2yu}
\nonumber
I_1(y,u)&\ffff & \int_{B(y,\pi-su\pi^{p-1})}
f'(x) \bigl(T_{y,s
u}^*[g_y](x)-g_y(x)
\bigr) \mu_\beta(dx),
\nonumber\\[-8pt]\\[-8pt]\nonumber
I_2(y,u)&\ffff &- \int_{B(y',su\pi^{p-1})} f'(x)g_y(x)
\mu_\beta(dx)
\end{eqnarray}
(recall from Sections~\ref{p0} and~\ref{p3} that for any measurable
function $g$,
$T^*_{y,s}[g]$ vanishes on $B(y',su\pi^{p-1})$).
The first integral is treated through the Cauchy--Schwarz inequality,
\begin{eqnarray*}
\bigl\llvert I_1(y,u)\bigr\rrvert&\leq& \sqrt{\int
\bigl(f'\bigr)^2 \,d\mu_\beta}\sqrt{
\int_{B(y,\pi-su\pi^{p-1})} \bigl(T_{y,s u}^*[g_y]-g_y
\bigr)^2 \mu_\beta}
\end{eqnarray*}
and Lemmas~\ref{Tstarpa0} and~\ref{Tstarpa3}, at least if $s\beta>0$
is smaller than a certain constant $\widetilde\sigma(p)\in(0,/12)$.
It follows that for a universal constant $k>0$,
we have
\begin{eqnarray*}
\int_{\TT\times[0,1]} \bigl\llvert I_1(y,u)\bigr\rrvert\nu
(dy)\,du&\leq& ks^2\beta^2 \biggl(\int(\partial
f)^2 \,d\mu_{\beta}+\int(f-1)^2 \,d\mu
_\beta\biggr)\int_0^1
u^2 \,du
\\
&=& \frac{k}{2}s^2\beta^2 \biggl(\int(\partial
f)^2 \,d\mu_{\beta
}+\int f^2 \,d\mu_\beta
\biggr)
\\
&\leq& \frac{k}{4}s\beta\biggl(\int(\partial f)^2 \,d\mu_{\beta
}+\int f^2 \,d\mu_\beta\biggr),
\end{eqnarray*}
bound going in the direction of (\ref{borned}).

Next, we turn to the integral $I_2(y,u)$. We cannot deal with it
uniformly over $y\in\TT$
but we get a convenient bound by integrating it with respect to $\nu(dy)$.
Recalling
that under the assumption~(\ref{aA}) the density of $\nu$ with
respect to $\lambda$
is bounded by $ 1+A\pi$, it appears that
\begin{eqnarray}
\label{J3}
\int\bigl\llvert I_2(y,u)\bigr\rrvert\nu(dy)
&\leq& \frac{1+A\pi }{2\pi
}\int_{-\pi}^{\pi} \bigl\llvert
I_2(y,u)\bigr\rrvert \,dy\nonumber
\\
&\leq& \frac{1+A\pi}{2\pi}\int_{\TT}\,dy \int
_{B(y',su\pi^{p-1})}\bigl\llvert f'(x)\bigr\rrvert\bigl
\llvert g_y(x)\bigr\rrvert\mu_\beta(dx)
\nonumber\\[-8pt]\\[-8pt]\nonumber
&\leq& \frac{1+A\pi}{2}\pi^{p-2}\int_\TT
\mu_\beta(dx)\bigl\llvert f'(x)\bigr\rrvert\bigl\llvert
f(x)-1\bigr\rrvert\int_{B(x',su\pi^{p-1})}\un\, dy\nonumber
\\
&=& (1+A\pi)\pi^{2p-3}su\int_\TT\bigl
\llvert f'\bigr\rrvert\llvert f-1\rrvert \,d\mu_\beta.\nonumber
\end{eqnarray}
The Cauchy--Schwarz inequality and
integration with respect to $\un_{[0,1]}(u) \,du$ lead again
to a bound contributing to (\ref{borned}).

It is time to consider the cases where $p\in[1,2)$.
We will rather decompose the left-hand side of (\ref{Tfprimeg}) into three parts.
Let us extend the notation $\widetilde u_{\pm}\ffff  y\pm
(su)^{\afrac{1}{2-p}}$ from Section~\ref{p2} to all $p\in[1,2)$.
Next, we modify the definition (\ref{gy}) by introducing $\widetilde
g_y(x)\ffff \un_{[y-\pi,\widetilde u_-]\sqcup[\widetilde u_+,y+\pi
]}(x)g_y(x)$.
Then we write
\begin{eqnarray*}
\int\bigl(T_{y,s u}\bigl[f'\bigr](x)-f'(x)
\bigr)g_y(x) \mu_\beta(dx) &=& \widetilde
I_1(y,u)+ I_2(y,u)+I_3(y,u),
\end{eqnarray*}
where
\begin{eqnarray*}
\widetilde I_1(y,u)&\ffff & \int_{B(y,\pi-su\pi^{p-1})}
f'(x) \bigl(T_{y,s
u}^*[\widetilde g_y](x)-g_y(x)
\bigr) \mu_\beta(dx),
\\
I_2(y,u)&\ffff &- \int_{B(y',su\pi^{p-1})} f'(x)g_y(x)
\mu_\beta(dx),
\\
I_3(y,u)&\ffff & \int_{[\widetilde u_-,\widetilde u_+]} T_{y,s
u}
\bigl[f'\bigr](x)g_y(x) \mu_\beta(dx).
\end{eqnarray*}
The treatment of $\widetilde I_1(y,u)$
is similar to that of $ I_1(y,u)$,
with Lemmas~\ref{Tstarpa1} and~\ref{Tstarpa2} (where a preliminary
integration with respect to $\nu(dy)$ was necessary)
replacing Lemmas~\ref{Tstarpa0} and~\ref{Tstarpa3}.

Concerning $I_2(y,u)$, it is bounded in the same manner as the
corresponding quantity defined in (\ref{I2yu}).

It seems that the most convenient way to deal with $I_3(y,u)$ is to
first integrate it with respect to $\un_{[0,1]}(u) \nu(dy)\,du$.
Taking into account that $\llVert \nu\rrVert _\infty\leq
(1+A\pi)$ and using
Cauchy--Schwarz inequality, we get
\begin{eqnarray*}
&& \int\bigl\llvert I_3(y,u)\bigr\rrvert\un_{[0,1]}(u)
\nu(y)\,du
\\
&&\quad  \leq \frac{1+A\pi}{2\pi} \int\bigl\llvert I_3(y,u)\bigr
\rrvert\un_{[0,1]}(u) \,dy\,du
\\
&&\quad \leq \frac{1+A\pi}{2\pi}\sqrt{\int\un_{[\widetilde
u_-,\widetilde
u_+]} (x) \bigl(T_{y,su}
\bigl[f'\bigr](x)\bigr)^2\un_{[0,1]}(u)
\mu_\beta(dx)\,dy\,du}
\\
&&\qquad{}\times \sqrt{\int\un_{[\widetilde u_-,\widetilde u_+]} (x)g_y^2(x)\un
_{[0,1]}(u) \mu_\beta(dx)\,dy\,du}.
\end{eqnarray*}
The last factor can be rewritten under the form
\begin{eqnarray}
\label{ssss}
\nonumber
&& \sqrt{\int\mu_\beta(dx)\int
\un_{[x-
s^{\afrac{1}{2-p}},x+
s^{\afrac{1}{2-p}} ]} (y)g_y^2(x) \,dy}\nonumber
\\
&&\quad  \leq
\pi^{p-1}\sqrt{\int\mu_\beta(dx) \bigl(f(x)-1
\bigr)^2\int_{x-
s^{\afrac{1}{2-p}}}^{x+
s^{\afrac{1}{2-p}}} \,dy}
\\
&&\quad =\pi\sqrt{2 s^{\afrac{1}{2-p}} }\sqrt{\int(f-1)^2 \,d\mu_\beta}.\nonumber
\end{eqnarray}
So it remains to consider the term
%
%
%e4.7 #&#
\begin{eqnarray}
\label{premsqrt}
&& \int\un_{[\widetilde u_-,\widetilde u_+]} (x) \bigl(T_{y,su}
\bigl[f'\bigr](x)\bigr)^2\un_{[0,1]}(u)
\mu_\beta(dx)\,dy\,du
\nonumber\\[-8pt]\\[-8pt]\nonumber
&&\quad =\frac{1}{2\pi}\int\un_{[\widetilde u_-,\widetilde u_+]} (x)T_{y,su}
\bigl[\bigl(f'\bigr)^2\bigr](x) \mu_\beta(x)
\un_{[0,1]}(u) \,dy\,du
\end{eqnarray}
(where as a function, $\mu_\beta$ stands for the density of the
measure $\mu_\beta$ with respect to $\lambda$).
Remember that for any measurable function $h$, we have\vspace*{1pt}
$T_{y,su}[h](x)\ffff  h(x+ sud^{p-1}(x,y)\* \dot{\gamma}(x,y,0))$. For
$x\in[\widetilde u_-,\widetilde u_+]$, we have $d(x,y)\leq
(su)^{\afrac{1}{2-p}}$
and it follows that $d(x,x+ sud^{p-1}(x,y)\dot{\gamma}(x,y,0))\leq
(su)^{\vafrac{3-p}{2-p}}$.
Taking into account that $\llVert U_p'\rrVert _\infty\leq\pi
^{p-1}$, we can
then a universal constant $k>0$ such that for $0\leq s\beta\leq
\widetilde
\sigma(p)$ (for an appropriate constant $\widetilde\sigma(p)\in(0,1/2)$)
and $x\in\TT$,
we have
$\mu_{\beta}(x)/\mu_{\beta}(x+ sud^{p-1}(x,y)\dot{\gamma
}(x,y,0))\leq k$.
This leads us to consider the function $h$ defined by
%
%
%e4.8 #&#
\begin{eqnarray}
\label{hfprime} \forall x\in\TT,\qquad h(x)&\ffff & \bigl(f'(x)
\bigr)^2\mu_\beta(x),
\end{eqnarray}
since up to a universal constant, we have to find an upper bound of
\begin{eqnarray*}
&& \int\un_{[\widetilde u_-,\widetilde u_+]} (x)T_{y,su}[h](x)\un
_{[0,1]}(u) \,dx \,dydu
\\
&&\quad \leq \int_{-\pi}^\pi dx\int_{x-
s^{\afrac{1}{2-p}}}^{x+
s^{\afrac{1}{2-p}}}dy
\int_{x-sd^{p-1}(x,y)}^{x+sd^{p-1}(x,y)}h(v) \frac
{dv}{sd^{p-1}(x,y)}
\\
&&\quad =\int_\TT H(v)h(v) \,dv,
\end{eqnarray*}
where for any fixed $v\in\TT$,
\begin{eqnarray*}
H(v)&\ffff & \frac{1}s\int_{\TT^2} \un_{ \{d(x,y)\leq
s^{\afrac{1}{2-p}}, d(v,x) \leq s d^{p-1}(x,y) \}}
\frac
{dx\,dy}{d^{p-1}(x,y)}.
\end{eqnarray*}
Let us furthermore fix $x\in\TT$,
\begin{eqnarray*}
\frac{1}s \int_{\TT} \un_{ \{(d(v,x)/s)^{\afrac{1}{p-1}}\leq
d(x,y)\leq
s^{\afrac{1}{2-p}} \}}
\frac{dy}{d^{p-1}(x,y)} &=& \frac{2}{(2-p)s} \biggl( s- \biggl(\frac{d(v,x)}{s}
\biggr)^{\vafrac{2-p}{p-1}} \biggr)_+.
\end{eqnarray*}
The integration of the last right-hand side with respect to $dx$ is
bounded above by
\begin{eqnarray*}
\frac{2}{2-p}\int_0^{(
s^{\afrac{1}{2-p}})^{p-1}s} \,dx &=&
\frac{2
}{2-p} s^{\afrac{1}{2-p}}.
\end{eqnarray*}
Thus, we have found a constant $k(p)>0$ depending on $p\in[1,2)$ such
that (\ref{premsqrt}) is bounded above
by $k(p) s^{\afrac{1}{2-p}}$ under our conditions on $s>0$ and $\beta
\geq1$.
In conjunction with (\ref{ssss}) and definition (\ref{hfprime}), it
enables to conclude to the existence of a constant $k(p,A)>0$,
depending on $p\in[1,2)$ and $A>0$,
such that
\begin{eqnarray*}
\int\bigl\llvert I_3(y,u)\bigr\rrvert\un_{[0,1]}(u)
\nu(y)\,du &\leq& k(p,A)s^{\afrac{1}{2-p}}\sqrt{\int(f-1)^2 \,d\mu_\beta}\sqrt{\int\bigl(f'\bigr)^2 \,d\mu_\beta}.
\end{eqnarray*}
Putting\vspace*{1pt} together all these estimates and taking into account that
$\beta\geq1$, $0<s\beta\leq\widetilde\sigma(p)$ and
$s^{2(p-1)}\geq s^{1/(2-p)}$,
it appears that
\begin{eqnarray*}
&& \biggl\llvert\int_{\TT\times[0,1]} \widetilde
I_1(y,u)+ I_2(y,u)+I_3(y,u) \nu(dy)\,du\biggr
\rrvert
\\
&&\quad \leq k'(p,A) \cases{ \beta s, &\quad if $p=1$ or $p\geq2$,
\vspace*{3pt}\cr
\beta s+ s^{2(p-1)}, &\quad if $p\in(1,2)$}
\\
&&\quad \leq 2k'(p,A) \cases{ \beta s, &\quad if $p=1$ or $p\geq3/2$,
\vspace*{3pt}\cr
\beta s+ s^{2(p-1)}, &\quad if $p\in(1,3/2)$,}
\end{eqnarray*}
for another constant $k'(p,A)>0$, depending on $p\in[1,2)$ and $A>0$.
This finishes the proof of~(\ref{borned}).
\end{pf}

To conclude the proof of Proposition~\ref{Iprime}, we must be able to compare,
for any $\beta\geq0$ and any $f\in\mathcal{C}^1(\TT)$, the energy
$\mu_\beta[(\partial f)^2]$
and the variance $\Var(f,\mu_\beta)$. This task was already done by~\cite{MR995752}, let us recall their result.
%

%pr4.2 #&#
\begin{pro}\label{HKS2}
Let $U_p$ be a $\mathcal{C}^1$ function on a compact Riemannian
manifold $M$ of dimension $m\geq1$.
Let $b(U_p)\geq0$ be the associated constant as in (\ref{bU}).
For any $\beta\geq0$, consider the Gibbs measure $\mu_\beta$
given in (\ref{mub}). Then there exists a constant $C_M>0$, depending
only on $M$,
such that the following Poincar\'e inequalities are satisfied:
\begin{eqnarray*}
\forall\beta\geq0, \forall f\in\mathcal{C}^1(M),\qquad \Var(f,\mu
_\beta) &\leq& C_M\bigl[1\vee\bigl(\beta\bigl\llVert
U_p'\bigr\rrVert_\infty\bigr)
\bigr]^{5m-2}\exp\bigl(b(U_p)\beta\bigr)\mu_\beta
\bigl[\llvert\nabla f\rrvert^2\bigr].
\end{eqnarray*}
\end{pro}

We can now come back to the study of the evolution of the quantity
$I_t=\Var(f_t,\mu_{\beta_t})$, for $t>0$. Indeed applying Lemma~\ref
{Lffm} and Proposition~\ref{HKS2} with $\alpha=\alpha_t$,
$\beta=\beta_t$ and $f=f_t$, we get at any time $t>0$
such that $\beta_t\geq1$ and $\alpha_t\beta_t^2\leq\varsigma(p)$,
\begin{eqnarray*}
&& \int L_{\alpha_t,\beta_t} [ f_t-1 ] (f_t-1) \,d\mu_{\beta_t}
\\
&&\quad \leq-c_4\beta_t^{-3}\exp\bigl(-b(U_p)
\beta_t\bigr) \bigl(1 -2c_3(p,A)\alpha_t^{\widetilde a(p)}
\beta_t^3 \bigr)I_t+c_3(p,A)
\alpha_t^{\widetilde
a(p)}\beta_t^3I_t
\\
&&\quad \leq -\bigl(c_4\beta_t^{-3}\exp
\bigl(-b(U_p)\beta_t\bigr)-c_5(p,A)
\alpha_t^{\widetilde
a(p)}\beta_t^3
\bigr)I_t,
\end{eqnarray*}
where $c_4\ffff (16\pi^3C_\TT)^{-1}$ and $c_5(p,A)\ffff  c_3(p,A)(1+2c_4)$.

Taking into account Lemma~\ref{palnmu}, the computations preceding
Lemma~\ref{Lffm}
and Proposition~\ref{mbl}, one can find constants
$c_1(p,A),c_2(p,A)>0$ and $\varsigma(p)\in(0,1/2)$ such that
Proposition~\ref{Iprime}
is satisfied.

This result leads immediately to conditions insuring the convergence
toward 0
of the quantity $I_t$ for large times $t>0$.
%

%pr4.3 #&#
\begin{pro}\label{condab}
Let
$\alpha\st\RR_+\rightarrow\RR_+^*$ and
$\beta\st\RR_+\rightarrow\RR_+$ be schemes as at the beginning of this
section and assume
\begin{eqnarray*}
\lim_{t\rightarrow+\infty}\beta_t&=&+\infty,
\\
\int_0^{+\infty}(1\vee\beta_t)^{-3}
\exp\bigl(-b(U_p)\beta_t\bigr) \,dt &=&+\infty
\end{eqnarray*}
and that for large times $t>0$,
\begin{eqnarray*}
\max\bigl\{ \alpha_t^{a(p)}\beta_t^4,
\alpha_t^{\widetilde
a(p)}\beta_t^3, \bigl
\llvert\beta_t'\bigr\rrvert\bigr\}&\ll& \exp
\bigl(-b(U_p)\beta_t\bigr)
\end{eqnarray*}
(where $a(p)>0$ and $\widetilde a(p)>0$ are defined in Propositions
\ref{mbl}
and~\ref{Iprime}). Then we are assured of
\begin{eqnarray*}
\lim_{t\rightarrow+\infty} I_t&=&0.
\end{eqnarray*}
\end{pro}

\begin{pf}
The differential equation of Proposition~\ref{Iprime} can be rewritten
under the form
%
%
%e4.9 #&#
\begin{eqnarray}
\label{Fee} F_t'&\leq& -\eta_t
F_t+\epsilon_t,
\end{eqnarray}
where for any $t>0$,
\begin{eqnarray*}
F_t&\ffff & \sqrt{I_t},
\\
\eta_t&\ffff & c_1(p,A) \bigl(\beta_t^{-3}
\exp\bigl(-b(U_p)\beta_t\bigr)-\alpha_t^{\widetilde a(p)}
\beta_t^3-\bigl\llvert\beta_t'
\bigr\rrvert\bigr)/2,
\\
\epsilon_t&\ffff & c_2(p,A) \bigl(\alpha_t^{a(p)}
\beta_t^4+\bigl\llvert\beta_t'
\bigr\rrvert\bigr)/2.
\end{eqnarray*}
The assumptions of the above proposition
imply that for $t\geq0$ large enough,
$\beta_t\geq1$ and $\alpha_t\beta_t^2\leq\varsigma(p)$, where
$\varsigma(p)\in(0, 1/2)$ is as in
Proposition~\ref{Iprime}.
This
ensures that there exists $T> 0$ such that (\ref{Fee})
is satisfied for any $t\geq T$ (and also $F_T<+\infty$). We deduce that
for any $t\geq T$,
%
%
%e4.10 #&#
\begin{eqnarray}
\label{Fprime} F_t&\leq& F_T\exp\biggl(-\int
_T^t \eta_s \,ds \biggr) +\int
_T^t \epsilon_s\exp\biggl(-\int
_s^t \eta_u \,du \biggr) \,ds.
\end{eqnarray}
It appears that $\lim_{t\rightarrow+\infty} F_t =0$ as soon as
\begin{eqnarray*}
\int_T^{+\infty}\eta_s \,ds&=&+\infty,
\\
\lim_{t\rightarrow+\infty} \epsilon_t/\eta_t&=& 0.
\end{eqnarray*}
The above assumptions were chosen to ensure these properties.
\end{pf}

In particular, remarking that $a(p)\leq\widetilde a(p)$ for any $p\geq1$,
the schemes given in (\ref{ab}) satisfy the hypotheses of the previous
proposition, so that under the conditions of Theorem~\ref{t1}, we get
\begin{eqnarray*}
\lim_{t\rightarrow+\infty} I_t&=&0.
\end{eqnarray*}
Let us deduce (\ref{mN1})
for any neighborhood $\cN$ of the set $\cM_p$ of the global minima of
$U_p$. From Cauchy--Schwarz inequality we have for any $t>0$,
\begin{eqnarray*}
\llVert m_t-\mu_{\beta_t}\rrVert_{\mathrm{tv}}&=& \int
\llvert f_t-1\rrvert\mu_{\beta_t}
\\
&\leq& \sqrt{I_t}.
\end{eqnarray*}
An equivalent definition of the total variation norm states that
\begin{eqnarray*}
\llVert m_t-\mu_{\beta_t}\rrVert_{\mathrm{tv}}&=&2\max
_{A\in\cT}\bigl\llvert m_t(A)-\mu_{\beta_t}(A)
\bigr\rrvert,
\end{eqnarray*}
where $\cT$ is the Borelian $\sigma$-algebra of $\TT$.
It follows that (\ref{mN1}) reduces to
\begin{eqnarray*}
\lim_{\beta\rightarrow+\infty} \mu_\beta(\cN)&=&1,
\end{eqnarray*}
for any neighborhood $\cN$ of $\cM_p$, property which is immediate
from the definition (\ref{mub}) of the Gibbs measures $\mu_\beta$
for $\beta\geq0$.
This finishes the proof of Theorem~\ref{t1}.
%

%re4.1 #&#
\begin{rem}\label{discu1}
Under mild conditions, the results of \cite{MR602391} enable to go
further, because
he identifies the weak limit $\mu_{\infty}$ of the Gibbs measures
$\mu
_\beta$
as $\beta$ goes to $+\infty$.
Thus, if one knows, as above, that
\begin{eqnarray*}
\lim_{t\rightarrow+\infty} \llVert m_t-\mu_{\beta_t}\rrVert
_{\mathrm
{tv}}&=&0,
\end{eqnarray*}
then one gets that $m_t$ also weakly converges toward $\mu_{\infty}$ for
large times $t>0$.
The weight given by $\mu_{\infty}$ to a point $x\in\cM_p$ is inversely
related to the value of $\sqrt{U_p''(x)}$
and in this respect Lemma~\ref{regU} is useful (still assuming that
$\nu$ admits a continuous density).

First note that for any $x\in\cM_p$, we have $U_p''(x)\geq0$, since
$x$ is a global minima of $U_p$,
and by consequence $\nu(x')\leq1$.
Next, assume that we have for any $x\in\cM_p$, $\nu(x')<1$. It
follows that $\cM_p$ is discrete
and by consequence finite, since $\TT$ is compact.
This property was already noted by \cite{2011arXiv11082141H}, among
other features of intrinsic means on the circle.
Then we deduce from \cite{MR602391} that
\begin{eqnarray*}
\mu_\infty&=&\frac{1}Z\sum_{x\in\cM_p}
\frac{1}{\sqrt{1-\nu
(x')}}\delta_{x},
\end{eqnarray*}
where $Z\ffff \sum_{x\in\cM_p}(1-\nu(x'))^{-1/2}$ is the normalizing
factor.

In this situation, $\cL(X_t)$ concentrates for large times $t>0$ on
all the $p$-means of $\nu$.
Thus, to find all of them with an important probability, one should
sample independently several
trajectories of $X$, for example, starting from a fixed point $X_0\in\TT$.
\end{rem}

%
%re4.2 #&#
\begin{rem}\label{discu2}
Similarly to the approach presented, for instance, in \cite
{MR1275365,MR1425361},
we could have studied the evolution of $(E_t)_{t>0}$, which are the
relative entropies of the time marginal laws with respect to the
corresponding instantaneous Gibbs measures, namely
\begin{eqnarray*}
\forall t>0,\qquad E_t&\ffff & \int\ln\biggl(\frac{m_t}{\mu_{\beta
_t}}
\biggr) \,dm_t.
\end{eqnarray*}
To get a differential
inequality satisfied by these functionals, the spectral gap estimate of
\cite{MR995752} recalled in Proposition~\ref{HKS2}
must be replaced by the corresponding logarithmic Sobolev constant estimate,
which is proven in the same article of \cite{MR995752}.
\end{rem}

%s5 #&#
\section{Extension to all probability measures \texorpdfstring{$\nu$}{$nu$}}\label{sec5}

Our main task here is to adapt the computations of the two previous
sections in order to prove Theorem~\ref{t2}.
As in the statement of this result, it is better for simplicity of the
exposition to restrict ourselves to the important and illustrative case $p=2$;
the general situation will be alluded to in the last remark of this section.

We begin by remarking that the algorithm $Z$ described in the
\hyperref[sec1]{Introduction} evolves similarly
to the process $X$, if we allow the probability measure $\nu$ to
depend on time.
More precisely, for any $\kappa>0$, consider
the probability measure $\nu_\kappa$ given by
%
%
%e5.1 #&#
\begin{eqnarray}
\label{nk} \forall z\in M,\qquad\nu_\kappa(dz)&\ffff &\int\nu(dy)
K_{y,\kappa
}(dz),
\end{eqnarray}
where the kernel on $M$, $(y,dz)\mapsto K_{y,\kappa}(dz)$ was defined
before the statement of Theorem~\ref{t2}.
For $\alpha>0$, $\beta\geq0$ and $\kappa>0$,
let us denote by $L_{\alpha,\beta,\kappa}$ the generator defined
in (\ref{Lab}), where $\nu$ is replaced by $\nu_\kappa$.
Then the law of $Z$ is solution of the time-inhomogeneous martingale problem
associated to the family of generators $(L_{\alpha_t,\beta_t,\kappa
_t})_{t\geq0}$.
This observation leads us to introduce the potentials
\begin{eqnarray*}
\forall\kappa>0, \forall x\in M,\qquad U_{2,\kappa}(x)&\ffff & \int
d^2(x,y) \nu_\kappa(dy),
\end{eqnarray*}
as well as the associated Gibbs measures:
\begin{eqnarray*}
\forall\beta\geq0, \forall\kappa>0,\qquad\mu_{\beta,\kappa
}(dx)&\ffff &
Z_{\beta,\kappa}^{-1} \exp\bigl(-\beta U_{2,\kappa}(x)\bigr)
\lambda(dx),
\end{eqnarray*}
where $Z_{\beta,\kappa}$ is the renormalization constant.

Denote by $m_t$ the law of $Z_t$ for any $t\geq0$.
The proof of Theorem~\ref{t2} is then similar to that of Theorem~\ref
{t1} and
relies on the investigation of the evolution of
%
%
%e5.2 #&#
\begin{eqnarray}
\label{Jt} \forall t>0,\qquad\cI_t&\ffff & \int\biggl(
\frac{m_t}{\mu_{\beta_t,\kappa_t}}-1 \biggr)^2 \,d\mu_{\beta_t,\kappa_t},
\end{eqnarray}
which play the role of the quantities defined in (\ref{It}).

While the above program was presented for a general compact Riemannian
manifold $M$,
we again restrict ourselves to the situation $M=\TT$.

We first need some estimates on the probability measures $
\nu_\kappa$, for $\kappa>0$.
%

%le5.1 #&#
\begin{lem}\label{nukappa}
For any $\kappa>0$, $\nu_\kappa$ admits a density with respect to
$\lambda$,
still denoted $\nu_\kappa$. Furthermore we have, for any $\kappa>
1/\pi$,
\begin{eqnarray*}
\llVert\nu_\kappa\rrVert_\infty&\leq&2\pi\kappa,
\\
\llVert\partial\nu_\kappa\rrVert_\infty&\leq&2\pi
\kappa^2,
\end{eqnarray*}
where $\partial\nu_\kappa$ stands for the weak derivative (so that
the last
norm $\llVert \cdot\rrVert _\infty$ is the essential
supremum norm with
respect to $\lambda$).
\end{lem}

\begin{pf}
When $M=\TT$, for any $\kappa>0$, the kernel $K_{\cdot,\kappa
}(\cdot)$ corresponds to the rolling around $\TT$
of the kernel defined on $\RR$ by $(y,dz)\mapsto\kappa(1-\kappa
\llvert z-y\rrvert )_+ \,dz$.
In particular for any $y\in\TT$, $K_{y,\kappa}(\cdot)$ is
absolutely continuous with respect to $\lambda$
and (\ref{nk}) shows that the same is true for $\nu_\kappa$.
If furthermore $\kappa> 1/\pi$, from this definition we can write for
any $z\in\TT$,
\begin{eqnarray*}
\nu_\kappa(dz)&=&\kappa\biggl(\int_{z-1/\kappa}^{z+1/\kappa}
\bigl(1-\kappa d(y,z)\bigr)_+ \nu(dy) \biggr) \,dz,
\end{eqnarray*}
namely, almost everywhere with respect to $\lambda(dz)$,
\begin{eqnarray*}
\nu_\kappa(z)&=&2\pi\kappa\int_{z-1/\kappa}^{z+1/\kappa
}
\bigl(1-\kappa d(y,z)\bigr)_+ \nu(dy)
\\
&\leq&2\pi\kappa\int_{z-1/\kappa}^{z+1/\kappa} \nu(dy)
\\
&\leq& 2\pi\kappa.
\end{eqnarray*}
Next, for almost every $x,y\in\TT$, we have
\begin{eqnarray*}
\bigl\llvert\nu_\kappa(x)-\nu_\kappa(y)\bigr\rrvert&\leq& 2
\pi\kappa\int_{\TT} \bigl\llvert\bigl(1-\kappa d(x,z)
\bigr)_+-\bigl(1-\kappa d(y,z)\bigr)_+\bigr\rrvert\nu(dz)
\\
&\leq&2\pi\kappa\int_{\TT} \bigl\llvert1-\kappa d(x,z)-1+
\kappa d(y,z)\bigr\rrvert\nu(dz)
\\
&\leq&2\pi\kappa^2\int_{\TT} \bigl\llvert
d(x,z)- d(y,z)\bigr\rrvert\nu(dz)
\\
&\leq& 2\pi\kappa^2 d(x,y).
\end{eqnarray*}
This proves the second bound.
\end{pf}

An immediate consequence of the last bound is that
for any $x\in\TT$, the map $(1/\pi,+\infty)\ni\kappa\mapsto
U_{2,\kappa}(x)$ is weakly
differentiable and for almost every $\kappa>1/\pi$, $\llvert
\partial
_\kappa U_{2,\kappa}(x)\rrvert \leq2\pi^4\kappa^2$; but one can
do better.
%

%le5.2 #&#
\begin{lem}\label{pakU}
For any $x\in\TT$ and any $\kappa>1/\pi$, we have
\begin{eqnarray*}
\bigl\llvert\partial_\kappa U_{2,\kappa}(x)\bigr\rrvert&\leq&
\frac
{3 \pi
^3}{\kappa}.
\end{eqnarray*}
\end{lem}

\begin{pf}
It is better to come back to the definition of $\nu_\kappa$, to get,
for $x\in\TT$ and $\kappa>1/\pi$ (where $\partial_\kappa$ stands for
weak derivative):
\begin{eqnarray*}
\partial_\kappa U_{2,\kappa}(x)&=&\partial_\kappa\biggl(
2\pi\kappa\int\lambda(dy) d^2(x,y) \int_{\TT}
\bigl(1-\kappa d(y,z)\bigr)_+ \nu(dz) \biggr)
\\
&=&2\pi\int\lambda(dy) d^2(x,y) \int_{\TT}
\nu(dz) \bigl(1-\kappa d(y,z)\bigr)_+
\\
&&{}-2\pi\kappa\int\lambda(dy) d^2(x,y) \int_{y-1/\kappa}^{y-1/\kappa}
\nu(dz) d(y,z).
\end{eqnarray*}
The first term of the right-hand side is equal to $U_{2,\kappa
}(x)/\kappa$
and is bounded by $\llVert U_{2,\kappa}\rrVert _{\infty
}/\kappa\leq\pi
^2/\kappa$.
In absolute value, the second term can be written under the form
\begin{eqnarray*}
2\pi\kappa\int\nu(dz) \int_{z-1/\kappa}^{z-1/\kappa} \lambda(dy)
d^2(x,y) d(y,z) &\leq&2 \pi^3\kappa\int\nu(dz) \int
_{z-1/\kappa}^{z-1/\kappa} \lambda(dy) \llvert y-z\rrvert
\\
&=&\frac{2 \pi^3}{\kappa}.
\end{eqnarray*}\upqed
\end{pf}

The improvement of the estimate of the previous lemma with respect to
the one given before
its statement is important for us, since it enables to obtain that if
$(\beta_t)_{t\geq0}$
and $(\kappa_t)_{t\geq0}$ are $\mathcal{C}^1$ schemes, then we have
%
%
%e5.3 #&#
\begin{eqnarray}
\label{palnmu2} \forall t\geq0,\qquad\bigl\llVert\partial_t \ln(
\mu_{\beta_t,\kappa_t})\bigr\rrVert_{\infty}&\leq& \pi^2\bigl
\llvert\beta_t'\bigr\rrvert+3 \pi^3
\beta_t \bigl\llvert\bigl(\ln(\kappa_t)
\bigr)'\bigr\rrvert.
\end{eqnarray}
This bound replaces that of Lemma~\ref{palnmu} in the present context.
Note that
for the schemes we have in mind and up to mild logarithmic corrections,
we recover a bound of order
$1/(1+t)$ for $\llVert \partial_t \ln(\mu_{\beta_t,\kappa
_t})\rrVert _{\infty
}$, which is compatible with our purposes.

In the same spirit, even if this cannot be deduced directly from Lemma
\ref{pakU}, we have the following.
%

%le5.3 #&#
\begin{lem}
As $\kappa$ goes to infinity, $U_{2,\kappa}$ converges uniformly
toward $U_2$. In particular,
if $b(\cdot)$ is the functional defined in (\ref{bU}), then we have
\begin{eqnarray*}
\lim_{\kappa\rightarrow+\infty}b(U_{2,\kappa})&=&b(U_2).
\end{eqnarray*}
\end{lem}

\begin{pf}
Since $\llVert \partial U_{2,\kappa}\rrVert _{\infty}\leq
2\pi$, for any $\kappa
>0$, it appears that
$(U_{2,\kappa})_{\kappa>0}$ is an equicontinuous family of mappings.
It is besides clear that $\nu_\kappa$ weakly converges toward $\nu$
as $\kappa$ goes to infinity,
so that $U_{2,\kappa}(x)$ converges toward $U_2(x)$ for any fixed
$x\in\TT$.
Compactness of $\TT$ and the Arzel\`{a}--Ascoli theorem then enable to conclude
to the uniform of $U_{2,\kappa}$ toward $U_2$ as $\kappa$ goes to infinity.
The second assertion of the lemma is an immediate consequence of
this convergence.
\end{pf}

Consider for the evolution of the inverse temperature the scheme
\begin{eqnarray*}
\forall t\geq0,\qquad\beta_t&\ffff & b^{-1}\ln(1+t),
\end{eqnarray*}
where $b>b(U_2)$ and denote $\rho\ffff (1+b(U_2)/b)/2<1$. Assume that
the scheme $(\kappa_t)_{t\geq0}$
is such that $\lim_{t\rightarrow+\infty}\kappa_t=+\infty$. Then
from the above
lemma and Proposition~\ref{HKS2}
(recall that $\llVert \partial U_{2,\kappa}\rrVert _{\infty
}\leq2\pi$, for any
$\kappa>0$), there exists a time $T>0$
such that for any $t\geq T$,
%
%
%e5.4 #&#
\begin{eqnarray}
\label{rho} \forall f\in\mathcal{C}^1(\TT),\qquad
\frac{2}{(1+t)^{\rho}}\Var(f,\mu_{\beta_t,\kappa_t})&\leq& \mu_{\beta
_t,\kappa_t}\bigl[(
\partial f)^2\bigr].
\end{eqnarray}
Like (\ref{palnmu2}), this crucial estimate for the investigation of
the evolution of the quantities (\ref{Jt})
still does not explain the requirement that $k\in(0,1/2)$ in Theorem
\ref{t2}.
Its justification comes from the next result, which replaces
Proposition~\ref{mbl0} in the present situation.
%

%pr5.1 #&#
\begin{pro}\label{mbl4}
For $\alpha>0$, $\beta\geq0$ and $\kappa>0$, let $L_{\alpha,\beta
,\kappa}^*$
be the adjoint operator of $L_{\alpha,\beta,\kappa}$ in $\LL^2(\mu
_{\beta,\kappa})$.
There exists a constant $C_1>0$ such that
for any $\beta\geq1$, $\kappa\geq1$ and $\alpha\in(0,(2\beta
)^{-1}\wedge(\beta^3(\beta+\kappa))^{-1/2})$, we have
\begin{eqnarray*}
\bigl\llVert L_{\alpha,\beta,\kappa}^*\un\bigr\rrVert_\infty&\leq
&C_1\alpha\beta^2\bigl(\beta^2+
\kappa^2\bigr).
\end{eqnarray*}
\end{pro}

\begin{pf}
It is sufficient to replace $U_2$ by $U_{2,\kappa}$ in the proofs of
Section~\ref{ri},
in particular note that (\ref{Lstarun}) still holds.
From Lemma~\ref{regU} and the first part of Lemma~\ref{nukappa}, it
appears that (\ref{Usec})
has to be replaced by
\begin{eqnarray*}
\forall\kappa\geq1,\qquad\bigl\llVert U_{2,\kappa}''
\bigr\rrVert_\infty&\leq& 4\pi\kappa.
\end{eqnarray*}
Instead of (\ref{expUU}), we deduce that for any $x,y\in\TT$ and
$\alpha$, $\beta$ and $\kappa$
as in the statement of the proposition,
\begin{eqnarray*}
&& \exp\biggl(\beta\biggl[U_{2,\kappa}(x)-U_{2,\kappa}
\biggl(x-\frac{\alpha\beta}{1-\alpha\beta}(y-x) \biggr) \biggr] \biggr)
\\
&&\quad =1+ \frac{ \alpha\beta^2}{1-\alpha\beta} U_{2,\kappa}'(x) (y-x)+\cO
\bigl(
\alpha^2\beta^3(\beta+\kappa)\bigr).
\end{eqnarray*}
Keeping following the computations of the same proof, we end up with
\begin{eqnarray*}
L_{\alpha,\beta,\kappa}^*\un(x) &=&\frac{\beta}{1-\alpha\beta} \frac
{1}{2\pi\alpha\beta}\int
_{x'-\alpha\beta\pi}^{x'+\alpha
\beta\pi} \nu_\kappa\bigl(x'
\bigr)-\nu_\kappa(y) \,dy +\cO\bigl(\alpha\beta^3(\beta+
\kappa)\bigr).
\end{eqnarray*}
To estimate the last integral, we resort to the second part of
Lemma~\ref{nukappa}: we get
\begin{eqnarray*}
\biggl\llvert\int_{x'-\alpha\beta\pi}^{x'+\alpha\beta\pi} \nu_\kappa
\bigl(x'\bigr)-\nu_\kappa(y) \,dy\biggr\rrvert&\leq& 2\pi
\kappa^2\int_{x'-\alpha\beta\pi}^{x'+\alpha\beta\pi
} \bigl\llvert
x'-y\bigr\rrvert \,dy
= 2\pi\kappa^2{(\alpha\beta\pi)^{2}}.
\end{eqnarray*}
This leads to the announced bound.
\end{pf}

Similar arguments transform Lemma~\ref{Lffm} into the following.
%

%le5.4 #&#
\begin{lem}
There exists a constant $C_2>0$, such that
for any $\alpha>0$, $\beta\geq1$ and $\kappa\geq1$ with $\alpha
\beta^2\leq1/2$, we have,
for any $f\in\mathcal{C}^2(\TT)$,
\begin{eqnarray*}
\int L_{\alpha,\beta,\kappa} [ f-1 ] (f-1) \,d\mu_{\beta}&\leq&- \biggl(
\frac{
1}2-C_2\alpha\beta^2(\beta+\kappa) \biggr)
\int(\partial f)^2 \,d\mu_{\beta}
\\
&&{} +C_2\alpha\beta^2(\beta+\kappa)\int(f-1)^2
\,d\mu_{\beta}.
\end{eqnarray*}
\end{lem}

\begin{pf}
The modifications with respect to the proof of Lemma~\ref{Lffm}
are very limited: one just needs to take into account the bounds
$\llVert U'_{p,\kappa}\rrVert _\infty\leq2\pi$ and $\llVert \nu
_\kappa\rrVert
_\infty\leq2\pi\kappa$
for $\kappa\geq1$.
Indeed, there are two main changes:

\begin{itemize}
\item in (\ref{LLR}), where the remaining operator has to be defined
by
\begin{eqnarray*}
R_{\alpha,\beta,\kappa}&\ffff & L_{\alpha,\beta
,\kappa} -\tfrac{1}2\bigl(
\partial^2-\beta U'_{p,\kappa}\partial\bigr),
\end{eqnarray*}

\item in (\ref{J3}), the factor $1+A\pi$ must be replaced by
$2\pi\kappa$,
by virtue of the first estimate of Lemma~\ref{nukappa}. It leads to
the supplementary term $\alpha\beta^2\kappa$ in the bound of the
above lemma.\quad\qed
\end{itemize}\noqed
\end{pf}

All the ingredients are collected together to get a differential
inequality satisfied by $(\cI_t)_{\geq0}$. More precisely, under the
requirement that (\ref{rho}) is true for $t\geq T>0$, as well as
$\beta_t\geq1$, $\kappa_t\geq1$ and
$\alpha_t\beta^2_t\sqrt{\kappa_t}\leq1/2$, we get that there
exists a constant $C_3>0$ such that
\begin{eqnarray*}
\forall t \geq T,\qquad\cI_t'&\leq& -
\eta_t \cI_t+\epsilon_t\sqrt{
\cI_t},
\end{eqnarray*}
where for any $t\geq T$,
\begin{eqnarray*}
\eta_t&\ffff &\frac{1}{(1+t)^{\rho}} -C_3\bigl(
\alpha_t\beta_t^2(\beta_t+
\kappa_t)+\bigl\llvert\beta_t'\bigr\rrvert
+\beta_t\bigl\llvert\bigl( \ln(\kappa_t)
\bigr)'\bigr\rrvert\bigr),
\\
\epsilon_t&\ffff & C_3\bigl(\alpha_t
\beta_t^2\bigl(\beta^2_t+
\kappa_t^2\bigr)+\bigl\llvert\beta_t'
\bigr\rrvert+\beta_t\bigl\llvert\bigl( \ln(\kappa_t)
\bigr)'\bigr\rrvert\bigr).
\end{eqnarray*}
Under the assumptions of Theorem~\ref{t2} (already partially used to
ensure the validity of (\ref{rho}) for some $\rho\in(0,1)$), it
appears that as $t$ goes to infinity,
\begin{eqnarray*}
\eta_t&\sim&\frac{1}{(1+t)^{\rho}},
\\
\epsilon_t&=&\cO\biggl(\frac{1}{1+t} \biggr)
\end{eqnarray*}
and this is sufficient to ensure that
\begin{eqnarray*}
\lim_{t\rightarrow+\infty} \cI_t&=&0.
\end{eqnarray*}

The proof of Theorem~\ref{t2} finishes by the arguments given at the
end of Section~\ref{proof}.

%
%
%re5.1 #&#
\begin{rem}\label{discu3}
As it was mentioned at the end of the \hyperref[sec1]{Introduction}, if one does not
want to waste rapidly
the sample $(Y_n)_{n\in\NN}$ (especially if it is not infinite$\ldots$),
one should take the exponent $c$ the smallest possible. From our
assumptions, we necessarily have $c>1$. But the limit case $c=1$ can be
attained: the above proof
shows that the convergence of Theorem~\ref{t2} is also valid for the schemes:
\begin{eqnarray*}
\forall t\geq0,\qquad\cases{ \alpha_t \ffff (1+t)^{-1},
\vspace*{3pt}\cr
\beta_t \ffff  b^{-1}\ln(1+t),
\vspace*{3pt}\cr
\kappa_t \ffff
\ln(2+t).}
\end{eqnarray*}
The drawback is that $\nu$ is not rapidly approached by $\nu_{\kappa_t}$
as $t$ goes to infinity and this may slow down the convergence of the
algorithm toward
$\cN$. Indeed, from the previous computations, it appears that the law
of $Z_t$
is rather close to the set of global minima of $U_{2,\kappa_t}$.
\end{rem}

%
%re5.2 #&#
\begin{rem}\label{t2p1}
The cases $p=1$ and $p\geq2$ can be treated in the same manner, but
for $p\in(1,2)$, one must follow the dependence on $A$
of the constants in the proof of Lemma~\ref{Tstarpa2}.
In the end it only leads to supplementary factors of $\kappa$, so that
Theorem~\ref{t2} is satisfied
with a sufficiently large constant $c$, depending on $p\geq1$ and on
the exponent $k$ entering in the definition
of the scheme $(\kappa_t)_{t\geq0}$.
But before going further in the direction of this generalization, it
would be more rewarding to first check if the
dependence on $p$ of $a_p$ in Theorem~\ref{t1} is just technical or
really necessary.
\end{rem}

%sA #&#
\begin{appendix}\label{append}
\section*{Appendix: Regularity of temporal marginal laws}\label{sec6}
\setcounter{equation}{0}
\setcounter{lem}{0}

Our goal is to see that at positive times, the marginal laws of the
considered algorithms
are absolutely continuous and that if furthermore $\nu\ll\lambda$, then
the corresponding densities belong to $\mathcal{C}^1(\TT)$.
We will also check that this is sufficient to justify the computations
made in Section~\ref{sec4}.

Let $X$ be the process described in the \hyperref[sec1]{Introduction}, for simplicity on
$\TT$,
but the following arguments could be extended to general connected and
compact Riemannian manifolds.
We are going to use the probabilistic construction of $X$
to obtain regularity results on $m_t$, which as usual stands for the
law of $X_t$, for any $t\geq0$.
So for fixed $t>0$, let $T_t$ be the largest jump time of $N^{(\alpha
)}$ in the interval $[0,t]$,
with the convention that $T_t=0$ if there is no jump time in this interval.
Denote by $\xi_t$ the law of $(T_t,X_{T_t})$ on $[0,t]\times\TT$.
Furthermore, let $P_s(x,dy)$ be the law at time $s\geq0$ of the
Brownian motion on $\TT$,
starting at $x\in\TT$.
From the construction given in the \hyperref[sec1]{Introduction}, we have
for any $t>0$,
%
%
%eA.1 #&#
\begin{eqnarray}
\label{mt} m_t(dx)&=&\int_{[0,t]\times\TT}
\xi_t(ds,dz) P_{t-s}(z,dx).
\end{eqnarray}
An immediate consequence is the following.
%

%leA.1 #&#
\begin{lem}\label{acmt}
Let $t>0$ be fixed. About the measurable evolutions $\alpha\st\RR
_+\rightarrow\RR^*_+$
and $\beta\st\RR_+\rightarrow\RR_+$, only assume that
$\inf_{s\in[0,t]} \alpha_s>0$. Then, whatever the probability
measure $\nu$
entering in the definition of $X$, we have that $m_t$ is absolutely continuous.
\end{lem}

\begin{pf}
By the hypothesis on $\alpha$, 0 is the unique atom of $\xi(\cdot
,\TT)$, the distribution of $T_t$ (its mass is
$\xi_t(\{0\},\TT)=\exp(-\int_0^t1/\alpha_s\,ds)$) and $\xi(\cdot
,\TT)$ admits a bounded density
on $(0,t]$. Since furthermore for any $s>0$ and $z\in\TT$,
$P_s(z,\cdot)$ is absolutely continuous, the same is true for $m_t$
due to~(\ref{mt}).
\end{pf}

To go further, we need to strengthen the assumption on $\nu$.
%

%leA.2 #&#
\begin{lem} \label{mtC1}
In addition to the hypotheses of the previous lemma, assume that
$\nu$ admits a bounded density and that $\inf_{s\in[0,t]} \beta
_s>0$. Then for any $t>0$, the density of $m_t$
belongs to~$\mathcal{C}^1(\TT)$.
\end{lem}

\begin{pf}
We begin by recalling a few bounds on the heat kernels $P_s(x,dy)$,
for $s>0$ and $x\in\TT$. We have already mentioned they admit a density,
namely they can be written under the form $p_s(x,y) \,dy$.
Since the Brownian motion on $\TT$ is just the rolling up of the usual
Brownian motion
on $\RR$, we have for any $x\in\TT$,
%
%
%eA.2 #&#
\begin{eqnarray}
\label{sumn} \forall y\in(x-\pi,x+\pi],\qquad p_s(x,y)&=&\sum
_{n\in\ZZ}\frac{\exp(-(y-x+2\pi n)^2/(2s))}{\sqrt
{2\pi s}}.
\end{eqnarray}
From a general bound due to \cite{MR1618694}, we deduce that there
exists a constant $C_0>0$ such that for any $s>0$
and $ y\in(x-\pi,x+\pi]$, we have
\begin{eqnarray*}
\bigl\llvert\partial_yp_s(x,y)\bigr\rrvert&
\leq&C_0 \biggl(\frac{d(x,y)}{s}+ \frac{1}{\sqrt{s}}
\biggr)p_s(x,y).
\end{eqnarray*}
To get an upper bound on $p_s(x,y)=p_s(0,y-x)$, consider separately in
(\ref{sumn})
the sums of $n\in\ZZ_{\sigma}$ and $n\in\ZZ_{-\sigma}\setminus\{
0\}$,
where $\sigma\in\{-,+\}$ is the sign of $y-x$. It appears that for
$s\in(0,t]$,
\begin{eqnarray*}
p_s(x,y)&\leq& 2\sum_{n\in\ZZ_{\sigma}}
\frac{\exp(-(y-x+2\pi
n)^2/(2s))}{\sqrt{2\pi s}}
\\
&\leq&2 \frac{\exp(-(y-x)^2/(2s))}{\sqrt{2\pi s}}\sum_{n\in\ZZ
_+}\exp\bigl(-(2
\pi n)^2/(2s)\bigr)
\\
&\leq&C_1(t) \frac{\exp(-d^2(x,y)/(2s))}{\sqrt{2\pi s}},
\end{eqnarray*}
where $C_1(t)\ffff \sum_{n\in\ZZ_+}\exp(-2(\pi n)^2/t)$.
Taking into account (\ref{mt}) and Lemma~\ref{acmt}, if we were
allowed to differentiate
under the sign integral, we would get for any $x\in\TT$,
%
%
%eA.3 #&#
\begin{eqnarray}
\label{pxmt} \partial_x m_t(x)&=&\int
_{[0,t]\times\TT}\xi_t(ds,dz) \partial_x
p_{t-s}(z,x)
\end{eqnarray}
(where the left-hand side stands for the density of $m_t$ with respect
to $2\pi\lambda$).
Unfortunately, the usual conditions do not apply here, so it is better
to consider the
approximation of the density $m_t$ by $m_{\epsilon,t}$, where for
$\epsilon\in(0,t)$,
\begin{eqnarray*}
\forall x\in\TT,\qquad m_{t,\epsilon}(x)&\ffff &\int_{[0,t-\epsilon
]\times\TT}
\xi_t(ds,dz) p_{t-s}(z,x).
\end{eqnarray*}
There is no difficulty in differentiating this expression under the
sign sum and in the end it appears to be smooth in $x$.
So to get the announced result, it is sufficient to see that $\partial_x
m_{\epsilon,t}(x)$
converges to the right-hand side of (\ref{pxmt}), uniformly in $x\in
\TT$ as $\epsilon$ goes to $0_+$.
Let us prove the stronger convergence
\begin{eqnarray*}
\lim_{\epsilon\rightarrow0_+}\sup_{x\in\TT}\int
_{[t-\epsilon
,t]\times
\TT}\xi_t(ds,dz) \bigl\llvert
\partial_ x p_{t-s}(z,x)\bigr\rrvert&=&0.
\end{eqnarray*}
The assumptions that $\inf_{s\in[0,t]}\alpha_s\beta_s>0$ and that
$\nu$ admits a bounded density imply that the latter is equally true
for $\xi_t(s,\cdot)$, the regular conditional law of $X_{T_t}$
knowing that $T_t=s$, for any $s>0$. We can even find $C_2(t)>0$ such
that $\xi_t(s,dz)\leq C_2(t) \,dz$, uniformly over $s\in(0,t]$ (but a
priori $C_2(t)$ may depend on $t>0$ through $\inf_{s\in[0,t]}\alpha
_s\beta_s$).
In the proof of Lemma~\ref{acmt}, we have already noticed that there
exists $C_3(t)>0$
such that $\xi_t(ds,\TT)\leq C_3(t) \,ds$, for $s\neq0$.
It follows that for $\epsilon\in(0,t)$,
\begin{eqnarray*}
&& \int_{[t-\epsilon,t]\times\TT}\xi_t(ds,dz) \bigl\llvert
\partial_ x p_{t-s}(z,x)\bigr\rrvert
\\
&&\quad \leq C_0 C_1(t)C_2(t)C_3(t)
\int_{[t-\epsilon,t]}ds\int_{\TT}\,dz \biggl(
\frac
{d(z,x)}{(t-s)^{3/2}}+\frac{1}{t-s} \biggr)\frac{\exp
(-d^2(z,x)/(2(t-s)))}{\sqrt{2\pi} }
\\
&&\quad = 2C_0C_1(t)C_2(t)C_3(t)
\int_{0}^\pi dz \int_0^\epsilon
ds \biggl(\frac{z}{s^{3/2}}+\frac{1}{s} \biggr)\frac{\exp
(-z^2/(2s))}{\sqrt{2\pi}}.
\end{eqnarray*}
This bound no longer depends on $x$ and to compute the latter integral,
consider the change of variable $u=z^2/s$, $z$ being fixed:
\begin{eqnarray*}
\int_{0}^\pi dz \int_0^\epsilon
ds \biggl(\frac{z}{s^{3/2}}+\frac{1}{s} \biggr){\exp
\bigl(-z^2/(2s)\bigr)} &=& \int_{0}^\pi
dz\int_{z^2/\epsilon}^{+\infty} du \biggl(\frac{1}{\sqrt{u}}+u
\biggr)\exp(-u/2).
\end{eqnarray*}
We conclude by remarking that by the dominated convergence theorem, the
latter term goes to zero with $\epsilon$.
\end{pf}

%
%reA.1 #&#
\begin{rem}
More generally, but still under the assumption that $\nu$ admits a
bounded density, the density $m_t$ is $\mathcal{C}^1$ at some time
$t>0$, if we can find
$\epsilon\in(0,t)$ such that
$\inf_{s\in[t-\epsilon,t]} \alpha_s>0$ and $\inf_{s\in[t-\epsilon
,t]} \beta_s>0$.
This comes
from the above proof or can be deduced directly from Lemma~\ref{mtC1}
and the Markov property of $X$.
\end{rem}

The same arguments cannot be used to prove that for $t>0$, the density
of $m_t$ belongs to $\mathcal{C}^2(\TT)$. A priori, this is annoying,
since in Section~\ref{proof},
to study the evolution of the quantity $I_t$ defined in~(\ref{It}), we
had to differentiate
it with respect to $t>0$ and the computations were justified only
if the densities $m_t$ were $\mathcal{C}^2$.
The classical way go around this apparent difficulty is to use a
mollifier.

Let $\rho$ be a smooth non-negative function on $\RR$ whose support
is included in $[-1,1]$
and satisfying $\int_\RR\rho(y) \,dy=1$. For any $\delta\in(0,1)$, define
\begin{eqnarray*}
\forall t\geq0, \forall x\in\TT,\qquad m^{(\delta)}_t(x)&\ffff &
\frac{1}\delta\int_\RR m_t(x+y)\rho
\biggl(\frac{y}{\delta} \biggr) \,dy
\end{eqnarray*}
(where functions on $\TT$ are naturally identified with $2\pi
$-periodic functions on $\RR$).
These functions are smooth and what is even more important
for Section~\ref{sec4},
the mapping $\RR_+^*\times\TT\ni(t,x)\mapsto\partial^2_x
m_t^{(\delta
)}(x)$ is continuous.
Furthermore, the $m^{(\delta)}_t$ are densities of probability
measures on
$\TT$.
More precisely, for any $t\geq0$, $m^{(\delta)}_t$ is the density of
$\cL
(X_t)$ when $\cL(X_0)=m^{(\delta)}_0$,
as a consequence of the linearity of the underlying evolution equation
(i.e., $\forall t\geq0, \partial_tm_t =m_tL_{\alpha_t,\beta
_t}$, in the
sense of distributions).
Thus, the computations of Section~\ref{sec4} are justified if we replace there
$(m_t)_{t>0}$ by
$(m^{(\delta)}_t)_{t>0}$, for any fixed $\delta\in(0,1)$.
In particular, the inequality (\ref{Fprime}) is satisfied for
$(m^{(\delta)}
_t)_{t>0}$ instead of
$(m_t)_{t>0}$. It remains to let $\delta$ go to $0_+$ to see that
the same bound is true for the flow $(m_t)_{t>0}$.
This proves Theorem~\ref{t1} for general initial distributions $m_0$,
for instance, Dirac masses.
In fact, one could pass to the limit $\delta\rightarrow0_+$ before
(\ref
{Fprime}), for instance,
already
in Proposition~\ref{Iprime}, to see that it is also valid.
\end{appendix}

\section*{Acknowledgements}
% \bigskip\bigskip
% \hskip5mm\textbf{\large Acknowledgments:} \sm\noindent
We would like to express our gratitude to our colleague S\'{e}bastien Gadat,
who has encoded and drawn all the figures of the present work.
We are also thankful to the referees whose suggestions
have helped us to improve the presentation of the paper.
Finally, we are grateful to
the support of
the Laboratoire de Math\'{e}matiques et Applications (UMR 7348)
of the Universit\'{e} de Poitiers, where much of this work took place.

%\section*{Acknowledgements}

%\begin{supplement}%[id=suppA]
%\slink[doi]{10.3150/00-BEJXXXXSUPP} %[doi,text={...}] - jei reikia
%suskaldyti doi
%\sfilename{BEJ000\_supp.pdf}
%\end{supplement}

% imsref loaded by linak, 2015-06-23 15:22:17
%

\printhistory
\end{document}